\documentclass[reqno,12pt]{amsart}

\usepackage{a4wide}
\usepackage[ps,all]{xy}
\usepackage[dvipsnames]{xcolor}
\usepackage[pdfusetitle]{hyperref}
\hypersetup{
 colorlinks=true,
 linkcolor=Blue,
 citecolor=Blue,
 urlcolor=Blue,
 }
 \usepackage[alphabetic,msc-links,initials]{amsrefs}
\usepackage{mathtools}
\usepackage{amsthm}
\usepackage{mathrsfs}
\usepackage{newtxtext}
\usepackage[varg]{newtxmath}
\usepackage{microtype}
\usepackage{tikz}
\usepackage{enumerate}
\usepackage{enumitem}
\setlist[enumerate,1]{label={\upshape(\arabic*)}}
\setlist[enumerate,2]{label={\upshape(\alph*)}}
\usetikzlibrary{arrows}
\usepackage{afterpage}
\theoremstyle{definition}
\newtheorem{theorem}{Theorem}
\newtheorem*{theorem*}{Theorem}
\newtheorem{definition}[theorem]{Definition}
\newtheorem*{definition*}{Definition}
\newtheorem{prop}[theorem]{Proposition}
\newtheorem*{prop*}{Proposition}
\newtheorem{lemma}[theorem]{Lemma}
\newtheorem*{lemma*}{Lemma}
\newtheorem{def-thm}[theorem]{Definition-Theorem}
\newtheorem*{def-thm*}{Definition-Theorem}
\newtheorem{def-pro}[theorem]{Definition-Proposition}
\newtheorem*{def-pro*}{Definition-Proposition}
\newtheorem{cor}[theorem]{Corollary}
\newtheorem*{cor*}{Corollary}
\newtheorem{example}[theorem]{Example}
\newtheorem*{example*}{Example}
\newtheorem{rem}[theorem]{Remark}
\newtheorem*{rem*}{Remark}

\newtheorem*{guide*}{Guide}

\numberwithin{theorem}{section}
\numberwithin{equation}{section}
\newtheorem*{fpa}{Framed preprojective algebras}
\newtheorem*{strata}{Parametrizations of nilpotent varieties}
\newtheorem*{mv}{Mirkovi\'c-Vilonen polytopes and generalized preprojective algebras}
\newtheorem*{org}{Organization}
\newtheorem*{conv}{Conventions and notation}
\newcommand{\textbm}[1]{\textit{#1}}
\newcommand{\sgn}[2]{\operatorname{sgn}\left(#1, #2\right)}
\newcommand{\eps}[1]{{\varepsilon}_{#1}}
\newcommand{\gcdd}[2]{\operatorname{gcd}(#1, #2)}
\newcommand{\diag}[1]{\operatorname{diag}\left(#1\right)}
\newcommand{\length}[1]{\mathit{l}( #1 )}
\newcommand{\summand}[1]{\left| #1 \right|}
\newcommand{\D}[1]{{\mathbb{D}} #1}
\newcommand{\Hom}[3]{\operatorname{Hom}_{#1}(#2, #3)}
\newcommand{\End}[2]{\operatorname{End}_{#1}( #2 )}

\newcommand{\Ext}[4]{\operatorname{Ext}^{#1}_{#2}(#3, #4)}
\newcommand{\Tor}[4]{\operatorname{Tor}_{#1}^{#2}\left(#3, #4\right)}
\newcommand{\rad}[1]{\operatorname{\mathrm{rad}}#1}

\newcommand{\fac}[1]{\operatorname{\mathrm{fac}}#1}

\newcommand{\sub}[1]{\operatorname{\mathrm{sub}}#1}
\newcommand{\add}[1]{\operatorname{\mathsf{add}}#1}
\newcommand{\Fac}[1]{\operatorname{\mathsf{Fac}}#1}
\newcommand{\Sub}[1]{\operatorname{\mathsf{Sub}}#1}
\newcommand{\Rep}[1]{\operatorname{\mathsf{Rep}}#1}
\newcommand{\rep}[1]{\operatorname{\mathsf{rep}}#1}
\newcommand{\tors}[1]{\operatorname{\mathsf{tors}}#1}
\newcommand{\torf}[1]{\operatorname{\mathsf{torf}}#1}

\newcommand{\Ker}[1]{\operatorname{\mathrm{Ker}}#1}
\newcommand{\Img}[1]{\operatorname{\mathrm{Im}}#1}
\newcommand{\Cok}[1]{\operatorname{\mathrm{Cok}}#1}
\newcommand{\fgpre}[2]{{\Pi}{\,}(#1, #2)}
\newcommand{\affPi}{\widetilde{\Pi}}
\newcommand{\hatPi}{\widehat{\Pi}}
\newcommand{\affI}[1]{\widehat{I}_{#1}}

\newcommand{\lfRep}[1]{\operatorname{\mathsf{Rep}_{\mathrm{l\mathchar`.f\mathchar`.}}}#1}
\newcommand{\lfrep}[1]{\operatorname{\mathsf{rep}_{\mathrm{l\mathchar`.f\mathchar`.}}}#1}
\newcommand{\proj}[1]{\operatorname{\mathsf{proj}}#1}

\newcommand{\sttilt}[1]{\mathsf{s\tau \mathchar`-tilt}{\,}#1}
\newcommand{\ind}{\operatorname{\mathsf{ind}}}

\newcommand{\pdim}[2]{\mathrm{proj.dim}_{#1}{\,}#2}

\newcommand{\rankvec}[1]{\operatorname{\underline{\mathrm{rank}}}#1}
\newcommand{\rk}{\operatorname{\mathrm{rank}}}
\newcommand{\Ltens}[1]{\otimes^{\mathbb{L}}_{#1}}
\newcommand{\dimvc}[1]{\operatorname{\underline{\mathrm{dim}}}#1}
\newcommand{\op}[1]{{#1}^\mathrm{op}}
\newcommand{\tens}[1]{{\otimes}_{#1}}
\newcommand{\Efilt}[1]{\operatorname{\mathbb{E}-\mathsf{filt}}(#1)}
\newcommand{\Irr}[1]{\operatorname{\sf{Irr}}\left( #1 \right)}
\newcommand{\Irrmax}[1]{\operatorname{\sf{Irr}}\left( #1 \right)^{\mathrm{max}}}
\newcommand{\djuni}[1]{\coprod_{\mathbf{r}\in {\mathbb{N}}^n}#1}
\newcommand{\br}{\mathbf{r}}
\newcommand{\wt}[1]{\operatorname{wt}( #1 )}
\newcommand{\Pol}[1]{\operatorname{Pol}\left( #1 \right)}
\newcommand{\T}{\mathfrak{T}}
\newcommand{\F}{\mathfrak{F}}
\newcommand{\calT}{\mathcal{T}}
\newcommand{\calF}{\mathcal{F}}
\newcommand{\calB}{\mathcal{B}}
\newcommand{\GL}{\mathrm{GL}}
\newcommand{\M}{\mathbf{M}}
\newcommand{\bN}{\mathbf{N}}
\newcommand{\N}{\mathbb{N}}
\newcommand{\B}{\mathfrak{B}}
\newcommand{\id}{\mathrm{id}}
\newcommand{\R}{\mathbb{R}}
\newcommand{\C}{\mathbb{C}}
\newcommand{\Z}{\mathbb{Z}}
\newcommand{\E}{\mathbb{E}}

\newcommand{\trp}{\mathsf{T}}
\newcommand{\nil}{\mathsf{nil}}
\newcommand{\ad}{\mathrm{ad}}
\newcommand{\Aut}{\mathrm{Aut}}
\newcommand{\spn}{\mathrm{Span}}
\newcommand{\Tmax}[1]{T^{\mathrm{max}}_{#1}}
\newcommand{\Tmin}[1]{T^{\mathrm{min}}_{#1}}
\newcommand{\MV}{\mathcal{MV}}
\newcommand{\ie}{\textit{i.e}.\,}
\newcommand{\resp}{\textit{resp}.\,}
\newcommand{\eg}{\textit{e.g.}~}
\newcommand{\cf}{\textit{cf.}~}
\newcommand{\etc}{\textit{etc.}}
\newcommand{\Gr}{\mathrm{Gr}}

\newcommand{\h}{\mathfrak{h}}
\newcommand{\g}{\mathfrak{g}}

\author[K. Murakami]{Kota Murakami}
\address{Department of Mathematics, Kyoto University, Kitashirakawa Oiwake-cho, Sakyo-ku, 
Kyoto 606-8502, Japan}
\email{k-murakami@math.kyoto-u.ac.jp}
\title[PBW parametrizations and generalized preprojective algebras]{PBW parametrizations and generalized preprojective algebras}
\thanks{This work was supported by the Kyoto Top Global University project,  Grant-in-Aid for JSPS Fellows (JSPS KAKENHI Grant Number JP21J14653) and JSPS bilateral program (Grant
Number JPJSBP120213210). This work was also supported by the Research Institute for Mathematical Sciences,
an International Joint Usage/Research Center located in Kyoto University.}
\begin{document}
\maketitle
\begin{abstract}
    Gei\ss-Leclerc-Schr\"oer [Invent.~Math.~{\bf{209}}~(2017)] has introduced a notion of generalized preprojective algebras associated with generalized Cartan matrices and their symmetrizers. These algebras realize crystal structures on the set of maximal dimensional irreducible components of the nilpotent varieties [Selecta~Math.~(N.S.)~{\bf{24}}~(2018)]. For general finite types, we give stratifications of these components via partial orders of torsion classes in module categories of generalized preprojective algebras in terms of Weyl groups. In addition, we realize Mirkovi\'c-Vilonen polytopes from generic modules of these components, and give an identification as crystals between the set of Mirkovi\'c-Vilonen polytopes and the set of maximal dimensional irreducible components. This generalizes results of Baumann-Kamnitzer [Represent.~Theory~{\bf{16}}~(2012)] and Baumann-Kamnitzer-Tingley [Publ.~Math.~Inst.~Hautes~{\'E}tudes~Sci.~{\bf{120}}~(2014)].
\end{abstract}
\tableofcontents
\section{Introduction}
Preprojective algebras are introduced by Gelfand-Ponomarev~\cite{GP} in order to understand representation categories of finite quivers. A preprojective algebra is defined as a double of a finite quiver with preprojective relations, which contains the path algebra of the original quiver as a subalgebra. These algebras have many applications to different contexts, for example tilting theory and McKay correspondences~\cite{AmIR,BIRS}, cluster structures of coordinate rings of algebraic groups or flag varieties~\cite{GLS7,Lec}, and geometric representation theory~\cite{Nak2,Nak1,VS}, \etc

Recently, Gei\ss-Leclerc-Schr\"oer~\cite{GLS1} introduced a certain class of $1$-Iwanaga-Gorenstein algebras and a notion of \emph{generalized preprojective algebras} for symmetrizable generalized Cartan matrices (=GCMs) and their symmetrizers (see Definition \ref{def of H&Pi}). These classes of algebras contain path algebras of acyclic quivers and their preprojective algebras as special cases. They generalize a part of links between Kac-Moody Lie algebras and preprojective algebras (or representations of quivers) from simply-laced types to symmetrizable types which contain $\mathsf{B, C, F, G}$ types (\textit{e.g.}~\cite{GLS3,GLS4}).
On the other hand, in studies of quantum enveloping algebras and algebraic groups, it is important to find bases of an algebra or its irreducible representations which have nice properties. Many bases are defined and are compared in this context, for example (dual) canonical bases~\cite{Lus4}, (dual) semi-canonical bases~\cite{Lus3}, Mirkovi\'c-Vilonen bases~\cite{MV} and theta bases~\cite{GHKK}. The canonical bases are described by some kinds of geometry, that is, affine Grassmannians and Nakajima quiver varieties. Kamnitzer \textit{et al.}~\cite{BKK,KTWWY}, \etc have compared these two kinds of geometry and have developed their relationships from a viewpoint of representation theory about a duality of conical symplectic singularities. 

Motivated from their works, we extend relationships between some combinatorics of affine Grassmannians and representation theory of quivers from a symmetric setting to our symmetrizable setting. In particular, we employ the representation theory of generalized preprojective algebras instead of employing Nakajima quiver varieties that are not available for $\mathsf{B, C, F, G}$ cases. Namely, we compare Kamnitzer's theory of Mirkovi\'c-Vilonen polytopes~\cite{Kam1} and representation theory of generalized preprojective algebras for symmetrizable Cartan matrices from a viewpoint of tilting theory.
\begin{fpa}
In representation theory of quivers, \emph{framing} is a useful technique for understanding some symmetries associated with quivers in terms of quantum groups and cluster algebras~\cite{Kel,Rei}. A framed quiver has original vertices $i$ and added corresponding vertices $i'$. Preprojective algebras of framed quivers have similar combinatorial features to Nakajima quiver varieties in the module categories, called “Nakajima's tricks” in the work of Baumann-Kamnitzer~\cite{BK}. Roughly speaking, modules over these algebras describe data of positive roots by dimensions of vector spaces on original vertices, and data of integral weights by those on added vertices. As a by-product of this construction, we obtain some constructible functions on varieties of nilpotent modules of preprojective algebras through a class of modules, called stable modules. Such modules are known to have a few kinds of constructions other than a use of framed preprojective algebras ({\eg}\cite{BIRS,GLS7}), and these types of constructible functions are often useful for relating preprojective algebras with contexts of canonical bases or cluster algebras. We generalize the definition and develop some algebraic properties of framed preprojective algebras and their module categories from symmetric cases to symmetrizable cases. Namely, we define framed preprojective algebras (see \S \ref{subsec-stab}) for generalized preprojective algebras. In particular, our algebras contain classical ones in \cite{BK} as special cases. By adopting the technique in the work of Baumann-Kamnitzer~\cite{BK} for classical cases, we prove that there are nice “generalized stable modules” as a generalization of \cite[Theorem 3.1]{BK}. The difference between the original work and our approach is that we must consider symmetrizers and work with suitable subcategories instead of the whole module categories.
\begin{theorem}[$\doteq$Theorem \ref{thm_stabmod}; see Definitions \ref{def of H_i&lfrep}, \ref{def_reflection} and \S \ref{pre:GLSalg},\S \ref{subsec:refl} for notation]\label{mainthm:1}
Let $\Pi\coloneqq\fgpre{C}{D}$ be the generalized preprojective algebra associated with a GCM $C$ and a symmetrizer $D$ of $C$. Let $\varpi_i$ be the fundamental weight with respect to $i\in Q_0$ and let $w\in W$, where $W$ is the Weyl group of $C$.
Then, there is a unique $\Pi$-module $N(\gamma)$ for each weight $\gamma$ in Tits cone, which satisfies the following:
	\begin{enumerate}
	\item If $\gamma\in \mathbb{Z}^n$ is an anti-dominant weight, then $N(\gamma)=0$;
	\item If $\varpi_i -w\varpi_i\neq 0$, then there exists unique $N(-w\varpi_i)\in \lfrep{\Pi}$ such that $\rankvec N(-w\varpi_i)=\varpi_i -w\varpi_i$ and $\sub N(-w\varpi_i)=E_i$. Conversely, if a module $M\in \lfrep{\Pi}$ satisfies $\rankvec M=\varpi_i -w\varpi_i$ and $\sub M=E_i$, then $M\cong N(-w\varpi_i)$;
	\item If $\gamma$ and $\delta$ belong to the same Weyl chamber, then $N(\gamma+ \delta)\cong N(\gamma)\oplus N(\delta)$.
	\end{enumerate}
\end{theorem}
This result makes us expect that there is nice geometry like Nakajima quiver varieties for Gei\ss-Leclerc-Schr\"oer's theory, via the existence of some nice combinatorial features of framed preprojective algebras. In fact, they enable us to develop a relationship between generalized preprojective algebras and canonical bases in Theorem~\ref{mainthm:2} and Theorem~\ref{mainthm:3}.
\end{fpa}
\begin{strata}
Kashiwara-Saito~\cite{KS} developed combinatorics of irreducible components and constructible functions on varieties of nilpotent modules over preprojective algebras, called \emph{nilpotent varieties}, in terms of Kashiwara's crystal basis~\cite{Kas1}. After their work, Lusztig~\cite{Lus3} constructed semi-canonical bases via generic values of constructible functions on irreducible components of nilpotent varieties. Gei\ss-Leclerc-Schr\"oer~\cite{GLS8} compared the multiplication of dual semi-canonical bases with that of dual canonical bases by using the decomposition theory of irreducible components in module varieties of Crawley-Boevey-Schr\"oer~\cite{CBS}.

Recently, Gei\ss-Leclerc-Schr\"oer~\cite{GLS4} defined nilpotent varieties for generalized preprojective algebras. They gave a crystal structure on the set of special irreducible components of nilpotent varieties, called \emph{maximal components}, and gave a conjectural description of dual semi-canonical bases as a generalization of the classical theory of Kashiwara-Saito~\cite{KS} and Lusztig~\cite{Lus3}. We give a parametrization of irreducible components of Gei\ss-Leclerc-Schr\"oer's varieties by studying a special class of modules, called \emph{crystal modules} by the $\tau$-tilting theory of Adachi-Iyama-Reiten~\cite{AIR}. This kind of results for symmetric GCMs can be found in the works of Geiss-Leclerc-Schr\"oer~\cite[Proposition 14.6]{GLS7} and Baumann-Kamnitzer-Tingley~\cite[paragraph after Proposition 5.23]{BKT}.
\begin{theorem}[$\doteq$Theorem~\ref{strata-thm}]\label{mainthm:2}
Let $\Pi\coloneqq \fgpre{C}{D}$  be the generalized preprojective algebra associated with a Cartan matrix $C$ and its symmetrizer $D$, and let $\mathbf{i}\coloneqq (i_1, \dots , i_r)$ be any reduced expression of an element $w$ of the Weyl group $W(C)$. We define
$V_{\mathbf{i}, k} := N\left (-s_{i_{1}}\cdots s_{i_{k}}
\varpi _{i_{k}}\right )$ and
$M_{\mathbf{i}, k}:= V_{\mathbf{i}, k}/V_{\mathbf{i}, k^{-}}$ for
$1\leq k \leq r$ where
$k^{-} =\max \{0, 1\leq s\leq k-1\mid i_{s}= i_{k}\}$. Then, there are constructible sets $\Pi_{\mathbf{i}}^{\mathbf{a}}$ of crystal modules $M$ which have a filtration $$M= T_0 \supsetneq T_1 \supsetneq \cdots \supsetneq T_r =0$$
such that 
$T_j/ T_{j+1} \cong (M_{\mathbf{i}, k}^{\oplus a_j})^*$ for $\mathbf{a} \coloneqq (a_1, \dots, a_r) \in \Z^r_{\geq 0}$, where $*$ is the natural involution $\operatorname{\mathsf{rep}}\Pi \rightarrow \operatorname{\mathsf{rep}}\Pi$ induced by taking the $K$-dual space. The constructible set $\Pi_{\mathbf{i}}^{\mathbf{a}}$ is irreducible. Moreover, $Z_{\mathbf{i}}^{\mathbf{a}}\coloneqq\overline{\Pi_{\mathbf{i}}^{\mathbf{a}}}$ is a maximal irreducible component of the nilpotent variety.
\end{theorem}
Note that any maximal irreducible component is given in this way because any maximal irreducible component has a crystal module as a generic point. In particular, each crystal module has a filtration in above Theorem \ref{mainthm:2} for arbitrary reduced expression of $w_0$. This theorem gives a geometric interpretation of the PBW type parametrizations of Kashiwara crystal bases by Saito~\cite{Sai}.
\end{strata}
\begin{mv}
Let $G$ be a complex reductive group and let $\Gr_{G}=G(\mathbb{C}(\!(t)\!))/G(\mathbb{C}[\![t]\!])$ be its affine Grassmannian. $\Gr_G$ has subvarieties $\Gr^{\lambda}_G$ whose intersection homologies are isomorphic to the irreducible highest weight representations $V(\lambda)$ of the highest weight $\lambda$ of the Langlands dual group $G^{\vee}$. Mirkovi\'c-Vilonen~\cite{MV} studied algebraic cycles of affine Grassmannian which form a basis of the intersection homology of $\Gr^{\lambda}_G$ mentioned above, called the Mirkovi\'c-Vilonen cycles (=MV cycles).

Kamnitzer~\cite{Kam1} defined the \emph{Mirkovi\'c-Vilonen polytopes} (=MV polytopes) in order to extract combinatorial features of MV cycles. The set of MV polytopes is associated with the root datum of $G$, and they have much information of the quantum enveloping algebra associated with the Lie algebra of the Langlands dual group $G^{\vee}$ as a result of the geometric Satake correspondence. In particular, the set of MV polytopes for $G$ bijectively corresponds to Kashiwara's crystal bases, and they have many numerical data about crystal bases which are developed by Lusztig~\cite{Lus4} and Berenstein-Zelevinsky~\cite{BZ}, including $\mathbf{i}$-Lusztig data via their edge lengths.

Baumann-Kamnitzer-Tingley~\cite{BK,BKT} compared combinatorial structures of module categories over preprojective algebras and MV polytopes along with the context of $\mathbf{i}$-Lusztig data. In particular, they realized MV polytopes in the Grothendieck group of the module category of a preprojective algebra by taking convex hull of the set of submodules of generic modules of irreducible components of the nilpotent variety. We generalize their arguments to our symmetrizable setting. This extends the results of Baumann-Kamnitzer-Tingley~\cite{BK,BKT} to our setting. Our approach is based on the $\tau$-tilting theory unlike their arguments.
\begin{theorem}[$\doteq$Theorem \ref{nilp_MV}, Theorem \ref{bij_MV&B}]\label{mainthm:3}
Let $\Pi\coloneqq \fgpre{C}{D}$  be the generalized preprojective algebra associated with a Cartan matrix $C$ and its symmetrizer $D$ except for type $\mathsf{G}_2$. Let $\Gamma\coloneqq \left\{ w\varpi_i \mid w\in W, i\in Q_0 \right\}$ be the set of chamber weights.
We consider each maximal irreducible component $Z$ of a variety $\Pi(\mathbf{r})$ of suitable modules over $\Pi$ and a generic module $T$ of this component. Then, there is a function $D_{\gamma}$ on $\Pi(\br)$, which defines a constructible function on each component $Z$. In addition, $P(T)=\{v\in \mathbb{R}^n \mid \langle \gamma, v \rangle \leq D_{\gamma}(T), \gamma \in \Gamma\}$ is an MV polytope associated with the Langlands dual root datum of $C$. In particular, the map $P{(-)}\colon \mathcal{B}\rightarrow \MV$ from the set of maximal irreducible components of a nilpotent variety coming from $C$ to the set of MV polytopes coming from the root datum of $C^\trp$ is an isomorphism.
\end{theorem}
Indeed, the function $D_{\gamma}(-)$ is given by $\dim \Hom{\Pi}{N(\gamma)}{-}$ (up to scalar) from Theorem \ref{mainthm:1}, and Theorem \ref{mainthm:2} gives the vertex description of our polytopes. Recently, Baumann-Kamnitzer-Knutson \textit{et al}.~\cite{BKK,HKW} try to relate the geometry of quiver Grassmannians of modules over preprojective algebras and the geometry of affine Grassmannian slices. In their observation, quiver Grassmannians from generic modules and MV cycles from these modules are related via representation theory of current Lie algebras. Our result suggests existence of some more geometric relationship between the theory of Gei\ss-Leclerc-Schr\"oer~\cite{GLS1,GLS2,GLS3,GLS4,GLS5,GLS6} and the geometry of affine Grassmannians (\eg the mathematical definition of Coulomb branches of quiver gauge theories with symmetrizers by Nakajima-Weekes~\cite{NW}).
\end{mv}
\begin{org}
In \S \ref{sec:pre}, we prepare our conventions about root systems, representation theory of quivers, tilting theory and reflection functors. In \S \ref{section-stab}, we develop two constructions about stability conditions, namely stable modules and $g$-vectors as generalizations of results of Baumann-Kamnitzer~\cite{BK} and Mizuno~\cite{Miz}. In \S \ref{section-MV}, we compare the data about stability conditions of module categories of generalized preprojective algebras and crystal basis. In particular, we realize Mirkovi\'c-Vilonen polytopes for general finite types as generalizations of results of Baumann-Kamnitzer~\cite{BK} and Baumann-Kamnitzer-Tingley~\cite{BKT}.
\end{org}
\begin{conv}
Throughout this paper, we refer to a field $K$ as a commutative algebraically closed field unless specified otherwise. However, the results obtained before \S~\ref{section-MV} are valid for any field. If we refer to a $K$-algebra $\Lambda$, then $\Lambda$ means a unital associative non-commutative $K$-algebra. If we refer to a module over a ring $\Lambda$, the module means a left $\Lambda$-module. The category $\Rep{\Lambda}$ (\resp $\rep{\Lambda}$) denotes the module category (\resp the category of finite dimensional modules) over a $K$-algebra $\Lambda$. For a finite dimensional $K$-algebra $\Lambda$, the duality $\D \coloneqq \Hom{K}{-}{K}\colon \rep{\Lambda}\rightarrow \rep{\op{\Lambda}}$ denotes the standard $K$-duality. For a finite dimensional $K$-algebra $\Lambda$, the covariant functor $\nu \coloneqq \D(\Lambda)\tens{\Lambda}(-)\colon \rep{\Lambda}\rightarrow \rep{\Lambda}$ denotes the Nakayama functor. For a module $M$ over a $K$-algebra $\Lambda$, the category $M^\bot$ denotes the full subcategory $\{N\in \rep{\Lambda} \mid \Hom{\Lambda}{M}{N}=0\}$ of $\rep{\Lambda}$. The category $\proj{\Lambda}$ denotes the full subcategory of finitely generated projective modules over a ring $\Lambda$. The commutative group $K_0(\mathcal{E})$ denotes the Grothendieck group of an exact category $\mathcal{E}$. For a commutative ring $R$, the set $M_n(R)$ denotes the set of $n\times n$-square matrices. For a matrix $A$, the matrix $A^\trp$ denotes the transpose matrix of $A$. We always consider a variety as a reduced separated scheme of finite type over an algebraically closed field $K$, whose topology is the Zariski topology. For an algebraic variety $V$, the set $\Irr{V}$ denotes the set of irreducible components of $V$.
\end{conv}

\section{Preliminary}\label{sec:pre}
\subsection{Cartan matrices, root systems and Weyl groups}\label{pre:rootsys}
We recall basic concepts about root systems. We refer to Kac~\cite{Kac}, Bj\"orner-Brenti~\cite{BjBr} and Humphreys~\cite{Hum}.

\begin{definition}\label{GCM}
	A matrix $C=(c_{ij}) \in M_n(\mathbb{Z})$ is called a \textbm{generalized Cartan matrix} (GCM), if it satisfies the following three conditions:
	\begin{description}
		\item[\sf{(C1)}] For each $1\leq i\leq n$, we have $c_{ii} =2$;
		\item[\sf{(C2)}] If $i \neq j$, then $c_{ij} \leq 0$;
		\item[\sf{(C3)}] We have $c_{ij} \neq 0$ if and only if we have $c_{ji} \neq 0$.
	\end{description}	
	In particular, if $C=(c_{ij})$ satisfies the following \textsf{(C4)}, then $C$ is called \textbm{symmetrizable}.
	\begin{description}
		\item[\sf{(C4)}] There is a diagonal matrix $D=\diag{c_1, \dots , c_n}{\,}(c_i \in \mathbb{Z}, c_i \geq 1)$ such that $DC$ is a symmetric matrix.
	\end{description}	
In the condition (C4), the matrix $D$ is called a \textbm{symmetrizer} of $C$.
\end{definition}
The following quadratic form $q_{DC}$ and graph $\Gamma(C)$ give a classification of GCMs.
\begin{definition}
Let $C=(c_{ij}) \in M_n(\mathbb{Z})$ be a symmetrizable GCM and let $D=\diag{c_1, \dots , c_n}$ be a symmetrizer of $C$.
\begin{enumerate}
	\item The graph $\Gamma(C)$ has vertices $1, \dots , n$. An edge between $i$ and $j$ exists in $\Gamma(C)$ if and only if $c_{ij}<0$. The edges are labeled by pairs $(|c_{ji}|, |c_{ij}|)$ as: 
$$\begin{tikzpicture}[auto]
	\node (a) at (0, 0) {$i$}; \node (b) at (2.4, 0) {$j$}; 
	\draw[-] (a) to node {$\scriptstyle (|c_{ji}|, |c_{ij}|)$} (b);
\end{tikzpicture}.$$
We call $\Gamma(C)$ the \textbm{valued graph} of $C$. If $\Gamma(C)$ is a connected graph, then $C$ is called \textbm{connected}.
\item We define the form $q_{DC} :\mathbb{Z}^n\longrightarrow \mathbb{Z}$ by
      $$q_{DC} =\sum_{i=1}^n c_i X_i^2 +\sum_{i<j}c_i c_{ij}X_iX_j.$$
     Since we have $c_i c_{ij} =c_j c_{ji}$, the form $q_{DC}$ is symmetric. If $q_{DC}$ is positive definite (\textit{resp}. positive semi-definite), then $C$ is called \textit{finite type} (\textit{resp}. \textit{Euclidean type}).
\end{enumerate}
\end{definition}
\begin{rem}
	If $C$ is a connected symmetrizable GCM, there exists a unique \textit{minimal} symmetrizer. That is, any symmetrizer of $C$ is equal to $mD$ for the minimal symmetrizer $D$ and some positive integer $m$.
\end{rem}
Next, we define data of roots and weights of complex semisimple Lie algebras. We put $Q_0 \coloneqq \{1, 2, \dots, n\}$. Let $\h$ be an $n$-dimensional $\C$-vector space, and let $\{\alpha_1, \dots, \alpha_n \} \subset \h^*$ and $\{ \alpha^{\vee}_1, \dots, \alpha^{\vee}_n \}\subset\h$ be the set of simple roots (\resp the set of simple co-roots).
We define bilinear forms $\langle -, -\rangle_C\colon \h^* \times \h^*\rightarrow \C$ and $(-, -)\colon \h^* \times \h^*\rightarrow \C$ by 
\begin{align*}
    &\langle \alpha_i, \alpha_j \rangle_{C} \coloneqq \alpha_i ( \alpha_j^{\vee}) =c_{ji},
     &&  ( \alpha_i, \alpha_j )_C \coloneqq c_jc_{ji}=c_i c_{ij}
\end{align*}
for any $i, j\in Q_0$.
We fix a basis $\{\varpi_j \mid 1\leq j \leq n\}$ of $\h^*$ such that $\varpi_j( \alpha_i^{\vee}) = \delta_{ij}$ for $1\leq i, j\leq n$. We call them fundamental weights. We have $\alpha_i = \sum_{j\in Q_0} c_{ji}\varpi_j$.

We define the \textbm{integral weight lattice} and the set of \textbm{dominant integral weights} by
\begin{align*}
    &P\coloneqq \left\{\nu \in \h^* \mid \left\langle \nu, \alpha_i \right\rangle \in \Z, i\in Q_0 \right\}= \bigoplus_{i\in Q_0}\Z \varpi_i;\\
    &P^+ \coloneqq \left\{\nu \in P \mid \left\langle \nu, \alpha_i \right\rangle \geq 0, i\in Q_0 \right\}
    =\bigoplus_{i \in Q_0} \Z_{\geq 0} \varpi_i.
    \end{align*}
    By the above definition, we can identify $\varpi_i$ with the $i$-th standard basis vector of $\Z^n$. Furthermore, we define the \textbm{root lattice} and the \textbm{positive root lattice} by
    \begin{align*}
        &R\coloneqq \bigoplus_{i\in Q} \Z \alpha_i,
        &&R^+\coloneqq \bigoplus_{i\in Q} \Z_{\geq 0} \alpha_i.
    \end{align*}
    By the above definition, we can also identify $\alpha_i$ with the $i$-th standard basis vector of $\Z^n$. Finally, we define the Weyl group and the root system. Define $s_i\in \Aut{(\Z^n)}$ by $s_i(\alpha_j ) \coloneqq \alpha_j - c_{ij} \alpha_i$ for each $i, j \in Q_0$.
    $W(C)$ denotes the group generated by $s_1, \dots, s_n$ and we call $W(C)$ the Weyl group of $C$.
     In general, the Weyl group $W(C)$ associated with the Kac-Moody Lie algebra $\g(C)$ is generated by $s_i\,\,(i\in Q_0)$ and only relations
    \begin{align*}
        s_i^2&=e\,\,(i\in Q_0);\\
        s_i s_j&=s_j s_i\,\,(c_{ij}c_{ji}=0);\\
        s_i s_j s_i&= s_j s_i s_j\,\,(c_{ij}c_{ji}=1);\\
        s_i s_j s_i s_j &= s_j s_i s_j s_i\,\,(c_{ij}c_{ji}=2);\\
        s_i s_j s_i s_j s_i s_j &= s_j s_i s_j s_i s_j s_i\,\,(c_{ij}c_{ji}=3).
    \end{align*}
    We define the set of (real) roots by $$\Delta\coloneqq\left\{w(\alpha_i)\mid w\in W(C), i\in Q_0 \right\}.$$
    
    Finally, we recall Weyl chambers.
    Let $\varpi_1, \dots , \varpi_n \in \mathbb{Z}^n$ be the standard basis of $\mathbb{Z}^n$. We define a group monomorphism $\sigma^* \colon W(C) \rightarrow \mathrm{GL}(\mathbb{Z}^n)$ by
    $$\sigma^*(s_i)(\varpi_j) \coloneqq \begin{cases} \varpi_i - \sum_{k\in Q_0}c_{kj}\varpi_k \,\,&(i=j)\\
    \varpi_j \,\,&(i \neq j).
    \end{cases}$$
    For $a, b \in \R^n$, we write $\langle a, b \rangle\coloneqq a\cdot b^\trp$. For each root $x \in \Delta$, we define a hyperplane $H_x \coloneqq \left\{y \in \R^n \mid \langle y, x \rangle =0\right\}$. 
    We refer to a connected component of $\R^n \backslash \bigcup_{x\in \Delta} H_x $ as a \textbm{Weyl chamber}. $C_0$ denotes the chamber
    $\{ a_1 \varpi_1 + \cdots + a_n \varpi_n \mid a_i \in \mathbb{R}_{>0} \}$. A Weyl chamber is given by $\sigma^* (w)(C_0)$ for some $w\in W(C)$ and there is a bijection between $W$ and the set of Weyl chambers by $w \mapsto \sigma^*(w)(C_0)$. We write $\Gamma\coloneqq \left\{ w\varpi_i \mid w\in W, i\in Q_0 \right\}$.
    
    For $w\in W(C)$, we take an expression $w=s_{i_1}\cdots s_{i_k}$. If $k$ is minimal among such expressions, then we write $\length{w}=k$ and refer to this expression as a reduced expression of $w$.
    If there are no elements $w'\in W(C)$ such that $\length{w'}>\length{w}$ for $w\in W(C)$, then we call $w$ the longest element of $W$ and $w_0$ denotes the longest element of $W(C)$.
    Finally, we define the right weak Bruhat order $\leq_R$ on $W(C)$. We say that $u\leq_{R} v$ if there exist simple reflections $s_{i_1}, \dots s_{i_k}$ such that $v= us_{i_1}\cdots s_{i_k}$ and $\length{us_{i_1}\cdots s_{i_j}}=\length{u} +j$ for $1\leq j \leq k$.
    
    We write with hats {\,\,}$\widehat{}${\,\,} concepts associated with an un-twisted affine type GCM $\widehat{C}$  for a finite type GCM $C$. (\textit{e.g}. affine Weyl groups $\widehat{W}(C)\coloneqq W(\widehat{C})$).
    
\subsection{Preprojective algebras associated with symmetrizable Cartan matrices}\label{pre:GLSalg}
In this subsection, we briefly review the representation theory of preprojective algebras with symmetrizable Cartan matrices introduced in  Gei\ss-Leclerc-Schr\"oer~\cite{GLS1}. In their work, the representation theory of acyclic quivers and of their preprojective algebras are generalized by using the data of symmetrizable GCMs and their symmetrizers. First, we define a quiver from a symmetrizable GCM.
\begin{definition}\label{parameta-orientation}
Let $C=(c_{ij}) \in M_n(\mathbb{Z})$ be a symmetrizable GCM and let $D=\diag{c_1, \dots , c_n}$ be its symmetrizer.
\begin{enumerate}
\item If $c_{ij} <0$, then we set
		$g_{ij}\coloneqq |\gcdd{c_{ij}}{c_{ji}}|$ and $f_{ij}\coloneqq |c_{ij}| /g_{ij}{\,}$.
\item If $\Omega \subset \{1, 2, \dots , n\} \times \{1, 2, \dots , n\}$ satisfies the following two conditions, then $\Omega$ is called an \textit{orientation} of $C$.
\begin{itemize}
	\item If $\{(i, j), (j, i)\} \cap \Omega \neq\phi$ holds, then $c_{ij} < 0$;
	\item Every sequence $((i_1, i_2), (i_2, i_3), \dots , (i_t, i_{t+1})){\,\,}(t \geq 1)$ in $\{1, 2, \dots , n\} \times \{1, 2, \dots , n\}$ such that each $(i_s, i_{s+1}) \in \Omega{\,\,\,}(1 \leq s \leq t)$ satisfies $i_1 \neq i_{t+1}$.
\end{itemize}
    We define the opposite orientation ${\Omega}^* \coloneqq\{(j, i){\,}\mid(i, j) \in \Omega \}$ and $\overline{\Omega}\coloneqq\Omega \cup {\Omega}^*$.
\end{enumerate}
\end{definition}
\begin{definition} \label{C-quiver}
	Under the setting of Definitions \ref{GCM}, \ref{parameta-orientation}, we define:
	\begin{enumerate}
	\item The quiver $Q= Q{\,}(C, \Omega)= \left(Q_0, Q_1, s, t\right)$ as follows:
	\begin{align*}	
		&Q_0= \{1, 2, \dots , n\},\\
		&Q_1= \{{\alpha}^{(g)}_{ij} :j \rightarrow i{\,}\mid(i, j) \in \Omega, 1 \leq g \leq g_{ij} \} \cup \{ {\varepsilon}_i :i \rightarrow i \mid 1 \leq i \leq n \},\\
		&s({\alpha}^{(g)}_{ij}) =j,{\,} t({\alpha}^{(g)}_{ij})=i,{\,} s({\varepsilon}_i) =t({\varepsilon}_i) =i\text{.}
	\end{align*}
	\item The \textit{double quiver} $\overline{Q} =(\overline{Q}_0, \overline{Q}_1, s, t)$ of $Q$ as follows:
	\begin{align*}	
		&\overline{Q}_0= Q_0= \{1, 2, \dots , n\},\\
		&\overline{Q}_1= \{{\alpha}^{(g)}_{ij} :j \rightarrow i{\,}\mid(i, j) \in \overline{\Omega}, 1 \leq g \leq g_{ij} \} 
		\cup \{ {\varepsilon}_i :i \rightarrow i \mid 1 \leq i \leq n \},\\
		&s({\alpha}^{(g)}_{ij}) =j,{\,} t({\alpha}^{(g)}_{ij})=i,{\,} s({\varepsilon}_i) =t({\varepsilon}_i) =i\text{.}
	\end{align*}
\end{enumerate}	
\end{definition}
Finally, we define $K$-algebras $H$ and $\Pi$ as quivers with relations.
\begin{definition} \label{def of H&Pi}
Under the setting of Definition \ref{C-quiver}, we define:
\begin{enumerate}
\item We define a $K$-algebra $H= H(C, D, \Omega)\coloneqq KQ/I$ by the quiver $Q$ with relations $I$ generated by \textsf{(H1), (H2)}:
		\begin{enumerate}
			\item[\textsf{(H1)}] $\eps{i}^{c_i}=0$\,\,\,($i\in Q_0$);
			\item[\textsf{(H2)}] For each $(i, j) \in \Omega$, we have $\eps{i}^{f_{ji}}\alpha^{(g)}_{ij}=\alpha^{(g)}_{ij}\eps{j}^{f_{ij}}$\, $(1 \leq g \leq g_{ij})$.
		\end{enumerate}
	\item	We define a $K$-algebra $\Pi= \fgpre{C}{D, \overline{\Omega}} \coloneqq K\overline{Q}/\overline{I}$ by the quiver $\overline{Q}$ with relations $\overline{I}$ generated by \textsf{(P1)-(P3)}:
		\begin{enumerate}
			\item[\textsf{(P1)}] $\eps{i}^{c_i}=0$ ($i\in Q_0$);
			\item[\textsf{(P2)}] For each $(i, j) \in \overline{\Omega}$, we have $\eps{i}^{f_{ji}}\alpha^{(g)}_{ij}=\alpha^{(g)}_{ij}\eps{j}^{f_{ij}}$ ($1 \leq g \leq g_{ij}$);
			\item[\textsf{(P3)}] For each $i \in Q_0$, we have
			$$\sum_{j \in \overline{\Omega}{\,}(-, i)} \sum_{g=1}^{g_{ij}} \sum_{f=0}^{{f_{ji}}-1} \sgn{i}{j}\eps{i}^f \alpha^{(g)}_{ij} \alpha^{(g)}_{ji} \eps{i}^{f_{ji}-1-f}=0, $$
			where we define
			\begin{align*}
	\overline{\Omega}{\,}(i, -)\coloneqq \{j\in Q_0\mid (i, j)\in \overline{\Omega}\}, &&\overline{\Omega}{\,}(-, j)\coloneqq\{i\in Q_0\mid (i, j)\in \overline{\Omega}\},
\end{align*}
and 
			$$\sgn{i}{j}\coloneqq 
			 \begin{cases}
				1 &(i, j)\in \Omega,\\
				-1 &(i, j)\in \Omega^*.
			\end{cases}$$
		\end{enumerate}
		\end{enumerate}
	 In the definition of $\Pi$, the set $\{ e_i\mid i\in Q_0\}$ is the complete set of primitive orthogonal idempotent elements corresponding to the set $Q_0$ of vertices. We note that $\Pi$ does not depend on the choice of orientation $\Omega$ up to isomorphism. We denote the above $\fgpre{C}{D, \overline{\Omega}}$ by $\fgpre{C}{D}$ in short. We refer to the above $\Pi=\fgpre{C}{D}$ as the \emph{generalized preprojective algebra} associated with the pair $(C, D)$.
\end{definition}
\begin{example} \label{exampreproj}
	Let $C=\left(
	\begin{array}{rr}
		2 & -1 \\
	   -2 & 2 \\
	\end{array}
	\right)
	$, 
	let $D=\diag{2d, d}\,(d\in \mathbb{Z}_{>0})$ and let $\Omega=\{(1, 2)\}$. We have $c_1=2d, c_2=d, g_{12}=g_{21}=1, f_{12}=1, f_{21}=2$. Then, $\Pi=\fgpre{C}{D}$ is isomorphic to the $K$-algebra defined as the quiver
	$$
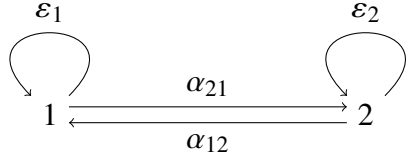
\begin{tikzpicture}[auto]
\node (a) at (0, 0) {$1$}; \node (b) at (4.2, 0) {$2$};
\draw[->, loop] (a) to node[swap] {$\eps{1}$} (a);
\draw[->, transform canvas={yshift=3pt}] (a) to node {$\alpha_{21}$} (b);
\draw[->, transform canvas={yshift=-3pt}] (b) to node {$\alpha_{12}$} (a);
\draw[->, loop] (b) to node[swap] {$\eps{2}$} (b);
\end{tikzpicture}$$
	with relations {(P1)} $\eps{1}^{2d}=0, \eps{2}^d=0$; {(P2)} $\eps{1}^2\alpha_{12}=\alpha_{12}\eps{2}, \eps{2}\alpha_{21}=\alpha_{21}\eps{1}^2$; {(P3)} $\alpha_{12}\alpha_{21}\eps{1}+\eps{1}\alpha_{12}\alpha_{21}=0, \alpha_{21}\alpha_{12}=0$.
\end{example}
As a $K$-algebra, a preprojective algebra of finite type is characterized by the following Proposition.
\begin{prop}[{Gei\ss-Leclerc-Schr\"oer~\cite[Corollary 11.3, 12.7]{GLS1}}]  \label{Piselfinj}
	Let $C$ be a connected GCM. Then, $\fgpre{C}{D}$ is a connected finite dimensional self-injective $K$-algebra if and only if $C$ is of finite type.
\end{prop}

In this paper, we mainly deal with a finite type GCM $C$ and its symmetrizer $D$. We put $D=mD'$ with the minimal symmetrizer $D'$. Let $\widehat{C}$ denote the corresponding un-twisted affine type GCM, and let $\widehat{D}=m\widehat{D'}$ denote a symmetrizer of $\widehat{C}$, where $\widehat{D'}$ is the minimal symmetrizer of $\widehat{C}$. Let $\hatPi\coloneqq \fgpre{\widehat{C}}{\widehat{D}}$ denote the generalized preprojective algebra associated with $\widehat{C}$ and $\widehat{D}$.

\begin{definition}[Locally free modules] \label{def of H_i&lfrep}
	Under the setting of Definition \ref{def of H&Pi}, we define:
	\begin{enumerate}
		\item $H_i \coloneqq e_i (K[\eps{i}]/(\eps{i}^{c_i})) e_i \cong  K[\eps{i}]/(\eps{i}^{c_i})$ for each $i \in Q_0$.
		\item A module $M$ over $H$ or $\Pi$ is \textit{locally free}, if $e_i M$ is a free $H_i$-module for each $i\in Q_0$.	
		\item For a locally free module $M$, we define $\rankvec{M}=(a_1, \dots, a_n)$ where $a_i\,\,(1\leq i \leq n)$ is the rank of $e_i M$ as a free $H_i$-module. For a free $H_i$-module $M$, we write the rank of $M$ as a free $H_i$-module by $\rk_i M$.
\end{enumerate}
\end{definition}
Let $\lfRep{\Lambda}$ (\resp $\lfrep{\Lambda}$) denote the full subcategory of $\Rep{\Lambda}$ consisting of locally free modules (\resp locally free modules such that each free $H_i$-module $e_i M$ is of finite rank) for $\Lambda = H$ or $\Pi$.
\begin{definition}
	For each $i\in Q_0$, we say that $E_i\in\lfrep{\Lambda}$ is a \textit{generalized simple module}, if we have $H_i$-module isomorphisms
	$$e_j E_i \cong\begin{cases}
		H_i &(j=i)\\
		0 &(j\neq i).
	\end{cases}$$
\end{definition}
By definition, $E_i$ is a uniserial module which has only simple modules $S_i$ as its composition factors. For $M, N\in \lfrep{H}$, we define
\begin{align*}
    \langle M, N \rangle_H &\coloneqq \dim \Hom{H}{M}{N} - \dim\Ext{1}{H}{M}{N},\\
    (M, N)_H&\coloneqq \langle M, N \rangle_H +\langle N, M \rangle_H,\\
    q_H(M)&\coloneqq \langle M, M \rangle_H.
\end{align*}
Then, we have the following proposition:
\begin{prop}[{\cite[Corollary 4.2]{GLS1}}]\label{prop:GLSprod}
 For $M,N\in \lfrep{H}$ with $\rankvec{M}= (a_1, \dots, a_n)$ and $\rankvec{N}= (b_1, \dots, b_n )$, we have
 $$\langle M, N\rangle_H =\sum_{i=1}^n c_i a_i b_i + \sum_{(j, i)\in \Omega} c_i c_{ij} a_i b_j.$$
 
Then, $(M, N)\mapsto \langle M, N \rangle_H$ descends to the Grothendieck group $\Z^n$ of $\lfrep{H}$ and induces a symmetric bilinear form on $(-, -)_H \colon \Z^n\times \Z^n\rightarrow \Z$. We have equalities $(-, -)_H=(-, -)_C$ and $q_H = q_{DC}$.
\end{prop}
We know the following properties about locally free modules.
\begin{prop}[{\cite[Proof of Lemma 3.8]{GLS1}}] \label{rep-closed}
	The category $\lfRep{\Pi}$ is closed under kernel of epimorphisms, cokernel of monomorphisms, and extensions. 
\end{prop}
Since each projective $\Pi$-module is locally free, we have the following:
\begin{prop}[{\cite[Corollary 2.7]{FG}}] \label{Pi_projdim}
	Let $M\in\Rep{\Pi}$. If $\pdim{\Pi}{M}<\infty$, then $M\in \lfRep{\Pi}$.
\end{prop}

\begin{prop}[{\cite[Theorem 12.6]{GLS1}}]\label{PiExtduality}
	For $M\in\lfRep{\Pi}$ and $N\in\lfrep{\Pi}$, we have the following:
	\begin{enumerate}		
\item We have a functorial isomorphism
$$\Ext{1}{\Pi}{M}{N} \cong \D{\Ext{1}{\Pi}{N}{M}}.$$
\item If $C$ does not contain any component of finite type, then we have more general functorial isomorphisms $$\Ext{2-i}{\Pi}{M}{N} \cong \D{\Ext{i}{\Pi}{N}{M}} {\,\,}(i=0, 1, 2) .$$
\item If $M$ is also of finite rank, we have the following equality for any type of symmetrizable GCM $C$:
$$\dim \Ext{1}{\Pi}{M}{N}=\dim\Hom{\Pi}{M}{N} + \dim \Hom{\Pi}{N}{M} - \left(M, N\right)_H.$$
We refer to this formula as the \textbm{generalized Crawley-Boevey formula}.
\end{enumerate}
\end{prop}

For classical preprojective algebras, a certain class of modules, called nilpotent modules, plays important roles in geometric representation theory about quantum groups. It is known that this class of modules can be characterized as the class of modules such that they have finite Jordan-H\"older filtrations with simple modules $\{S_i\}_{i\in Q_0}$ as composition factors. Gei\ss-Leclerc-Schr\"oer~\cite{GLS4} adapted this characterization as a first step of a generalization of nilpotent modules.
\begin{definition}\label{def:Efilt}
    A module $M\in \lfrep{\Pi}$ is $\E$-filtered if and only if $M$ has a filtration
    $$M=M_n \supsetneq M_{n-1} \supsetneq \cdots \supsetneq M_0 =0$$
    such that each $M_{i+1}/M_i$ is isomorphic to a generalized simple module $E_j$ $(j\in Q_0)$. Let $\Efilt{\Pi}$ denote the category of $\E$-filtered modules.
    \end{definition}
    \begin{rem}
    There is a locally free module which is not $\E$-filtered already for $\mathsf{B}_2$-type.
    \end{rem}
\subsection{Module varieties for generalized preprojective algebras}
We briefly recall representation varieties for generalized preprojective algebras from \cite[\S 2.2--\S 2.5]{GLS4}. Let $\Pi\coloneqq \fgpre{C}{D}$. For a dimension vector $d=(d_1, \dots, d_n)$, we define
$$\overline{H}(d)\coloneqq \prod_{a\in \overline{Q}_1}\Hom{K}{K^{d_{s(a)}}}{K^{d_{t(a)}}},$$ and let $\rep{(\Pi, d)}$ be the representation variety of $\Pi$-modules with dimension vector $d$. Then, by definition, a point in $\rep{(\Pi, d)}$ is a tuple $(M(a))_{a\in \overline{Q}_1} \in \overline{H}(d)$ which satisfies equations
$$M(\varepsilon_i)^{c_i}=0,\,\, M(\eps{i})^{f_{ji}}M(\alpha^{(g)}_{ij})=M(\alpha^{(g)}_{ij})M(\eps{j})^{f_{ij}},$$
$$\sum_{j \in \overline{\Omega}{\,}(-, i)} \sum_{g=1}^{g_{ij}} \sum_{f=0}^{{f_{ji}}-1} \sgn{i}{j}M(\eps{i})^f M(\alpha^{(g)}_{ij}) M(\alpha^{(g)}_{ji}) M(\eps{i})^{f_{ji}-1-f}=0$$
for all $i\in Q_0$, $(i, j)\in \overline{\Omega}$, and $1\leq g \leq g_{ij}$. The group $G(d)\coloneqq \prod_{i\in Q_0} \GL_K(d_i)$ acts on $\rep{(\Pi, d)}$ by conjugation. For a module $M\in \rep{(\Pi, d)}$, let $\mathcal{O}(M)\coloneqq G(d)M$ denote its $G(d)$-orbit. The set of $G(d)$-orbits in $\rep{(\Pi, d)}$ bijectively corresponds to the set of isoclasses of modules in $\rep{(\Pi, d)}$. Let $\lfrep{(\Pi, d)}\subseteq \rep{(\Pi, d)}$ be the subvarieties of locally free modules, and let $\nil_{
\E}(\Pi, d)\subseteq \lfrep{(\Pi, d)}$ be the constructible subset of $\E$-filtered modules in $\lfrep{(\Pi, d)}$.

We interpret above module varieties in terms of species (or sometimes called modulated graphs). We define $H_i$-$H_j$ bimodules as follows:
\begin{align*}
    {}_i H_j &\coloneqq H_i \tens{K} \spn_K 
    \{\alpha_{ij}^{(g)}\mid 1\leq g\leq g_{ij}\}\tens{K} H_j\\
    &= \spn_K \{ \varepsilon_i^{f_i}\tens{K} \alpha_{ij}^{(g)} \tens{K}\varepsilon_j^{f_j} \mid f_i, f_j \geq 0, 1\leq g \leq g_{ij} \}\\
    &= \bigoplus_{g=1}^{g_{ij}}\bigoplus_{f=0}^{f_{ij}-1}H_i \tens{K}(\alpha_{ij}^{(g)}\tens{K}\varepsilon_j^f) 
    = \bigoplus_{g=1}^{g_{ij}}\bigoplus_{f=0}^{f_{ji}-1}(\varepsilon_i^f \tens{K}\alpha_{ij}^{(g)}) \tens{K}H_j.
\end{align*}
In the below, we omit $\tens{K}$ for simplicity. In particular, we have isomorphisms as $H_i$-modules and as $H_j$-modules:
$${}_{H_i}({}_i H_j) \cong H_i^{|c_{ij}|} \cong ({}_j H_i) {}_{H_i};$$
$${}_{H_j}({}_j H_i) \cong H_j^{|c_{ji}|} \cong ({}_i H_j) {}_{H_j}.$$

Let  $M_i$ be a finite dimensional $H_i$-module for each $i \in Q_0$, and let $M_{ij}$ be $H_i$-linear map
$$M_{ij}\colon {}_i H_j\tens{H_j}M_j \rightarrow M_i$$
for each $(i,j)\in \overline{\Omega}$. Let $\otimes_{i}$ denote $\otimes_{H_i}$ for simplicity in the below.
Let $\M \coloneqq (M_1, \dots, M_n)$ denote a tuple of $H_i$-modules $M_i$
and let
\begin{align*}
H(\M)\coloneqq \prod_{(i,j)\in \Omega} &\Hom{H_i}{{}_iH_j\tens{j}M_j}{M_i},
&&H^*(\M)\coloneqq \prod_{(j,i)\in \Omega^*} \Hom{H_j}{{}_jH_i\tens{i}M_i}{M_j}.
\end{align*}
Let $M=(\M, \M_{ij}, \M_{ji})$ denote an object determined by a tuple of a pair of modules $\M$, elements $\M_{ij}\in H(\M)$ and $\M_{ji}\in H^*(\M)$. We define a morphism between two tuples $\M=(\M, \M_{ij}, \M_{ji})$ and $\bN=(\bN, \bN_{ij}, \bN_{ji})$ by a tuple $(f_i)_i$ of $H_i$-linear maps $f_i\colon M_i \rightarrow N_i$ such that there is a commutative diagram
$$\begin{tikzpicture}[auto]
\node (a) at (0, 0) {$M_i$}; 
\node (b) at (0, 1.8) {${}_i H_j\tens{H_j}M_j$};  
\node (c) at (4, 0) {$N_i$};
\node (d) at (4, 1.8) {${}_i H_j\tens{H_j}N_j$};
\draw[->] (b) to node[swap] {\footnotesize{$M_{ij}$}} (a);
\draw[->] (b) to node {\footnotesize{$1\tens{j} f_j$}} (d);
\draw[->] (a) to node {\footnotesize{$f_i$}} (c);
\draw[->] (d) to node {\footnotesize{$N_{ij}$}} (c);
\end{tikzpicture}$$
for each $i\in Q_0$.

Let $\rep{(C, D, \overline{\Omega})}$ denote the category defined by above data.
Now, we have a functorial isomorphism
$$\ad_{ji}\colon \Hom{H_j}{{}_jH_i\tens{i}M_i}{M_j}\rightarrow \Hom{H_i}{M_i}{{}_iH_j\tens{j}M_j}$$
$$\psi\mapsto \left(\psi^{\vee}\colon m\mapsto \sum_{g=1}^{g_{ij}}\sum_{f=0}^{f_{ji}-1}\varepsilon_i^f \alpha_{ij}^{(g)}\tens{j}  \psi(\alpha_{ji}^{(g)}\varepsilon_i^{f_{ji}-1-f}\tens{i} m )\right).$$
We write $M_{ji}^{\vee}\coloneqq \ad_{ji}(M_{ji})$. We define two maps $M_{i, \mathrm{in}}$ and $M_{i, \mathrm{out}}$ by
$$M_{i, \mathrm{in}} \coloneqq (\sgn{i} {j}M_{ij})_j\colon \bigoplus_{j\in \overline{\Omega}(i, -)}{}_iH_j\tens{j}M_j\rightarrow M_i$$
and
$$M_{i, \mathrm{out}}\coloneqq (M_{ji}^{\vee})_j\colon M_i \rightarrow \bigoplus_{j\in \overline{\Omega}(i, -)}{}_iH_j\tens{j}M_j.$$
Then, we have the following:
\begin{prop}[{\cite[Proposition 5.2]{GLS1}}]\label{prop:equivspecies}
Keep the above notation. The category $\rep{\Pi}$ is equivalent to the full subcategory of $\rep{(C, D, \overline{\Omega})}$ with objects $M=(\M, \M_{ij}, \M_{ji})$ such that $M_{i, \mathrm{in}}\circ M_{i, \mathrm{out}}=0$ for each $i\in Q_0$.
\end{prop}
By Proposition~\ref{prop:equivspecies}, we often identify objects in $\rep{\Pi}$ with those in $\rep{(C, D, \overline{\Omega})}$. Note that the category $\rep{\Pi}$ has the natural involution $*$ induced by taking the $K$-dual space $(M_i)^*$ of each $M_i$ and the $K$-dual map $(M_{ij})^*$ for each $M_{ij}$.

We write
$$\overline{H}(\M)\coloneqq\prod_{(i,j)\in \overline{\Omega}} \Hom{H_i}{{}_iH_j\tens{j}M_j}{M_i}.$$
From now on, we write the morphism in $\Hom{H_i}{{}_iH_j\tens{j}M_j}{M_i}$  by $M_{ij}\colon {}_iH_j\tens{j}M_j\rightarrow M_i$ for $M\in \overline{H}(\M)$. The group $G(\M)\coloneqq \prod_{i\in Q_0}\GL_{H_i}(M_i)$ acts by conjugation on $\overline{H}(\M)$, where $\GL_{H_i}(M_i)$ is the group of $H_i$-automorphisms of $M_i$. 
Now, we can think of $\rep{(\Pi, \M)}\coloneqq \{M\in \overline{H}(\M) \mid M_{i, \mathrm{in}}\circ M_{i, \mathrm{out}}=0 \}$ as an affine variety of $\Pi$-modules $M$ with equalities $e_i M= M_i$ for $i\in Q_0$. This $\rep{(\Pi, \M)}\subseteq \overline{H}(\M)$ is $G(\M)$-stable. The set of $G(\M)$-orbits in $\rep{(\Pi, \M)}$ bijectively corresponds to the set of isoclasses of $\Pi$-modules $M$ with equalities $e_i M=M_i$ for $i\in Q_0$. We call $\M$ \textit{locally free} if each $M_i$ is a free $H_i$-module. In this case, we write $\rankvec{\M}\coloneqq \left(\rk_1(M_1), \dots, \rk_n(M_n)\right)$ and call it the \textit{rank vector} of $\M$. For a locally free $\M$, we define $\nil_{\E}(\Pi, \M)= \Pi(\M)$ be the constructible subset of $\E$-filtered modules in $\rep{(\Pi, \M)}$. For a rank vector $\br=(r_1, \dots, r_n)$, we define
\begin{align*}
    \M(\br)&\coloneqq (H_1^{\oplus r_1}, \dots, H_n^{\oplus r_n});\\
    H(\br)&\coloneqq H(\M(\br));\\
    \overline{H}(\br)&\coloneqq \overline{H}(\M(\br));\\
    \rep{(\Pi, \br)}&\coloneqq\rep{(\Pi, \M(\br))};\\
    \nil_{\E}(\Pi, \br)&\coloneqq \Pi(\br)\coloneqq \Pi(\M(\br));\\
    G(\br)&\coloneqq G(\M(\br)).
\end{align*}
Then, each variety $\Pi(\M)$ is isomorphic to $\Pi(\br)$ for each $\br=\rankvec{\M}$ and we identify these varieties $\Pi(\br)$ and $\Pi(\M)$.

Finally, we compare these two constructions. We have the natural projection $\rep{\left(\Pi, d\right)}\xrightarrow{p_{\Pi}} \prod_{i\in Q_0} \rep{(H_i, d_i)}$.
Then, we have $p^{-1}_{\Pi}(\M)\cong \rep{(\Pi, \M)}$ for $\M = \left(M_1, \dots, M_n \right)$. In this setting, we find that the $G(\M)$-orbit of $M$ on this fiber can be written as $G(\M)M=G(d)M \cap \rep{(\Pi, \M)}$. Now, a closed $G(\M)$-stable subset $Z$ of $\rep{(\Pi, \M)}$ has $\dim Z =\dim G(\M)+m$ for some $m\in \Z$ if and only if the corresponding subset $G(d)Z$ of $\rep{(\Pi, d)}$ has $\dim G(d) +m$.
\subsection{Partial orders on module categories}
In this subsection, we briefly review a notion of $\tau$-rigid modules introduced in Auslander-Smal{\o}~\cite{AS} and later developed as the $\tau$-tilting theory due to Adachi-Iyama-Reiten~\cite{AIR}. Let $\Lambda$ be a basic finite dimensional $K$-algebra and $\tau$ be the Auslander-Reiten translation for $\rep{\Lambda}$. Let $\summand{M}$ denote the number of non-isomorphic indecomposable direct summands of $M$, and let $\langle e\rangle$ denote the two-sided ideal generated by an element $e\in\Lambda$.
\begin{definition} \label{ttilt-theory}
	Let $M\in\rep{\Lambda}$ and let $P\in\proj{\Lambda}$. We define:
	\begin{enumerate}
		\item $M$ is a \textit{$\tau$-rigid} $\Lambda$-module, if $\Hom{\Lambda}{M}{\tau M} =0$;
		\item $M$ is a \textit{$\tau$-tilting} $\Lambda$-module, if $M$ is $\tau$-rigid and $\summand{M} = \summand{\Lambda}$;
		\item $M$ is an \textit{almost complete $\tau$-tilting} $\Lambda$-module, if $M$ is $\tau$-rigid and $\summand{M} = \summand{\Lambda}-1$;
		\item $M$ is a \textit{support $\tau$-tilting} module, if there exists an idempotent $e \in \Lambda$ such that $M$ is a $\tau$-tilting $(\Lambda / \langle e \rangle )$-module;
		\item $(M, P)$ is a \textit{$\tau$-rigid pair}, if $M$ is $\tau$-rigid and $\Hom{\Lambda}{P}{M} =0$;
		\item $(M, P)$ is a \textit{support $\tau$-tilting pair} (\resp \textit{almost complete $\tau$-tilting pair}), if $(M, P)$ is a $\tau$-rigid pair and $\summand{M} + \summand{P} = \summand{\Lambda}$ (\resp $\summand{M} + \summand{P} = \summand{\Lambda}-1$);
		\item $(M, P)$ is a \textit{direct summand of $(M', P')$}, if $(M, P)$ and $(M', P')$ are $\tau$-rigid pairs and $M$ (\resp $P$) is a direct summand of $M'$  (\resp a direct summand of $P'$);
		\item $(M, P)$ is \textit{basic}, if $M$ and $P$ are basic (\ie each direct summand of $M\oplus P$ is multiplicity free).
	\end{enumerate}
	\end{definition}
 Let $\sttilt{\Lambda}$ denote the full subcategory of basic support $\tau$-tilting $\Lambda$-modules. We can think of $\tau$-rigid modules as a generalization of classical partial tilting modules in the sense of the classical Bongartz's lemma (see \cite[Theorem 2.10]{AIR}).
	Let $\add{M}$ denote the full subcategory of finite direct summands of finite direct sums of $M$. A characterization of $\tau$-rigid pairs and support $\tau$-tilting pairs is given as the following theorem:
	\begin{theorem}[{\cite[Proposition 2.3]{AIR}}] \label{tau-pair} Let $M \in \rep{\Lambda}$, let $P \in \proj{\Lambda}$, and let $e \in \Lambda$ be an idempotent element such that $\add{P} =\add{\Lambda e}$. 
	\begin{enumerate}
		\item $(M, P)$ is a $\tau$-rigid pair, if and only if $M$ is a $\tau$-rigid $(\Lambda /\langle e \rangle)$-module;
		\item  $(M, P)$ is a support $\tau$-tilting pair if and only if $M$ is a $\tau$-tilting $(\Lambda /\langle e \rangle)$-module;
		\item $(M, P)$ is an almost complete support $\tau$-tilting pair, if and only if $M$ is an almost complete $\tau$-tilting $(\Lambda /\langle e \rangle)$ module;
		\item If $(M, P)$ and $(M, Q)$ are support $\tau$-tilting pairs in $\rep{\Lambda}$, then $\add{P} = \add{Q}$.
	\end{enumerate}
\end{theorem}
By Theorem \ref{tau-pair}, we can identify basic support $\tau$-tilting modules with basic support $\tau$-tilting pairs. 
\begin{definition}[{{\cf} \cite[Proposition 1.1]{AIR}}]
\begin{enumerate}
	\item A full subcategory $\mathcal{T}$ in $\rep{\Lambda}$ (\resp $\mathcal{F}$ in $\rep{\Lambda}$) is a \textit{torsion class} (\resp a \textit{torsion-free class}), if $\mathcal{T}$ (\resp $\mathcal{F}$) is closed under extensions and taking a factor module of objects (\resp taking a submodule of objects).
	\item $\Fac{M}$ (\resp $\Sub{M}$) denotes the full subcategory of factor modules (\resp submodules) of finite direct sums of $M \in \rep{\Lambda}$.
\end{enumerate}
\end{definition}
Let $\tors{\Lambda}$ (\resp $\torf{\Lambda}$) denote the set of torsion classes in $\rep{\Lambda}$  (\resp torsion-free classes). In the $\tau$-tilting theory, one of the most important classes of algebras is $\tau$-tilting finite algebras:
\begin{definition}[\cite{DIJ}]\label{ttilt-finite}
	An algebra $\Lambda$ is called \textit{$\tau$-tilting finite}, if $\sttilt{\Lambda}$ is a finite set.
\end{definition}
\begin{theorem}[{\cite{DIJ}, \cite[Theorem 2.7]{AIR}}] \label{sttilt-tors} For a $\tau$-tilting finite algebra, we have a bijection between $\sttilt{\Lambda}$ and $\tors{\Lambda}$:
	\begin{align*}
	\sttilt{\Lambda} &\longrightarrow \tors{\Lambda}\\
	M &\longmapsto \Fac{M} \text{.}
	\end{align*}
\end{theorem}
We can define a partial order of $\sttilt{\Lambda}$ by Theorem \ref{sttilt-tors} as follows: 
\begin{definition}\label{poofsttilt}
	For $T, T' \in \sttilt{\Lambda}$, we define a partial order $\leq$ on $\sttilt{\Lambda}$ by $T \leq T'$ if $\Fac{T} \subseteq \Fac{T'}$.\end{definition}
	The above partial order is understood in terms of mutations:
	\begin{def-pro}[{\cite[Theorem 2.18]{AIR}}]\label{defofmutation}
		Any basic almost complete $\tau$-tilting pair $(U, Q)$ is a direct summand of precisely two different basic support $\tau$-tilting pairs $(T, P)$ and $(T', P')$. In addition, these $T, T'\in \sttilt{\Lambda}$ satisfy $T'<T$ or $T'>T$. In this setting, if $T'<T$ (\resp $T'>T$), we say that $(T', P')$ is a \textit{left} (\resp \textit{right}) \textit{mutation} of $(T, P)$. For this $T\in \sttilt{\Lambda}$ and the indecomposable summand $X$ of $T$ such that $T=X\oplus U$, we say that $T$ is the \textit{left} (\resp \textit{right}) \textit{mutation} of $T'$ at $X$, if $T'<T$ (\resp $T'>T$).
	\end{def-pro}
	\begin{prop}[{\cite[Definition-Proposition 2.28]{AIR}}]\label{equivcondofmutation}
	Under the setting of Definition-Proposition \ref{defofmutation}, $T'$ is the left mutation of $T$ at $X$ if and only if $X\notin \Fac{U}$.
	\end{prop}
	We know that for a $\tau$-tilting finite algebra $\Lambda$, the finite set $\sttilt{\Lambda}$ forms a finite complete lattice by Definition-Proposition \ref{ttilt-finite} and Definition \ref{poofsttiltofpi}. In particular, we find that $M, N \in \sttilt{\Lambda}$ are related by a mutation if and only if one is next to the other in the finite complete lattice of $\sttilt{\Lambda}$ by the following Theorem \ref{mutation&order}:
	\begin{theorem}[{\cite[Theorem 2.33]{AIR}\label{mutation&order}}]
	Let $T, U \in\sttilt{\Lambda}$. The following conditions are equivalent:
	\begin{enumerate}
		\item $U$ is a left mutation of $T$;
		\item $T$ is a right mutation of $U$;
		\item $T$ and $U$ satisfy $T>U$, and there does not exist $V\in\sttilt{\Lambda}$ such that $T>V>U$.
	\end{enumerate}
\end{theorem}
\subsection{Tilting theory for preprojective algebras}
In this subsection, we briefly review the work of Fu-Geng~\cite{FG} about a relationship between the generalized preprojective algebras and the Weyl groups of the Kac-Moody Lie algebras. Let $\Pi =\fgpre{C}{D}$ be the generalized preprojective algebra associated with a symmetrizable GCM $C$ and its symmetrizer $D$. We define an \textit{idempotent ideal} $I_i$ for each $i\in Q_0$ by $I_i\coloneqq \Pi (1-e_i)\Pi$. Note that we have an isomorphism $\Pi/ I_i \cong E_i$ for each $i \in Q_0$.
In \cite{FG}, some results of Buan-Iyama-Reiten-Scott~\cite{BIRS} and Mizuno~\cite{Miz} are generalized for our situation.
\begin{definition}
	\begin{enumerate}
		\item We say that a two-sided ideal $T$ of $\Pi$ is a \textit{tilting ideal} if $T$ is a left tilting $\Pi$-module and a right tilting $\Pi$-module in the sense of Miyasita~\cite{Miy}. 
		\item We assume $\Pi$ is finite dimensional over $K$. Then, we say that a two-sided ideal $T$ of $\Pi$ is a \textit{support $\tau$-tilting ideal} if $T$ is a left support $\tau$-tilting $\Pi$-module and a right support $\tau$-tilting $\Pi$-module.
	\end{enumerate}
\end{definition}
\begin{theorem}[{\cite[Lemma 3.2, 3.9, Theorem 3.12, 4.7, 5.14, 5.17]{FG}}] \label{bijection_W_tilt}
	 Let $C\in M_n(\mathbb{Z})$ be a symmetrizable GCM, and let $D$ be any symmetrizer of $C$.
	\begin{enumerate}
	\item There is a bijection $\psi$ from $W(C)$ to the monoid $$\langle I_1, I_2, \dots, I_n \rangle\coloneqq \left\{I_{i_1} I_{i_2} \cdots I_{i_k}\mid i_1, i_2, \dots, i_k \in Q_0, k\geq 0 \right\}$$
	given by $\psi{\,}(w) = I_w = I_{i_1} I_{i_2} \cdots I_{i_k}$, where $w = s_{i_1} s_{i_2} \cdots s_{i_k}$ is a reduced expression of $w\in W(C)$.
	Here, $\psi$ does not depend on the choice of reduced expression of $w$.
		\item Assume that $C$ has no components of finite type. Then, $T \in \langle I_1, I_2, \dots, I_n \rangle$ if and only if $T$ is a tilting ideal of $\Pi$ such that $\Pi/T$ has finite length. Any object in $\langle I_1, I_2, \dots, I_n \rangle \subseteq \Rep{\Pi}$ has projective dimension at most $1$.
		\item Assume that $C$ is of finite type. Then, $T \in \langle I_1, I_2, \dots, I_n \rangle$ if and only if $T$ is a basic support $\tau$-tilting ideal. In particular, $\psi\colon w\mapsto I_w$ gives a bijection between $W(C)$ and $\sttilt{\Pi}$.
	\end{enumerate}
	\end{theorem}
	\begin{example}
	    Let $C=\left(
	\begin{array}{rrr}
		2 & -1 & -1\\
	   -1 & 2 & -1\\
	   -1 & -1 & 2
	\end{array}
	\right)
	$, 
	let $D=\diag{1, 1, 1}$ and let $\Pi\coloneqq \fgpre{C}{D}$. Namely, we consider the following double quiver with preprojective relation:
	$$
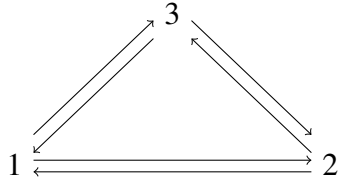
\begin{tikzpicture}[auto]
\node (a) at (0, 0) {$1$};
\node (b) at (4.2, 0) {$2$};
\node (c) at (2.1, 2) {$3$};
\draw[->, transform canvas={yshift=1.4pt}] (a) to (b);
\draw[->, transform canvas={yshift=-3.0pt}] (b) to (a);
\draw[->, transform canvas={yshift=4.2pt}] (c) to (b);
\draw[->, transform canvas={yshift=-2.6pt}] (b) to (c);
\draw[->, transform canvas={yshift=4.2pt}] (a) to (c);
\draw[->, transform canvas={yshift=-2.6pt}] (c) to (a);
\end{tikzpicture}.
$$
Then, the Loewy series of indecomposable projective modules $\Pi e_1$, $\Pi e_2$ and $\Pi e_3$ are respectively given by
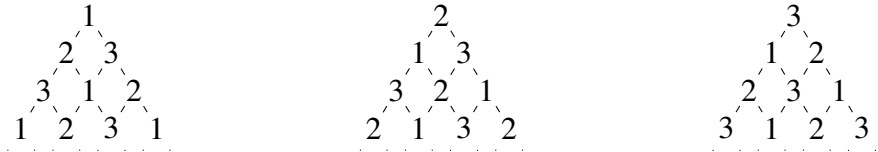
\begin{align*}
    &\begin{tikzpicture}[auto]
    \node (1) at (0, 0) {$1$};
    \node (2) at (-0.3, -0.5) {$2$};
    \node (3) at (0.3, -0.5) {$3$};
    \node (4) at (-0.6, -1) {$3$};
    \node (5) at (0, -1) {$1$};
    \node (6) at (0.6, -1) {$2$};
    \node (7) at (-0.9, -1.5) {$1$};
    \node (8) at (-0.3, -1.5) {$2$};
    \node (9) at (0.3, -1.5) {$3$};
    \node (10) at (0.9, -1.5) {$1$};
    \node (11) at (-1.2, -2) {};
    \node (12) at (-0.6, -2) {};
    \node (13) at (0, -2) {};
    \node (14) at (0.6, -2) {};
    \node (15) at (1.2, -2) {};
    \draw[-] (1) to (2);
    \draw[-] (1) to (3);
    \draw[-] (2) to (4);
    \draw[-] (2) to (5);
    \draw[-] (3) to (5);
    \draw[-] (3) to (6);
    \draw[-] (4) to (7);
    \draw[-] (4) to (8);
    \draw[-] (5) to (8);
    \draw[-] (5) to (9);
    \draw[-] (6) to (9);
    \draw[-] (6) to (10);
    \draw[dotted] (7) to (11);
    \draw[dotted] (7) to (12);
    \draw[dotted] (8) to (12);
    \draw[dotted] (8) to (13);
    \draw[dotted] (9) to (13);
    \draw[dotted] (9) to (14);
    \draw[dotted] (10) to (14);
    \draw[dotted] (10) to (15);
    \end{tikzpicture}
    &&\begin{tikzpicture}[auto]
    \node (1) at (0, 0) {$2$};
    \node (2) at (-0.3, -0.5) {$1$};
    \node (3) at (0.3, -0.5) {$3$};
    \node (4) at (-0.6, -1) {$3$};
    \node (5) at (0, -1) {$2$};
    \node (6) at (0.6, -1) {$1$};
    \node (7) at (-0.9, -1.5) {$2$};
    \node (8) at (-0.3, -1.5) {$1$};
    \node (9) at (0.3, -1.5) {$3$};
    \node (10) at (0.9, -1.5) {$2$};
    \node (11) at (-1.2, -2) {};
    \node (12) at (-0.6, -2) {};
    \node (13) at (0, -2) {};
    \node (14) at (0.6, -2) {};
    \node (15) at (1.2, -2) {};
    \draw[-] (1) to (2);
    \draw[-] (1) to (3);
    \draw[-] (2) to (4);
    \draw[-] (2) to (5);
    \draw[-] (3) to (5);
    \draw[-] (3) to (6);
    \draw[-] (4) to (7);
    \draw[-] (4) to (8);
    \draw[-] (5) to (8);
    \draw[-] (5) to (9);
    \draw[-] (6) to (9);
    \draw[-] (6) to (10);
    \draw[dotted] (7) to (11);
    \draw[dotted] (7) to (12);
    \draw[dotted] (8) to (12);
    \draw[dotted] (8) to (13);
    \draw[dotted] (9) to (13);
    \draw[dotted] (9) to (14);
    \draw[dotted] (10) to (14);
    \draw[dotted] (10) to (15);
    \end{tikzpicture}
    &&\begin{tikzpicture}[auto]
    \node (1) at (0, 0) {$3$};
    \node (2) at (-0.3, -0.5) {$1$};
    \node (3) at (0.3, -0.5) {$2$};
    \node (4) at (-0.6, -1) {$2$};
    \node (5) at (0, -1) {$3$};
    \node (6) at (0.6, -1) {$1$};
    \node (7) at (-0.9, -1.5) {$3$};
    \node (8) at (-0.3, -1.5) {$1$};
    \node (9) at (0.3, -1.5) {$2$};
    \node (10) at (0.9, -1.5) {$3$};
    \node (11) at (-1.2, -2) {};
    \node (12) at (-0.6, -2) {};
    \node (13) at (0, -2) {};
    \node (14) at (0.6, -2) {};
    \node (15) at (1.2, -2) {};
    \draw[-] (1) to (2);
    \draw[-] (1) to (3);
    \draw[-] (2) to (4);
    \draw[-] (2) to (5);
    \draw[-] (3) to (5);
    \draw[-] (3) to (6);
    \draw[-] (4) to (7);
    \draw[-] (4) to (8);
    \draw[-] (5) to (8);
    \draw[-] (5) to (9);
    \draw[-] (6) to (9);
    \draw[-] (6) to (10);
    \draw[dotted] (7) to (11);
    \draw[dotted] (7) to (12);
    \draw[dotted] (8) to (12);
    \draw[dotted] (8) to (13);
    \draw[dotted] (9) to (13);
    \draw[dotted] (9) to (14);
    \draw[dotted] (10) to (14);
    \draw[dotted] (10) to (15);
    \end{tikzpicture}.
\end{align*}
The Loewy series of $I_{s_1s_2}$ is given by the direct sum of the following three modules
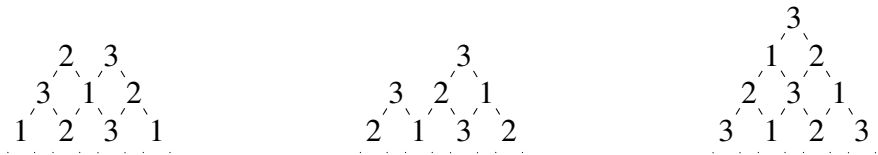
\begin{align*}
    &\begin{tikzpicture}[auto]
    \node (1) at (0, 0) {};
    \node (2) at (-0.3, -0.5) {$2$};
    \node (3) at (0.3, -0.5) {$3$};
    \node (4) at (-0.6, -1) {$3$};
    \node (5) at (0, -1) {$1$};
    \node (6) at (0.6, -1) {$2$};
    \node (7) at (-0.9, -1.5) {$1$};
    \node (8) at (-0.3, -1.5) {$2$};
    \node (9) at (0.3, -1.5) {$3$};
    \node (10) at (0.9, -1.5) {$1$};
    \node (11) at (-1.2, -2) {};
    \node (12) at (-0.6, -2) {};
    \node (13) at (0, -2) {};
    \node (14) at (0.6, -2) {};
    \node (15) at (1.2, -2) {};
    \draw[-] (2) to (4);
    \draw[-] (2) to (5);
    \draw[-] (3) to (5);
    \draw[-] (3) to (6);
    \draw[-] (4) to (7);
    \draw[-] (4) to (8);
    \draw[-] (5) to (8);
    \draw[-] (5) to (9);
    \draw[-] (6) to (9);
    \draw[-] (6) to (10);
    \draw[dotted] (7) to (11);
    \draw[dotted] (7) to (12);
    \draw[dotted] (8) to (12);
    \draw[dotted] (8) to (13);
    \draw[dotted] (9) to (13);
    \draw[dotted] (9) to (14);
    \draw[dotted] (10) to (14);
    \draw[dotted] (10) to (15);
    \end{tikzpicture}
    &&\begin{tikzpicture}[auto]
    \node (1) at (0, 0) {};
    \node (2) at (-0.3, -0.5) {};
    \node (3) at (0.3, -0.5) {$3$};
    \node (4) at (-0.6, -1) {$3$};
    \node (5) at (0, -1) {$2$};
    \node (6) at (0.6, -1) {$1$};
    \node (7) at (-0.9, -1.5) {$2$};
    \node (8) at (-0.3, -1.5) {$1$};
    \node (9) at (0.3, -1.5) {$3$};
    \node (10) at (0.9, -1.5) {$2$};
    \node (11) at (-1.2, -2) {};
    \node (12) at (-0.6, -2) {};
    \node (13) at (0, -2) {};
    \node (14) at (0.6, -2) {};
    \node (15) at (1.2, -2) {};
    \draw[-] (3) to (5);
    \draw[-] (3) to (6);
    \draw[-] (4) to (7);
    \draw[-] (4) to (8);
    \draw[-] (5) to (8);
    \draw[-] (5) to (9);
    \draw[-] (6) to (9);
    \draw[-] (6) to (10);
    \draw[dotted] (7) to (11);
    \draw[dotted] (7) to (12);
    \draw[dotted] (8) to (12);
    \draw[dotted] (8) to (13);
    \draw[dotted] (9) to (13);
    \draw[dotted] (9) to (14);
    \draw[dotted] (10) to (14);
    \draw[dotted] (10) to (15);
    \end{tikzpicture}
    &&\begin{tikzpicture}[auto]
    \node (1) at (0, 0) {$3$};
    \node (2) at (-0.3, -0.5) {$1$};
    \node (3) at (0.3, -0.5) {$2$};
    \node (4) at (-0.6, -1) {$2$};
    \node (5) at (0, -1) {$3$};
    \node (6) at (0.6, -1) {$1$};
    \node (7) at (-0.9, -1.5) {$3$};
    \node (8) at (-0.3, -1.5) {$1$};
    \node (9) at (0.3, -1.5) {$2$};
    \node (10) at (0.9, -1.5) {$3$};
    \node (11) at (-1.2, -2) {};
    \node (12) at (-0.6, -2) {};
    \node (13) at (0, -2) {};
    \node (14) at (0.6, -2) {};
    \node (15) at (1.2, -2) {};
    \draw[-] (1) to (2);
    \draw[-] (1) to (3);
    \draw[-] (2) to (4);
    \draw[-] (2) to (5);
    \draw[-] (3) to (5);
    \draw[-] (3) to (6);
    \draw[-] (4) to (7);
    \draw[-] (4) to (8);
    \draw[-] (5) to (8);
    \draw[-] (5) to (9);
    \draw[-] (6) to (9);
    \draw[-] (6) to (10);
    \draw[dotted] (7) to (11);
    \draw[dotted] (7) to (12);
    \draw[dotted] (8) to (12);
    \draw[dotted] (8) to (13);
    \draw[dotted] (9) to (13);
    \draw[dotted] (9) to (14);
    \draw[dotted] (10) to (14);
    \draw[dotted] (10) to (15);
    \end{tikzpicture}.
\end{align*}
The Loewy series of $I_{s_2s_1}$ is given by the direct sum of
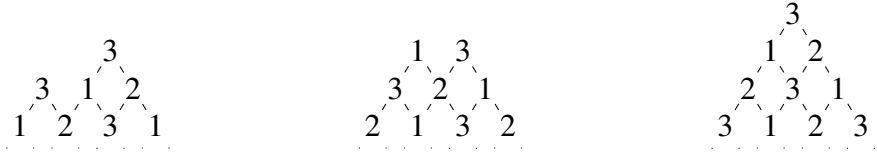
\begin{align*}
    &\begin{tikzpicture}[auto]
    \node (1) at (0, 0) {};
    \node (2) at (-0.3, -0.5) {};
    \node (3) at (0.3, -0.5) {$3$};
    \node (4) at (-0.6, -1) {$3$};
    \node (5) at (0, -1) {$1$};
    \node (6) at (0.6, -1) {$2$};
    \node (7) at (-0.9, -1.5) {$1$};
    \node (8) at (-0.3, -1.5) {$2$};
    \node (9) at (0.3, -1.5) {$3$};
    \node (10) at (0.9, -1.5) {$1$};
    \node (11) at (-1.2, -2) {};
    \node (12) at (-0.6, -2) {};
    \node (13) at (0, -2) {};
    \node (14) at (0.6, -2) {};
    \node (15) at (1.2, -2) {};
    \draw[-] (3) to (5);
    \draw[-] (3) to (6);
    \draw[-] (4) to (7);
    \draw[-] (4) to (8);
    \draw[-] (5) to (8);
    \draw[-] (5) to (9);
    \draw[-] (6) to (9);
    \draw[-] (6) to (10);
    \draw[dotted] (7) to (11);
    \draw[dotted] (7) to (12);
    \draw[dotted] (8) to (12);
    \draw[dotted] (8) to (13);
    \draw[dotted] (9) to (13);
    \draw[dotted] (9) to (14);
    \draw[dotted] (10) to (14);
    \draw[dotted] (10) to (15);
    \end{tikzpicture}
    &&\begin{tikzpicture}[auto]
    \node (1) at (0, 0) {};
    \node (2) at (-0.3, -0.5) {$1$};
    \node (3) at (0.3, -0.5) {$3$};
    \node (4) at (-0.6, -1) {$3$};
    \node (5) at (0, -1) {$2$};
    \node (6) at (0.6, -1) {$1$};
    \node (7) at (-0.9, -1.5) {$2$};
    \node (8) at (-0.3, -1.5) {$1$};
    \node (9) at (0.3, -1.5) {$3$};
    \node (10) at (0.9, -1.5) {$2$};
    \node (11) at (-1.2, -2) {};
    \node (12) at (-0.6, -2) {};
    \node (13) at (0, -2) {};
    \node (14) at (0.6, -2) {};
    \node (15) at (1.2, -2) {};
    \draw[-] (2) to (4);
    \draw[-] (2) to (5);
    \draw[-] (3) to (5);
    \draw[-] (3) to (6);
    \draw[-] (4) to (7);
    \draw[-] (4) to (8);
    \draw[-] (5) to (8);
    \draw[-] (5) to (9);
    \draw[-] (6) to (9);
    \draw[-] (6) to (10);
    \draw[dotted] (7) to (11);
    \draw[dotted] (7) to (12);
    \draw[dotted] (8) to (12);
    \draw[dotted] (8) to (13);
    \draw[dotted] (9) to (13);
    \draw[dotted] (9) to (14);
    \draw[dotted] (10) to (14);
    \draw[dotted] (10) to (15);
    \end{tikzpicture}
    &&\begin{tikzpicture}[auto]
    \node (1) at (0, 0) {$3$};
    \node (2) at (-0.3, -0.5) {$1$};
    \node (3) at (0.3, -0.5) {$2$};
    \node (4) at (-0.6, -1) {$2$};
    \node (5) at (0, -1) {$3$};
    \node (6) at (0.6, -1) {$1$};
    \node (7) at (-0.9, -1.5) {$3$};
    \node (8) at (-0.3, -1.5) {$1$};
    \node (9) at (0.3, -1.5) {$2$};
    \node (10) at (0.9, -1.5) {$3$};
    \node (11) at (-1.2, -2) {};
    \node (12) at (-0.6, -2) {};
    \node (13) at (0, -2) {};
    \node (14) at (0.6, -2) {};
    \node (15) at (1.2, -2) {};
    \draw[-] (1) to (2);
    \draw[-] (1) to (3);
    \draw[-] (2) to (4);
    \draw[-] (2) to (5);
    \draw[-] (3) to (5);
    \draw[-] (3) to (6);
    \draw[-] (4) to (7);
    \draw[-] (4) to (8);
    \draw[-] (5) to (8);
    \draw[-] (5) to (9);
    \draw[-] (6) to (9);
    \draw[-] (6) to (10);
    \draw[dotted] (7) to (11);
    \draw[dotted] (7) to (12);
    \draw[dotted] (8) to (12);
    \draw[dotted] (8) to (13);
    \draw[dotted] (9) to (13);
    \draw[dotted] (9) to (14);
    \draw[dotted] (10) to (14);
    \draw[dotted] (10) to (15);
    \end{tikzpicture}.
\end{align*}
The Loewy series of $I_{s_1s_2s_1}=I_{s_2s_1s_2}$ is given by the direct sum of
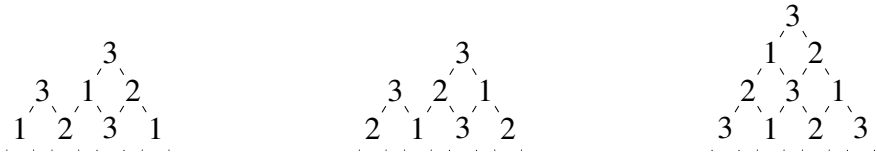
\begin{align*}
    &\begin{tikzpicture}[auto]
    \node (1) at (0, 0) {};
    \node (2) at (-0.3, -0.5) {};
    \node (3) at (0.3, -0.5) {$3$};
    \node (4) at (-0.6, -1) {$3$};
    \node (5) at (0, -1) {$1$};
    \node (6) at (0.6, -1) {$2$};
    \node (7) at (-0.9, -1.5) {$1$};
    \node (8) at (-0.3, -1.5) {$2$};
    \node (9) at (0.3, -1.5) {$3$};
    \node (10) at (0.9, -1.5) {$1$};
    \node (11) at (-1.2, -2) {};
    \node (12) at (-0.6, -2) {};
    \node (13) at (0, -2) {};
    \node (14) at (0.6, -2) {};
    \node (15) at (1.2, -2) {};
    \draw[-] (3) to (5);
    \draw[-] (3) to (6);
    \draw[-] (4) to (7);
    \draw[-] (4) to (8);
    \draw[-] (5) to (8);
    \draw[-] (5) to (9);
    \draw[-] (6) to (9);
    \draw[-] (6) to (10);
    \draw[dotted] (7) to (11);
    \draw[dotted] (7) to (12);
    \draw[dotted] (8) to (12);
    \draw[dotted] (8) to (13);
    \draw[dotted] (9) to (13);
    \draw[dotted] (9) to (14);
    \draw[dotted] (10) to (14);
    \draw[dotted] (10) to (15);
    \end{tikzpicture}
    &&\begin{tikzpicture}[auto]
    \node (1) at (0, 0) {};
    \node (2) at (-0.3, -0.5) {};
    \node (3) at (0.3, -0.5) {$3$};
    \node (4) at (-0.6, -1) {$3$};
    \node (5) at (0, -1) {$2$};
    \node (6) at (0.6, -1) {$1$};
    \node (7) at (-0.9, -1.5) {$2$};
    \node (8) at (-0.3, -1.5) {$1$};
    \node (9) at (0.3, -1.5) {$3$};
    \node (10) at (0.9, -1.5) {$2$};
    \node (11) at (-1.2, -2) {};
    \node (12) at (-0.6, -2) {};
    \node (13) at (0, -2) {};
    \node (14) at (0.6, -2) {};
    \node (15) at (1.2, -2) {};
    \draw[-] (3) to (5);
    \draw[-] (3) to (6);
    \draw[-] (4) to (7);
    \draw[-] (4) to (8);
    \draw[-] (5) to (8);
    \draw[-] (5) to (9);
    \draw[-] (6) to (9);
    \draw[-] (6) to (10);
    \draw[dotted] (7) to (11);
    \draw[dotted] (7) to (12);
    \draw[dotted] (8) to (12);
    \draw[dotted] (8) to (13);
    \draw[dotted] (9) to (13);
    \draw[dotted] (9) to (14);
    \draw[dotted] (10) to (14);
    \draw[dotted] (10) to (15);
    \end{tikzpicture}
    &&\begin{tikzpicture}[auto]
    \node (1) at (0, 0) {$3$};
    \node (2) at (-0.3, -0.5) {$1$};
    \node (3) at (0.3, -0.5) {$2$};
    \node (4) at (-0.6, -1) {$2$};
    \node (5) at (0, -1) {$3$};
    \node (6) at (0.6, -1) {$1$};
    \node (7) at (-0.9, -1.5) {$3$};
    \node (8) at (-0.3, -1.5) {$1$};
    \node (9) at (0.3, -1.5) {$2$};
    \node (10) at (0.9, -1.5) {$3$};
    \node (11) at (-1.2, -2) {};
    \node (12) at (-0.6, -2) {};
    \node (13) at (0, -2) {};
    \node (14) at (0.6, -2) {};
    \node (15) at (1.2, -2) {};
    \draw[-] (1) to (2);
    \draw[-] (1) to (3);
    \draw[-] (2) to (4);
    \draw[-] (2) to (5);
    \draw[-] (3) to (5);
    \draw[-] (3) to (6);
    \draw[-] (4) to (7);
    \draw[-] (4) to (8);
    \draw[-] (5) to (8);
    \draw[-] (5) to (9);
    \draw[-] (6) to (9);
    \draw[-] (6) to (10);
    \draw[dotted] (7) to (11);
    \draw[dotted] (7) to (12);
    \draw[dotted] (8) to (12);
    \draw[dotted] (8) to (13);
    \draw[dotted] (9) to (13);
    \draw[dotted] (9) to (14);
    \draw[dotted] (10) to (14);
    \draw[dotted] (10) to (15);
    \end{tikzpicture}.
\end{align*}
The Loewy series of $I_{s_1s_2s_3s_2s_1}$ is given by the direct sum of
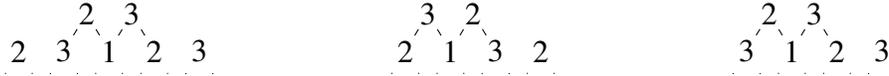
\begin{align*}
    &\begin{tikzpicture}[auto]
    \node (1) at (0, 0) {};
    \node (2) at (-0.3, -0.5) {};
    \node (3) at (0.3, -0.5) {};
    \node (4) at (-0.6, -1) {};
    \node (5) at (0, -1) {};
    \node (6) at (0.6, -1) {};
    \node (7) at (-0.9, -1.5) {};
    \node (8) at (-0.3, -1.5) {$2$};
    \node (9) at (0.3, -1.5) {$3$};
    \node (10) at (0.9, -1.5) {};
    \node (11) at (-1.2, -2) {$2$};
    \node (12) at (-0.6, -2) {$3$};
    \node (13) at (0, -2) {$1$};
    \node (14) at (0.6, -2) {$2$};
    \node (15) at (1.2, -2) {$3$};
    \node (16) at (-1.5, -2.5) {};
    \node (17) at (-0.9, -2.5) {};
    \node (18) at (-0.3, -2.5) {};
    \node (19) at (0.3, -2.5) {};
    \node (20) at (0.9, -2.5) {};
    \node (21) at (1.5, -2.5) {};
    \draw[-] (8) to (12);
    \draw[-] (8) to (13);
    \draw[-] (9) to (13);
    \draw[-] (9) to (14);
    \draw[dotted] (11) to (16);
    \draw[dotted] (11) to (17);
    \draw[dotted] (12) to (17);
    \draw[dotted] (12) to (18);
    \draw[dotted] (13) to (18);
    \draw[dotted] (13) to (19);
    \draw[dotted] (14) to (19);
    \draw[dotted] (14) to (20);
    \draw[dotted] (15) to (20);
    \draw[dotted] (15) to (21);
    \end{tikzpicture}
    &&\begin{tikzpicture}[auto]
    \node (1) at (0, 0) {};
    \node (2) at (-0.3, -0.5) {};
    \node (3) at (0.3, -0.5) {};
    \node (4) at (-0.6, -1) {$3$};
    \node (5) at (0, -1) {$2$};
    \node (6) at (0.6, -1) {};
    \node (7) at (-0.9, -1.5) {$2$};
    \node (8) at (-0.3, -1.5) {$1$};
    \node (9) at (0.3, -1.5) {$3$};
    \node (10) at (0.9, -1.5) {$2$};
    \node (11) at (-1.2, -2) {};
    \node (12) at (-0.6, -2) {};
    \node (13) at (0, -2) {};
    \node (14) at (0.6, -2) {};
    \node (15) at (1.2, -2) {};
    \draw[-] (4) to (7);
    \draw[-] (4) to (8);
    \draw[-] (5) to (8);
    \draw[-] (5) to (9);
    \draw[dotted] (7) to (11);
    \draw[dotted] (7) to (12);
    \draw[dotted] (8) to (12);
    \draw[dotted] (8) to (13);
    \draw[dotted] (9) to (13);
    \draw[dotted] (9) to (14);
    \draw[dotted] (10) to (14);
    \draw[dotted] (10) to (15);
    \end{tikzpicture}
    &&\begin{tikzpicture}[auto]
    \node (1) at (0, 0) {};
    \node (2) at (-0.3, -0.5) {};
    \node (3) at (0.3, -0.5) {};
    \node (4) at (-0.6, -1) {$2$};
    \node (5) at (0, -1) {$3$};
    \node (6) at (0.6, -1) {};
    \node (7) at (-0.9, -1.5) {$3$};
    \node (8) at (-0.3, -1.5) {$1$};
    \node (9) at (0.3, -1.5) {$2$};
    \node (10) at (0.9, -1.5) {$3$};
    \node (11) at (-1.2, -2) {};
    \node (12) at (-0.6, -2) {};
    \node (13) at (0, -2) {};
    \node (14) at (0.6, -2) {};
    \node (15) at (1.2, -2) {};
    \draw[-] (4) to (7);
    \draw[-] (4) to (8);
    \draw[-] (5) to (8);
    \draw[-] (5) to (9);
    \draw[dotted] (7) to (11);
    \draw[dotted] (7) to (12);
    \draw[dotted] (8) to (12);
    \draw[dotted] (8) to (13);
    \draw[dotted] (9) to (13);
    \draw[dotted] (9) to (14);
    \draw[dotted] (10) to (14);
    \draw[dotted] (10) to (15);
    \end{tikzpicture}.
\end{align*}
	\end{example}
	 Note that the multiplication gives an isomorphism
 $I_{i_1}\tens{\Pi}\cdots\tens{\Pi}I_{i_n}\cong I_{i_1}\cdots I_{i_n}$
if $w=s_{i_1}\cdots s_{i_n}$ is a reduced expression by the above Theorem \ref{bijection_W_tilt}.
		We note that we can obtain a similar result in $\Rep{\op{\Pi}}$ to Theorem \ref{bijection_W_tilt}. From now till the end of this subsection, let $C\in M_n(\mathbb{Z})$ be a symmetrizable GCM of finite type, and let $D$ be any symmetrizer $D$ of $C$.
	\begin{prop}[{\cite[Lemma 3.10, Theorem 3.11]{Mur}}]\label{I_w-Efilt}
	For each $w\in W$, the two sided ideal $I_w$ of $\fgpre{C}{D}$ is an $\mathbb{E}$-filtered module. In particular, any basic support $\tau$-tilting $\Pi$-module is $\mathbb{E}$-filtered. We have an isomorphism $I_w^{*} \cong \Pi / I_{w w_0}$ in $\lfrep{\Pi}$.
	\end{prop}
	 Since $C$ is of finite type, we have a relationship between the right weak Bruhat order $\leq_R$ on $W(C)$ and the mutation in $\sttilt{\Pi}$.
\begin{theorem}[{\cite[Lemma 5.11, Proposition 5.13 and its proof]{FG}}]\label{Tleftmutation}
	Let $T\in\langle I_1, \dots , I_n \rangle$. If $TI_i\neq T$, then $T$ has the left mutation $TI_i$ at $Te_i$ in $\sttilt{\Pi}$.
\end{theorem}
\begin{theorem}[{\cite[Theorem 5.16]{FG}}]\label{mutationI}
	For $i \in Q_0$ and $w \in W$, we write the corresponding $\tau$-tilting pair $(I_w, P_w)$.
	\begin{enumerate}
	\item $I_w, I_{ws_i}\in \sttilt{\Pi}$ are related by a right or left mutation. In particular, if $\length{ws_i}>\length{w}$, then we have
	$$I_{w s_i}=\begin{cases}
		I_w(1-e_i)\, &(I_wI_i=0)\\
		I_wI_ie_i\oplus I_w(1-e_i)\, &(I_wI_i\neq 0).
	\end{cases}$$
	\item We have the following isomorphism:
	$$P_w \cong \bigoplus_{i\in Q_0, I_w e_i=0} \Pi e_{\sigma(i)},$$
	where $\sigma$ is the Nakayama permutation determined by the self-injectivity of $\Pi$, \textit{i.e.} $\Pi e_i \cong \D(e_{\sigma(i)}\Pi)$.
	\end{enumerate}
\end{theorem}
Note that the Nakayama permutation $\sigma$ of $\Pi$ in the above Theorem \ref{mutationI} has a description in terms of the Weyl group $W(C)$, as we can understand by using $g$-matrices in the \S \ref{subsec:gvec}. We have the following generalization of a result of Mizuno~\cite{Miz}:
\begin{theorem}[{\cite[Theorem 2.31]{Mur}}]\label{poofsttiltofpi}
	Let $w\in W$ and $i\in Q_0$. The following are equivalent:
	\begin{itemize}
		\item $\length{w}<\length{ws_i}=\length{w}+1$;
		\item $I_w I_i \neq I_w$;
		\item $I_w$ has a left mutation $I_{ws_i}$ at $I_w e_i$.
	\end{itemize}
	\end{theorem}
	\begin{theorem}[{\cite[Theorem 2.32]{Mur} }]\label{ppoofsttiltofpi}
		Let $w\in W$ and $i\in Q_0$. The following are equivalent:
	\begin{itemize}
		\item $\length{w}>\length{ws_i}=\length{w}-1$;
		\item $I_w I_i = I_w$;
		\item $I_{ws_i}$ has a left mutation $I_w$ at $I_{ws_i}e_i$.
	\end{itemize}
	\end{theorem}
In particular, we find that $u \leq_R v$ in $W$ if and only if $I_u \geq I_v$ in $\sttilt{\Pi}$. That is, $(W, \leq_R)$ can be identified with $\op{(\sttilt{\Pi}, \leq )}$ as a poset. We note that we can consider the similar poset structure of $\sttilt{\op{\Pi}}$, if we consider the left weak order instead of the right weak Bruhat order. Finally, we recall homological properties about $I_w$.
\begin{prop}[{\cite[Proposition 3.8, 3.9]{Mur}}]\label{prpty-I_w}
We have the following equalities:
\begin{enumerate}
 \item If $\length{ws_i}> \length{w}$, then we have an equality $\Tor{1}{\Pi}{I_w}{E_i}=0$.
 \item If $\length{s_i w}> \length{w}$, then we have an equality $\Ext{1}{\Pi}{I_w}{E_i}=0$.
 \end{enumerate}
\end{prop}
\subsection{Reflection functors and idempotent ideals}\label{subsec:refl}
In the work of Baumann-Kamnitzer~\cite{BK} and Gei\ss-Leclerc-Schr\"oer~\cite{GLS1}, reflection functors on module categories of generalized preprojective algebras are developed in terms of species. We briefly review the definitions of reflection functors and their descriptions through idempotent ideals. 
\begin{definition}\label{def_reflection}
    Keep the notation of Proposition \ref{prop:equivspecies}. Let $M\in \rep{\Pi}$. We identify $M$ with $(\M, \M_{ij}, \M_{ji})$ such that $M_{i, \mathrm{in}}\circ M_{i, \mathrm{out}}=0$ for each $i\in Q_0$. We construct a new $\Pi$-module $\Sigma_i (M)$ by replacing the diagram
    $$\widetilde{M_i}\xrightarrow{M_{i, \mathrm{in}}}M_i\xrightarrow{M_{i, \mathrm{out}}}\widetilde{M_i}$$
    with
    $$\widetilde{M_i}\xrightarrow{\overline{M}_{i, \mathrm{out}}M_{i, \mathrm{in}}}\Ker M_{i, \mathrm{in}}\xrightarrow{\mathsf{can}}\widetilde{M_i},$$
    where $\widetilde{M}_i \coloneqq \bigoplus_{\overline{\Omega}(i,j)}{}_i H_j \tens{j} M_j$, and $\overline{M}_{i, \mathrm{out}}\colon M_i\longrightarrow\Ker M_{i, \mathrm{in}}$ is induced by $M_{i, \mathrm{out}}$, and $\mathsf{can}$ is the canonical injection. Similarly, we construct a new module $\Sigma_i^{-}(M)$ by replacing the diagram
    $$\widetilde{M_i}\xrightarrow{M_{i, \mathrm{in}}}M_i\xrightarrow{M_{i, \mathrm{out}}}\widetilde{M_i}$$
    with
    $$\widetilde{M_i}\xrightarrow{\mathsf{can}}\Cok M_{i, \mathrm{out}}\xrightarrow{M_{i, \mathrm{out}}\overline{M}_{i, \mathrm{in}}}\widetilde{M_i},$$
    where $\overline{M}_{i, \mathrm{in}}\colon \Cok M_{i, \mathrm{out}}\longrightarrow M_i$ is induced by $M_{i, \mathrm{in}}$, and $\mathsf{can}$ is the canonical projection.  Then, $(\Sigma_i(M))^* \cong \Sigma^-_i (M^*)$ and $(\Sigma_i^-(M))^*\cong \Sigma_i(M^*)$ for each $i\in Q_0$.
\end{definition}
\begin{rem}\label{rem:canmap}
By Definition \ref{def_reflection}, we have a diagram

$$\begin{tikzpicture}[auto]
\node (a) at (0, 0) {$\widetilde{M}_i$}; 
\node (b) at (4, 0) {$\Ker{M_{i, \mathrm{in}}}$};
\node (c) at (8, 0) {$\widetilde{M}_i$};
\node (d) at (0, 2) {$\widetilde{M}_i$};
\node (e) at (4, 2) {$M_i$};
\node (f) at (8, 2) {$\widetilde{M}_i$};
\node (g) at (0, 4) {$\widetilde{M}_i$};
\node (h) at (4, 4) {$\Cok{M_{i, \mathrm{out}}}$};
\node (i) at (8, 4) {$\widetilde{M}_i$};
\draw[->] (a) to node {\footnotesize{$\overline{M}_{i, \mathrm{out}}M_{i, \mathrm{in}}$}} (b);
\draw[->] (b) to node {\footnotesize{$\mathsf{can}$}} (c);
\draw[->] (d) to node {\footnotesize{$M_{i, \mathrm{in}}$}} (e);
\draw[->] (e) to node {\footnotesize{$M_{i, \mathrm{out}}$}} (f);
\draw[->] (g) to node {\footnotesize{$\mathsf{can}$}} (h);
\draw[->] (h) to node {\footnotesize{$M_{i, \mathrm{out}}\overline{M}_{i, \mathrm{in}}$}} (i);
\draw[double distance=2pt] (d) to (a);
\draw[double distance=2pt] (f) to (c);
\draw[double distance=2pt] (g) to (d);
\draw[double distance=2pt] (i) to (f);
\draw[->] (e) to node {\footnotesize{$\overline{M}_{i, \mathrm{out}}$}} (b);
\draw[->] (h) to node {\footnotesize{$\overline{M}_{i, \mathrm{in}}$}} (e);
\end{tikzpicture}.$$
By this diagram, we have canonical morphisms $\Sigma_i^- M \rightarrow M \rightarrow \Sigma_i M$.
\end{rem}
For $M\in \rep{\Pi}$ and $i\in Q_0$, let $\sub_i{M}$ (\textit{resp}. $\fac_i{M}$) be the largest submodule (\textit{resp}. factor module) of $M$ such that each composition factor of $\sub_i{M}$ (\textit{resp}. $\fac_i{M}$) is isomorphic to $S_i$. We say that $M$ has a trivial $i$-sub (\textit{resp}. $i$-factor) if $\sub_i{M}=0$ (\textit{resp}. $\fac_i{M}=0$). We define $\sub{M}\coloneqq\bigoplus_{i\in Q_0}\sub_i{M}$ and $\fac M\coloneqq\bigoplus_{i\in Q_0} \fac_i{M}$. The following properties about reflection functors are important.
\begin{theorem}[{\cite[Proposition 9.1]{GLS1} , K\"ulshammer~\cite[Lemma 5.2]{Kul}}]\label{refl}
For each $i\in Q_0$, the following statements \ref{refl1}--\ref{refl5} hold:
\begin{enumerate}
\item \label{refl1} The pair $(\Sigma_i^{-}, \Sigma_i)$ is a pair of adjoint functors, \ie there is a functorial isomorphism
$$\Hom{\Pi}{\Sigma_i^- M}{N}\cong \Hom{\Pi}{M}{\Sigma_i N}.$$
\item The adjunction morphisms $\mathrm{id}\rightarrow \Sigma_i\Sigma_i^{-}$ and $\Sigma_i^{-}\Sigma_i\rightarrow\mathrm{id}$ can be inserted in functorial short exact sequences
$$0\rightarrow \sub_i \rightarrow\mathrm{id}\rightarrow \Sigma_i \Sigma_i^{-}\rightarrow 0$$
and
$$0\rightarrow \Sigma_i^{-}\Sigma_i\rightarrow\mathrm{id}\rightarrow\fac_i\rightarrow 0.$$
\item We define subcategories $\mathcal{S}_i \coloneqq \{M \in\rep {\Pi} \mid \sub_i{M}=0 \}$ and $\mathcal{T}_i \coloneqq \{ M\in \rep{\Pi}\mid \fac_i{M}=0\}$. Then, we have mutual inverse equivalences of categories $\Sigma_i\colon \mathcal{T}_i \rightarrow \mathcal{S}_i$ and $\Sigma_i^- \colon \mathcal{S}_i \rightarrow \mathcal{T}_i$.
\item If $M$ is a locally free module, then we have an equality $\rankvec{\Sigma_i M} = s_i (\rankvec{M})$ (\resp $\rankvec{\Sigma_i^- M} = s_i (\rankvec{M})$) if and only if $\fac_i M= 0$ (\resp $\sub_i M=0$). 
\item \label{refl5} For any $i\in Q_0$, we have $\Sigma_i \cong \Hom{\Pi}{I_i}{\bullet}$ and $\Sigma_i^- \cong I_i \tens{\Pi}(\bullet)$.
\end{enumerate}
\end{theorem}

The following formulas proved in the work of Gei\ss-Leclerc-Schr\"oer~\cite{GLS4} are useful:
	\begin{prop}[\cite{GLS4}]\label{dimform1}
	    Let $M\in\rep{\Pi}$ and $i\in Q_0$. Then, we have the following:
	    \begin{enumerate}
	    \item $\dim\Hom{\Pi}{E_i}{M}=\dim \Ker M_{i, \mathrm{out}}= \dim\sub_i(M)$;
	    \item $\dim\Hom{\Pi}{M}{E_i}
	    =\dim\Cok{M_{i, \mathrm{in}}}
	    =\dim\fac_i(M)$;
	    \item If $M$ is locally free, then we have
	    \begin{align*}
	    \dim\Ext{1}{\Pi}{M}{E_i}&=\dim \widetilde{M}_i- \dim\Img{M_{i, \mathrm{in}}} -\dim\Img{M_{i, \mathrm{out}}}\\
	    &=\dim(\Ker M_{i, \mathrm{in}}/\Img{M_{i, \mathrm{out}}}).
	    \end{align*}
	    \end{enumerate}
	\end{prop}

Finally, we recall a notion of torsion pairs, which play important roles in this paper.
\begin{definition}
    Let $\mathcal{A}$ be an Abelian category. Then, a pair $(\calT, \calF)$ of extension closed subcategories is a torsion pair if and only if we have the following two equalities:
    \begin{itemize}
        \item $\calT=\{M\in \mathcal{A} \mid \Hom{\mathcal{A}}{M}{N}=0\,\,\textrm{for any}\, N\in \calF \}$;
        \item $\calF=\{N\in \mathcal{A} \mid \Hom{\mathcal{A}}{M}{N}=0\,\,\textrm{for any}\, M\in \calT \}$.
    \end{itemize}
\end{definition}

For a torsion pair $(\calT, \calF)$ and $M\in \rep{\Pi}$, we have the \textbm{canonical short exact sequence} $0\rightarrow tM \rightarrow M \rightarrow M/tM \rightarrow 0$ with torsion submodule $tM$ of $M$. Any short exact sequence of the form $0\rightarrow L \rightarrow M \rightarrow N \rightarrow 0$ with $L\in \calT$ and $N \in \calF$ is isomorphic to the canonical sequence (see \cite[\S VI Proposition 1.5]{ASS}).

For the module category over $\Pi$, we describe torsion pairs in terms of the Weyl group $W=W(C)$.
\begin{definition}
    Let $w\in W$. We define subcategories of $\rep{\Lambda}$ as follows:
    \begin{align}\label{eq:tors^w}
    &\calT^w=\{T \in \rep{\Pi} \mid \Ext{1}{\Pi}{I_w}{T}=0 \},
        &&\calF^w= \{T\in \rep{\Pi} \mid \Hom{\Pi}{I_w}{T}=0 \};
\end{align}
\begin{align}\label{eq:tors_w}
    &\calT_w=\{T \in \rep{\Pi} \mid I_w \tens{\Pi}{T}=0 \},
        &&\calF_w= \{T\in \rep{\Pi} \mid \Tor{1}{\Pi}{I_w}{T}=0  \}.
\end{align}
\end{definition}
\begin{prop}[{\cite[Corollary 3.12]{Mur}}]\label{prop:dualtors}
Let $C$ be a symmetrizable GCM of finite type.
\begin{enumerate}
    \item \label{eq:dualtors1} The natural $K$-duality $*$ on $\rep{\Pi}$ gives a bijection $*\colon \tors{\Pi}\rightarrow \torf{\Pi}$ by $\Fac{I_w}\mapsto \Sub{\Pi/I_{w w_0}}$.
    \item \label{eq:dualtors2} We have a bijection between $W$ and the set of torsion pairs given by
    $w\mapsto \left(\Fac{I_w}, \Sub{\Pi/I_{w}}\right)$
    and we can write
    \begin{align*}
        &\Fac{I_w}=\left\{I_w\tens{\Pi} T\mid T\in \rep{\Pi} \right\},
        &&\Sub{\Pi/I_w}= \left\{T\in \rep{\Pi} \mid \Hom{\Pi}{I_w}{T}=0 \right \}.
    \end{align*}
\end{enumerate}
\begin{proof}
For the reader's convenience, we give a sketch of the proof of \ref{eq:dualtors2}. The assertion $\Fac{I_w}=\left\{I_w\tens{\Pi} T\mid T\in \rep{\Pi} \right\}$ is a consequence of the Brenner-Butler theory. Namely, we have an algebra epimorphism $\hatPi \rightarrow \hatPi/\langle e_0 \rangle =\Pi$ for the generalized preprojective algebra of the un-twisted affine type associated with $\Pi$, so that we have an exact inclusion of categories $\rep{\Pi}\rightarrow \rep{\hatPi}$. Now, the $\hatPi$-modules $\affI{w}$ are tilting modules for $w\in W(C)$ by Theorem \ref{bijection_W_tilt}. So, we have functorial isomorphisms $\Ext{1}{\hatPi}{\affI{w}}{M} \cong \Ext{1}{\Pi}{I_w}{M}$ and $\Hom{\hatPi}{\affI{w}}{N} \cong \Hom{\Pi}{I_w}{N}$ for $M, N \in \rep{\Pi}$ \textit{etc}. In particular, we have two forms of torsion pairs $(\calT^w, \calF^w)$ and $(\calT_w, \calF_w)$. Note that the well-known fact in the $\tau$-tilting theory~\cite[Proposition 2.16]{AIR} tells us that $(\Fac{I_w}, I_w^\bot)$ is a torsion pair in $\rep{\Pi}$ for any support $\tau$-tilting module $I_w$, so that we have $\Fac{I_w}= \calT^w$. Any module $T$ has an epimorphism $\Pi^{\oplus m} \rightarrow T$ from a free $\Pi$-module $\Pi^{\oplus m}$, so that we have an epimorphism $I_w^{\oplus m}\rightarrow I_w \tens{\Pi} T$ by applying the functor $I_w \tens{\Pi}(-)$. So, any module of the form $I_w\tens{\Pi} T$ belongs to $\Fac{I_w}$. Now, we have a pair of exact functors whose compositions are identities
$$\begin{tikzpicture}[auto]
\node (a) at (0, 0) {$\calF_w$}; \node (b) at (4.2, 0) {$\calT^w$};
\draw[->, transform canvas={yshift=3pt}] (a) to node {\footnotesize{$I_w \tens{\Pi}(-)$}} (b);
\draw[->, transform canvas={yshift=-3pt}] (b) to node {\footnotesize{$\Hom{\Pi}{I_w}{-}$}} (a);
\end{tikzpicture}$$
by our exact inclusion of categories $\rep{\Pi}\rightarrow \rep{\hatPi}$.
So, $\calT^w$ is the essential image of $I_w \tens{\Pi}(-)$. 
Next, by the self-injectivity of $\Pi$, any ideal of $\Pi$ is basic if $\Pi$ is basic. So, the torsion submodule of $\Pi$ with respect to $\Fac{I_w}$ is isomorphic to $I_w$. Thus, the canonical short exact sequence of $\Pi$ associated with a torsion pair $(\calT^w, \calF^w)$ is isomorphic to
$$0 \rightarrow I_w \rightarrow \Pi \rightarrow \Pi/ I_w \rightarrow 0.$$
Since the bijection $W(C)\ni w \mapsto \Sub{\Pi/ I_w} \in \torf{\Pi}$ in Theorem \ref{sttilt-tors} and \ref{eq:dualtors1} implies that $(\Fac{I_w}, \Sub{\Pi/I_w})$ is a torsion pair by the finiteness of Weyl group $W(C)$. Thus, we obtain our assertion.
\end{proof}
\end{prop}

\begin{rem}\label{rem:torfF_w}
By a dual argument of Proposition \ref{prop:dualtors}, we find that the torsion pair $(\calT_w, \calF_w)$ in the Proof of Proposition \ref{prop:dualtors} coincides with the torsion pair $(\Fac{I_{w^{-1}w_0}}, \Sub{\Pi/ I_{w^{-1}w_0}})$ for each $w \in W$ and 
$$\calF_w= \left\{\Hom{\Pi}{I_w}{T} \mid T\in \rep{\Pi} \right \}.$$
\end{rem}

\begin{prop}\label{prop:intersecTF}
 Let $u, v \in W$. Suppose that  $u$ and $v$ satisfy $\length{uv}= \length{u}+ \length{v}$. Then, we have mutually inverse equivalences
 
 $$\begin{tikzpicture}[auto]
\node (a) at (0, 0) {$\calF^{v}$}; \node (b) at (6.2, 0) {$\calF^{uv} \cap \calT^{u}$};
\draw[->, transform canvas={yshift=3pt}] (a) to node {\footnotesize{$I_u \tens{\Pi}(-)$}} (b);
\draw[->, transform canvas={yshift=-3pt}] (b) to node {\footnotesize{$\Hom{\Pi}{I_u}{-}$}} (a);
\end{tikzpicture}.$$
\begin{proof}
If $T\in \calF^v$, then we have $T\in \calF_u$ by our choice of $u$ and $v$ by Remark \ref{rem:torfF_w}. In fact, we have $\calF^v= \Sub{\Pi/ I_v} \subset \Sub{\Pi/ I_{uv}}\subset \Sub{\Pi/ I_{uw_0}}= \calF_u$. Thus, we obtain
$$T\cong \Hom{\Pi}{I_u}{I_u\tens{\Pi}T}$$
by applying an equivalence $\Hom{\Pi}{I_u}{-}\colon \calT^u \rightarrow \calF_u$ for $I_u\tens{\Pi}T \in \calT^u$ by Proposition \ref{prop:dualtors}. We obtain
\begin{align*}
    \Hom{\Pi}{I_{uv}}{I_u\tens{\Pi}T}&\cong \Hom{\Pi}{I_v}{\Hom{\Pi}{I_u}{I_u\tens{\Pi}T}}\\
    &\cong\Hom{\Pi}{I_v}{T}=0.
\end{align*}
In particular, we obtain $I_u\tens{\Pi}T\in \calF^{uv}$. Conversely, if $T\in \calF^{uv}$, then
$$\Hom{\Pi}{I_v}{\Hom{\Pi}{I_u}{T}}=\Hom{\Pi}{I_{uv}}{T}=0.$$
Thus, we obtain $\Hom{\Pi}{I_u}{T}\in \calF^u$. Thus, we obtain our assertion by the equivalence between $\calF_u$ and $\calT^u$ induced by $\Hom{\Pi}{I_u}{-}$ and $I_u\tens{\Pi}(-)$.
\end{proof}
\end{prop}

The evaluation map $I_w \tens{\Pi} \left( \Hom{\Pi}{I_w}{M}\right) \rightarrow M$ is injective and the image is the torsion submodule of $M$ with respect to the torsion pair $(\Fac{I_w}, \Sub{\Pi/I_w})$. Let $T_w \coloneqq  I_w \tens{\Pi} \left( \Hom{\Pi}{I_w}{M}\right)$ denote the torsion submodule of $T$ with respect to a torsion pair $(\Fac{I_w}, \Sub{\Pi/I_w})$. Finally, we refer to a relationship between rank vectors and reflection functors.
\begin{prop}\label{prop:I_wrankvec}
 Let $w\in W$ and $M\in \lfrep{\Pi}$. If $M\in \calF_w$ (\resp $M\in \calT^w$), then $I_w\tens{\Pi} M$ (\resp $\Hom{\Pi}{I_w}{M}$) is a locally free module with $\rankvec{\,(I_w\tens{\Pi} M)} = w (\rankvec{M})$ (\resp $\rankvec{\,(\Hom{\Pi}{I_w}{M})} = w (\rankvec{M})$), where we identify $\rankvec{M}$ with some element of the root lattice.
 \begin{proof}
 Since we can prove the second assertion about $M\in \calT^w$ dually as the first assertion about $M\in \calF_w$, we only prove the first assertion. By Remark \ref{rem:torfF_w}, we put $M=\Hom{\Pi}{I_w}{T}$ with some $T\in \rep{\Pi}$. We prove by induction on $\length{w}$. If $\length{w}=1$ and we put $w = s_i$, then we have an equality 
 $$\dim\sub_i{M}=\dim\Hom{\Pi}{E_i}{\Hom{\Pi}{I_{s_i}}{T}}=\dim\Hom{\Pi}{I_{s_i}\tens{\Pi}E_i}{T}=0.$$
 Thus, the assertion is proved by Proposition \ref{dimform1} and \cite[Proposition 9.4]{GLS1}.
 
 We choose a reduced expression $w=s_{i_1}\cdots s_{i_k}$. Then, we have an isomorphism $\Hom{\Pi}{I_{w}}{(-)}\cong \Hom{\Pi}{I_{s_{i_k}}}{(\cdots (\Hom{\Pi}{I_{s_{i_1}}}{(-)}))}$ by $I_w \cong I_{s_{i_1}}\tens{\Pi}\cdots \tens{\Pi}{I_{s_{i_k}}}$. Now, the module $$I_{s_{i_1}w}\tens{\Pi}\Hom{\Pi}{I_{w}}{T} \cong I_{s_{i_1}w} \tens{\Pi} \Hom{\Pi}{I_{s_{i_1}w}}{\Hom{\Pi}{I_{s_{i_1}}}{T}}$$ is a torsion submodule of $\Hom{\Pi}{I_{s_{i_1}}}{T}$. We take a locally free module $N\coloneqq \Hom{\Pi}{I_w}{T}\cong \Hom{\Pi}{I_{s_{i_1}w}}{\Hom{\Pi}{I_{s_{i_1}}}{T}}$, so that $I_{s_{i_1}w} \tens{\Pi} N$ belongs to $\calT^{s_{i_1}w}$ and this module is locally free by induction hypothesis. Then, $I_{s_{i_1}w} \tens{\Pi} N$ is a submodule of an object in $\calF_{s_{i_1}}$. Since $\calF_{s_{i_1}}$ is a torsion-free class, any submodule of an object in $\calF_{s_{i_1}}$ belongs to $\calF_{s_{i_1}}$. Thus, we deduce our assertion inductively from the $\length{w}=1$ case because $I_w \tens{\Pi}N \cong I_w\tens{\Pi} \Hom{\Pi}{I_w}{T} \cong  I_{s_{i_1}}\tens{\Pi}I_{s_{i_1}w}\tens{\Pi}N$.
 \end{proof}
\end{prop}
Theorem \ref{refl} ensures that $\E$-filtered property is stable under taking reflection functors in some subcategories of $\lfrep{\Pi}$ in the following sense:
\begin{prop}
Let $M\in\Efilt{\Pi}$ such that $\fac_i M=0$ (\resp $\sub_i M=0$). Then, the locally free module $\Sigma_i M$ (\resp $\Sigma_i^- M$) is an $\E$-filtered module.
\begin{proof}
By Definition \ref{def_reflection}, we find that $\Sigma_i^- E_j \neq 0$ if and only if $i \neq j$. Now, we take an $\E$-filtered module $M$ and a filtration $0 = M_0 \subsetneq M_1 \subsetneq \cdots \subsetneq M_n =M$ such that $M_{i}/ M_{i-1} \cong E_{k_i}$ for each $1\leq i \leq n$ and some $k_i \in Q_0$. Let $m(M) \coloneqq n$ the length of this filtration. We prove our assertion by induction on $n$. If $m(M)=1$, then $M= E_j$ for some $j\in Q_0$. By the definition of $\Sigma_i^-$, we deduce $\Sigma_i^- E_j$ is $\E$-filtered by inspection. Next, we assume that our assertion is true for each $m(M)\leq l$ and take $m(M)=l+1$. Now, we write $M_{l+1}/M_l \cong E_{k_{l+1}}$. If $k_{l+1}\neq i$, then we have an exact sequence
$$0\rightarrow \Sigma_i^- M_l \rightarrow \Sigma_i^- M_{l+1} \rightarrow \Sigma_i^- E_{k_{l+1}} \rightarrow 0$$
by Proposition \ref{prpty-I_w}, Theorem~\ref{refl} and our choice of $k_{l+1}$.
Then, $\Sigma_i^- M_l$ is $\E$-filtered by our induction hypothesis. So, $\Sigma_i^- M_{l+1}$ is also $\E$-filtered by the first step of our induction and the above exact sequence. If $k=i$, then we have an epimorphism $f\colon \Sigma_i^- M_l \rightarrow \Sigma_i^- M_{l+1}$. Since $\lfrep{\Pi}$ is closed under taking kernels of epimorphisms, $\Ker{f}$ is a locally free module. Furthermore, $\Ker{f}$ is a direct sum of $E_i$ by the definition of $\Sigma_i^-$. So, we can deduce the assertion by the exact sequence
$$0\rightarrow \Ker{f} \rightarrow \Sigma_i^- M_l \rightarrow \Sigma_i^- M_{l+1}\rightarrow 0$$
and the induction hypothesis. Thus, we obtain our assertion by induction on $m(M)$.
\end{proof}
\end{prop}
\section{On stability of module categories}\label{section-stab}
Keep the setting of the previous section.
\subsection{Framed preprojective algebras and stable modules}\label{subsec-stab}
In this subsection, we develop representation theory of \textbm{framed preprojective algebras} as an analogue of the work of Baumann-Kamnitzer~\cite{BK}. Namely, we attach a symmetrizable GCM $\widetilde{C}$ and its symmetrizer $\widetilde{D}$
\begin{align*}
&\widetilde{C}=\left(
	\begin{array}{rr}
		C & -I_n \\
	   -I_n & 2I_n \\
	\end{array}
	\right)
	&&\widetilde{D}=\left(
	\begin{array}{rr}
		D & O \\
	   O & D \\
	\end{array}
	\right)
\end{align*}
	to any symmetrizable GCM $C=(c_{ij})$ and its symmetrizer $D$, and consider a new generalized preprojective algebra $\widetilde{\Pi}\coloneqq\fgpre{\widetilde{C}}{\widetilde{D}}$. A $\affPi$-module $M$ is a direct sum of $Q_0$-graded vector spaces $\bigoplus_{i\in Q_0}M_i$ and $\bigoplus_{i\in Q_0} M_{i'}$ where we set $\widetilde{Q}=(Q_0 \coprod Q_0', \widetilde{Q}_1)$ as the quiver associated with $(\widetilde{C}, \widetilde{D})$. The rank-vector of a locally free module $M$ is defined as a pair $\rankvec{M}=(\nu, \lambda)\in \mathbb{N}^{Q_0}\times \mathbb{N}^{Q_0}$, where $\nu_i= \rk_i M_i$ and $\lambda_{i'} =\rk_{i'}M_{i'}$. We will identify $\nu=(\nu_i)_{i\in Q_0}$ with a positive root $\sum_{i\in Q_0}\nu_i \alpha_i \in R$ and identify $\lambda=(\lambda_{i'})_{i\in Q_0}$ with an anti-dominant weight $-\sum_{i\in Q_0}\lambda_{i'} \varpi_i \in P$ (see \S \ref{pre:rootsys}). Under this identification and our choice of $\widetilde{C}$, we have an equality
	$$\langle (\nu, \lambda), (\alpha_i, 0) \rangle_{\widetilde{C}}=\langle \nu + \lambda, \alpha_i\rangle$$
	where the right hand side is the standard duality pairing between $P$ and $R$.
	We develop properties of stable modules defined as follows:
	\begin{definition}
	    Let $M\in \Efilt{\affPi}$. We say that $M$ is a stable module if $\sub_i {M}=0$ for each $i\in Q_0$.
	\end{definition}
	
	\begin{lemma}\label{lem:cineq}
	Let $M\in\lfrep{\affPi}$ be a stable module with $\rankvec{M}=(\nu, \lambda)$ and $c \coloneqq (\dim\fac_i M)/c_i$. Then, we have equalities $\dim \Hom{\affPi}{M}{E_i} = cc_i$ and $\dim \Ext{1}{\affPi}{M}{E_i} = c c_i-c_i\left\langle \nu+\lambda, \alpha_i\right\rangle$. In particular, we have an inequality $c\geq\max\left\{0, \left\langle\nu+\lambda, \alpha_i\right\rangle\right\}$.
	\begin{proof}
	By Proposition \ref{PiExtduality}, we have the following equality:
	$$\dim\Ext{1}{\affPi}{M}{E_i}=\dim\Hom{\affPi}{M}{E_i}+\dim\Hom{\affPi}{E_i}{M}-c_i\langle \rankvec{M}, \alpha_i\rangle_{\widetilde{C}}.$$
	Since $M$ is a stable module, Proposition \ref{dimform1} yields our statement by identifying $\langle \rankvec{M}, \alpha_i\rangle_{\widetilde{C}} = \langle \nu + \lambda, \alpha_i \rangle$.
	\end{proof}
	\end{lemma}
	By Lemma~\ref{lem:cineq}, we have the inequality $c \geq \max\{0, \langle\nu+\lambda, \alpha_i\rangle \}$. We define the following:
	\begin{definition}[Small $i$-factor]
	    Let $M\in\lfrep{\affPi}$ be stable. We say that $M$ has a \textit{small $i$-factor} if and only if $\rk_i M=\max\{0, \langle\nu+\lambda, \alpha_i\rangle\}$. 
	\end{definition}
	\begin{lemma}\label{stab-lem}
	Let $M\in\lfrep{\affPi}$ be stable. Then, we have the following:
	\begin{enumerate}
	\item \label{stablem1} If $M$ has $\rankvec M=(\nu, \lambda)$ and $\dim\fac_i M=cc_i$, then there exists a stable $\affPi$-module $M'$ with $\rankvec M'=(\nu +k\alpha_i, \lambda)$ for any $k \in \mathbb{Z} \cap [-c, c-\langle \nu+\lambda, \alpha_i\rangle]$.
	\item \label{stablem2} If $\langle \nu +\lambda, \alpha_i\rangle\leq 0$, then we have the following equivalence of categories:
	\begin{equation*}
\left\{\begin{gathered}
\text{stable $\affPi$-modules with}\\
\text{rank-vector $(\nu,\lambda)$}\\
\text{with small $i$-factor}
\end{gathered}\right\}
\xymatrix{\ar@<1ex>[r]^{\Sigma_i}&\ar@<1ex>[l]^{\Sigma_i^{-}}}
\left\{\begin{gathered}
\text{stable $\affPi$-modules with}\\
\text{rank-vector $(\nu-\langle\nu+\lambda,\alpha_i\rangle \alpha_i,\lambda)$}\\
\text{with small $i$-factor}
\end{gathered}\right\}.
\end{equation*}
	\item \label{stablem3} If there exists some stable module $M$ with $\rankvec M=(\nu, \lambda)$, then $\nu + \lambda \in \bigcap_{w\in W}w(\lambda + Q_+)$.
	\end{enumerate}
	\begin{proof}
	For $M\in\rep{\affPi}$, we consider this module in terms of species. Namely, we consider the following data based on Proposition \ref{prop:equivspecies}:
	\begin{align*}
	(\bigoplus_{j\in\overline{\Omega}(i, -)}{}_i H_j\tens{j}M_j)\oplus ({}_i H_{i'}\tens{i'} M_{i'})\xrightarrow{\sum \mathrm{sign}(i, j)M_{ij}-M_{ii'}} M_i\\ \xrightarrow{ (M^\vee_{ji}, M_{i'i}^{\vee})^\trp}(\bigoplus_{j\in\overline{\Omega}(i, -)}{}_i H_j\tens{j}M_j)\oplus ({}_i H_{i'}\tens{i'}M_{i'}).
	\end{align*}
	For simplicity, we write this diagram by $\widetilde{M_i}\xrightarrow{M_{i, \mathrm{in}}}M_i\xrightarrow{M_{i, \mathrm{out}}}\widetilde{M_i}$.
	We remark that $M$ is stable if and only if $M_{i, \mathrm{out}}$ is injective for any $i$ by Proposition \ref{dimform1}. If we put $u \coloneqq M_{i, \mathrm{out}}\circ M_{i, \mathrm{in}}\in\End{\affPi}{\widetilde{M}_i}$, then we have
	$$\Img{u} \subseteq \Img{M_{\mathrm{out}(i)}} \subseteq \Ker u$$
	and
	$u(\Img{M_{i, \mathrm{out}}})=0$ by a preprojective relation. Thus, we obtain the following equalities because $M$ is stable:
	$$\dim\Img{M_{i, \mathrm{out}}}-\dim\Img{u}=\dim\Cok M_{i, \mathrm{in}}=cc_i$$and
	\begin{align*}
	    &\dim\Ker{u} -\dim\Img{M_{i, \mathrm{out}}}\\=& \dim \widetilde{M}_i -\dim\Img{u}-\dim\Img{M_{i, \mathrm{out}}}\\
	    =&\dim\widetilde{M}_i+ (\dim\Img{M_{i, \mathrm{out}}}-\dim \Img{u})-2 \dim\Img{M_{i, \mathrm{out}}}\\
	    =& cc_i + \dim\widetilde{M}_i-2 \dim M_i\\
	    =& cc_i - c_i\langle \nu +\lambda, \alpha_i \rangle. 
	\end{align*}
	By these equalities, we obtain a submodule $V\subseteq \widetilde{M}_i$ of $\rk_i{\widetilde{M}_i}=k+\rk_i{M_i}$ such that $\Img{u}\subseteq V \subseteq \Ker{u}$. Thus, we can construct a module required in \ref{stablem1} by taking $\widetilde{M_i}\xrightarrow{u}V\xrightarrow{M_{i, \mathrm{out}}}\widetilde{M_i}$ instead of $\widetilde{M_i}\xrightarrow{M_{i, \mathrm{in}}}M_i\xrightarrow{M_{i, \mathrm{out}}}\widetilde{M_i}$.
	
	We prove \ref{stablem2}. We find that $M$ has a small $i$-factor if and only if $c=0$ or $c-\langle \nu +\lambda, \alpha_i \rangle =0$ in the Proof of \ref{stablem1}. This is equivalent to $\Img{M_{i, \mathrm{out}}} =\Img{u}$ or $\Ker{u}$. 
    
    If we assume $c=0$, then we have $\Img{M_{i, \mathrm{out}}}= \Img{u}$ and $\Ker{M_{i, \mathrm{in}}} =\Ker{u}$ by injectivity of $M_{i, \mathrm{out}}$. Then, by definition, the functor $\Sigma_i M$ is given by replacing the diagram
    $$\widetilde{M_i}\xrightarrow{M_{i, \mathrm{in}}}M_i\xrightarrow{M_{i, \mathrm{out}}}\widetilde{M_i}$$
    with
    $$\widetilde{M_i}\xrightarrow{\overline{M}_{i, \mathrm{out}}M_{i, \mathrm{in}}}\Ker u \xrightarrow{\mathsf{can}}\widetilde{M_i}.$$
    
    On the other hand,
    if we assume $c- \langle \nu +\lambda, \alpha_i \rangle =0$, then we have $\Img M_{i, \mathrm{out}} =\Ker u$ and $\Cok{\overline{M}_{i, \mathrm{out}}} =\widetilde{M}_i/\Ker{u} \cong \Img{u}$. Then,
    $\Sigma_i^{-}M$ is given by replacing the diagram
    $$\widetilde{M_i}\xrightarrow{M_{i, \mathrm{in}}}M_i\xrightarrow{M_{i, \mathrm{out}}}\widetilde{M_i}$$
    with
    $$\widetilde{M_i}\xrightarrow{\mathsf{can}} \Img{u} \xrightarrow{M_{i, \mathrm{out}}\overline{M}_{i, \mathrm{in}}}\widetilde{M_i}.$$
    
    This yields the assertion \ref{stablem2} from our equalities of dimensions.
    
    Finally, if we fix $\lambda$ and allow to vary $\nu$, then \ref{stablem1} implies that the set
    $$\{\nu +\lambda \in P \mid \text{there exists a stable module with rank-vector ($\nu$, $\lambda$)}\}\subseteq \lambda + Q_+$$
    is $W$-invariant. Therefore, the set is contained in $\bigcap_{w\in W}w(\lambda + Q_+)$.
	\end{proof}
	\end{lemma}
	By these lemmas, we have the following theorem as a generalization of \cite[Theorem 3.4]{BK}:
	\begin{theorem}\label{thm:stabexist}
	Let $\lambda\in \mathbb{Z}^n$ be an anti-dominant weight and $w\in W$. If $M$ is a stable $\affPi$-module with $\rankvec{M}=(w\lambda-\lambda, \lambda)$, then $M$ has a small $i$-factor for any $i\in Q_0$. Moreover, there exists a unique stable $\affPi$-module $\widetilde{N}(w\lambda)$ with $\rankvec \widetilde{N}(w\lambda)=(w\lambda-\lambda, \lambda)$. 
	\begin{proof}
	We prove the existence and the uniqueness later. Let $M$ be a stable module with $\rankvec{M}=(w\lambda -\lambda, \lambda)$ and $i\in Q_0$. First, we prove that $M$ has a small $i$-factor for any $i\in Q_0$.
	
	If $\length{s_i w} > \length{w}$, then there exists a stable $\affPi$-module $M'$ with $\rankvec{M'}=(w\lambda-\lambda -c \alpha_i, \lambda)$ by Lemma \ref{stab-lem} \ref{stablem1}. Now, we have $(w\lambda -\lambda -c\alpha_i)  +\lambda =w\lambda -c\alpha_i\in w(\lambda +R^+)$ by Lemma \ref{stab-lem} \ref{stablem3}. Thus, $-cw^{-1}(\alpha_i) \in R^+$. By $\length{s_i w}> \length{w}$, it is found that $w^{-1}(\alpha_i)$ is a positive root. So, $c$ must be equal to $0$. Since $\lambda$ is anti-dominant, we have $\langle w\lambda, \alpha_i \rangle \leq 0$.
	
	If $\length{s_i w}<\length{w}$, then there exists a stable module $M''$ with $\rankvec{M''}= (w\lambda -\lambda +(c-\langle w\lambda, \alpha_i \rangle )\alpha_i, \lambda))$ by Lemma \ref{stab-lem} \ref{stablem1}. Thus, $(w\lambda -\lambda +(c-\langle w\lambda, \alpha_i \rangle )\alpha_i)+ \lambda =w\lambda +(c-\langle w\lambda, \alpha_i \rangle )\alpha_i \in w(\lambda +R^+)$ by Lemma \ref{stab-lem} \ref{stablem3}. Since $w^{-1}(\alpha_i)$ is a negative root, we have $c-\langle w\lambda, \alpha_i \rangle \leq 0$. On the other hand, we have an inequality $c-\langle w\lambda, \alpha_i \rangle \geq 0$ by the existence of $M$ and Lemma \ref{lem:cineq}. Thus, we obtain $c=\langle w\lambda, \alpha_i \rangle$. So, $M$ has a small $i$-factor.
	
	We prove the existence and the uniqueness by induction on $\length{w}$. If $\length{w}=0$, the assertion follows from the fact that any two vertices in $Q'_0$ are not connected by arrows. Namely, a stable module $M$ with $\rankvec{M}=(0, \lambda)$ has only generalized simple modules for some vertices in $Q'_0$ as its direct summands. If $\length{w}>0$, then we can find $i \in Q_0$ such that $\length{s_i w}<\length{w}$. Then, the following
	\begin{equation*}
\left\{\begin{gathered}
\text{stable $\affPi$-modules with}\\
\text{rank-vector $(s_i w\lambda -\lambda,\lambda)$}\\
\end{gathered}\right\}
\xymatrix{\ar@<1ex>[r]^{\Sigma_i}&\ar@<1ex>[l]^{\Sigma_i^{-}}}
\left\{\begin{gathered}
\text{stable $\affPi$-modules with}\\
\text{rank-vector $(w\lambda-\lambda,\lambda)$}\\
\end{gathered}\right\}
\end{equation*}
gives an equivalence of categories by our first half argument and Lemma \ref{stab-lem} \ref{stablem2}. Thus, the assertion for $w$ is equivalent to that for $s_i w$. So, we obtain the assertion for any $w$ by induction.
	\end{proof}
	\end{theorem}
	By the above lemmas, we have a description of these $\widetilde{N}(\lambda)$ under the actions of reflection functors and dimension equalities as follows:
	\begin{cor}\label{cor:dimNhomext}
	Let $\gamma$ be a weight and $i\in Q_0$. 
	\begin{enumerate}
	  \item \label{dimhomext1} If $\langle \gamma, \alpha_i\rangle \leq 0$, then we have
	  \begin{align*}
	  &\Sigma_i\widetilde{N}(\gamma) \cong \widetilde{N}(s_i\gamma); &&\Sigma_i\widetilde{N}(s_i\gamma)\cong \widetilde{N}(s_i\gamma); &\Sigma_i^- \widetilde{N}(s_i\gamma)\cong \widetilde{N}(\gamma).
	  \end{align*}
	  \item \label{dimhomext2} We have the following equalities:
	  \begin{align*}
	      \dim \Hom{\affPi}{\widetilde{N}(\gamma)}{E_i}=\max\{0, c_i\langle \gamma, \alpha_i\rangle\};\\
	      \dim\Ext{1}{\affPi}{\widetilde{N}(\gamma)}{E_i}=\max\{0, -c_i\langle \gamma, \alpha_i\rangle\}.
	  \end{align*}
	  \end{enumerate}
	\begin{proof}
	\begin{itemize}
	\item[(1)] Let $w\in W$ such that $\length{s_i w}>\length{w}$ and let $\gamma =w\lambda$ with $\lambda$ anti-dominant. The first isomorphism follows from Theorem \ref{thm:stabexist}. 
	The second isomorphism follows from
	$$\Sigma_i \widetilde{N}(s_i \gamma)\cong \Sigma_i^2 \widetilde{N}(\gamma) \cong \Sigma_i \widetilde{N}(\gamma)\cong \widetilde{N}(s_i \gamma)\,\,\quad\textit{(cf. Proposition \ref{refl})}.$$
	Since $\widetilde{N}(\gamma)$ has a trivial $i$-factor, Theorem \ref{refl} yields
	$$\Sigma_i^- \widetilde{N}(s_i \gamma) \cong \Sigma^-_i \Sigma_i \widetilde{N}(\gamma) \cong \widetilde{N}(\gamma).$$
	
	\item[(2)] These equalities follow from Lemma~\ref{lem:cineq} and Theorem~\ref{thm:stabexist}. \qedhere
	\end{itemize}
	\end{proof}
	\end{cor}
	By examining the rank-vectors, we find that $\widetilde{N}(\lambda)$ is contained in $\widetilde{N}(w\lambda)$ as a submodule, and the canonical inclusion $\widetilde{N}(\lambda) \hookrightarrow \widetilde{N}(w\lambda)$ is an isomorphism at vertices in $Q'_0$. We define a new $\affPi$-module $N(w\lambda)$ from $\widetilde{N}(w\lambda)$ by $N(w\lambda)\coloneqq \widetilde{N}(w\lambda)/\widetilde{N}(\lambda)$. Since $\lfrep{\affPi}$ is closed under taking cokernels of epimorphisms, $N(w\lambda)$ is a locally free $\affPi$-module. This $N(w\lambda)$ has the natural $\Pi$-module structure because $N(w\lambda)\in\rep{\affPi}$ has dimension $0$ on any $Q'_0$ component. Thus, $N(w\lambda)$ is an $\E$-filtered module by our construction. The following theorem which is a generalization of \cite[Theorem 3.1]{BK} plays an important role in order to realize MV polytopes:
	\begin{theorem}\label{thm_stabmod}
	Let $\varpi_i$ be the fundamental weight with respect to $i\in Q_0$ and $w\in W$.
	\begin{enumerate}
	\item If $\gamma\in \mathbb{Z}^n$ is an anti-dominant weight, then we have $N(\gamma)=0$.
	\item If $\varpi_i -w\varpi_i\neq 0$, then there exists a unique $\Pi$-module $N(-w\varpi_i)\in \Efilt{\Pi}$ such that $\rankvec N(-w\varpi_i)=\varpi_i -w\varpi_i$ and $\sub N(-w\varpi_i)=E_i$. Conversely, if a module $M\in \rep{\Pi}$ satisfies $\rankvec M=\varpi_i -w\varpi_i$ and $\sub M=E_i$, then we have $M\cong N(-w\varpi_i)$.
	\item If $\gamma$ and $\delta$ belong to the same Weyl chamber, then we have $N(\gamma+ \delta)\cong N(\gamma)\oplus N(\delta)$.
	\end{enumerate}
	\begin{proof}
	\begin{itemize}
	\item[(1)] Since $\rankvec{\widetilde{N}(\gamma)}=(0, \gamma)$, the assertion is clear.
	
	\item[(2)] Let $w\in W$ and $i\in Q_0$. The module $\widetilde{N}(-w\varpi_i)$ is stable with $\rankvec{\widetilde{N}(-w\varpi_i)}=(\varpi_i - w\varpi_i, -\varpi_i)$. Now, a graded vector space $\bigoplus_{j\in Q'_0} \widetilde{N}(-w\varpi_i)_j$ coincides with $\widetilde{N}(-\varpi_i)\cong E_{i'}$. On the other hand, $N(-w\varpi_i)$ has $\rankvec{N(-w\varpi_i)}=\varpi_i -w\varpi_i$. For each $j\in Q_0$, we have $\sub_j{\widetilde{N}}(-w\varpi_i) = \Ker{M_{j'j}} \cap (\sub_j{N(-w\varpi_i)})$. Since ${\widetilde{N}}(-w\varpi_i)$ is stable, this equality yields that the inequality $\rk_{i} \sub_i{N(-w\varpi_i)} \leq 1$ and $\sub_j{N(-w\varpi_i)}=0$ for $i\neq j$. Since $N(-w\varpi_i)\neq 0$, we obtain $\sub{N(-w\varpi_i)}\cong E_i$. 
	
	Conversely, we assume that $M\in \rep{\Pi}$ satisfies $\rankvec{M}=\varpi_i -w\varpi_i$ and $\sub{M}=E_i$. Now, we can construct a stable module $N\in \rep{\affPi}$ such that $\rankvec{N}=(\varpi_i-w\varpi_i, -\varpi_i)$. In fact, we take
	\begin{align*}
	    &N_j= M_j\,\,\quad j\in Q_0\\
	    &N_h=M_h\,\,\quad h\in \overline{Q}_1\\
	    &N_{i'}=H_{i'}^{\oplus r}\,\,\quad(r\coloneqq \rk_i(\sub_i{M}))\\
	    &N_{i'i}\colon {}_{i'}H_i\tens{i} M_i \rightarrow N_{i'}
	\end{align*}
	where $N_{i'i}$ is any $H_{i'}$-isomorphism and the other arrows and spaces are zero. By the uniqueness of stable modules, we deduce $N\cong \widetilde{N}(-w\varpi_i)$. Thus, we obtain $M \cong N(-w\varpi_i)$.
	
	\item[(3)] By the uniqueness of stable modules, we have an isomorphism $\widetilde{N}(w\lambda)\oplus \widetilde{N}(w\mu)\cong \widetilde{N}(w(\lambda + \mu))$. Then, we take quotients to obtain $N(w\lambda)\oplus N(w\mu)\cong N(w(\lambda + \mu))$. \qedhere
	\end{itemize}
	\end{proof}
	\end{theorem}
	\begin{example}
	    Let $C =\left(
	\begin{array}{rr}
		2 & -1 \\
	   -2 & 2 \\
	\end{array}
	\right)
	$, 
	$D=\diag{2, 1}$. The double of the framed quiver $\widetilde{Q}$ is illustrated in Figure \ref{fig:my_label}.
	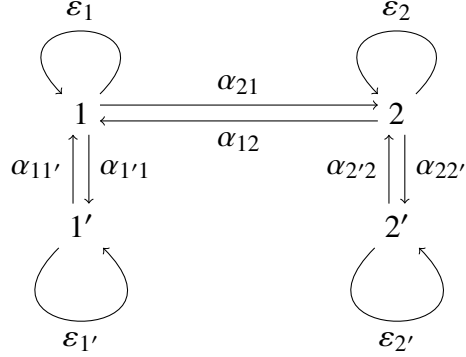
\begin{figure}
	\begin{tikzpicture}[auto]
\node (a) at (0, 0) {$1$}; \node (b) at (4.2, 0) {$2$};
\node (c) at (0, -1.5) {$1'$};
\node (d) at (4.2, -1.5) {$2'$};
\draw[->, out=45, in=135, loop] (a) to node[swap] {$\eps{1}$} (a);
\draw[->, transform canvas={yshift=3pt}] (a) to node {$\alpha_{21}$} (b);
\draw[->, transform canvas={yshift=-3pt}] (b) to node {$\alpha_{12}$} (a);
\draw[->, out=45, in=135, loop] (b) to node[swap] {$\eps{2}$} (b);
\draw[->, out=225, in=315, loop] (c) to node[swap] {$\eps{1'}$} (c);
\draw[->, out=225, in=315, loop] (d) to node[swap] {$\eps{2'}$} (d);
\draw[->, transform canvas={xshift=3pt}] (a) to node {$\alpha_{1'1}$} (c);
\draw[->, transform canvas={xshift=-3pt}] (c) to node {$\alpha_{11'}$} (a);
\draw[->, transform canvas={xshift=3pt}] (b) to node {$\alpha_{22'}$} (d);
\draw[->, transform canvas={xshift=-3pt}] (d) to node {$\alpha_{2'2}$} (b);
\end{tikzpicture}
 \centering
	    \caption{The framed quiver of type $\mathsf{B}_2$}
	    \label{fig:my_label}
	\end{figure}
The $\affPi$-module $\widetilde{N}(-\varpi_1)$ (\resp $\widetilde{N}(-\varpi_2)$) is isomorphic to a generalized simple module $E_{1'}$ (\resp $E_{2'}$). The $\Pi$-module $N(-s_1\varpi_1)$ is isomorphic to a generalized simple module $E_1$. The $\Pi$-modules $N(-s_2s_1\varpi_1)$, $N(-s_1s_2s_1\varpi_1)$ and $N(-s_2s_1s_2s_1\varpi_1 )$ are respectively illustrated as Figure \ref{fig:stab}. In particular, we have an isomorphism $N(-s_1s_2s_1\varpi_1) \cong N(-s_2s_1s_2s_1\varpi_1 )$. The $\Pi$-modules $N(-s_2s_1\varpi_2)$, $N(-s_1s_2s_1\varpi_2)$ and $N(-s_2s_1s_2s_1\varpi_2)$ are also illustrated in Figure \ref{fig:stab}.
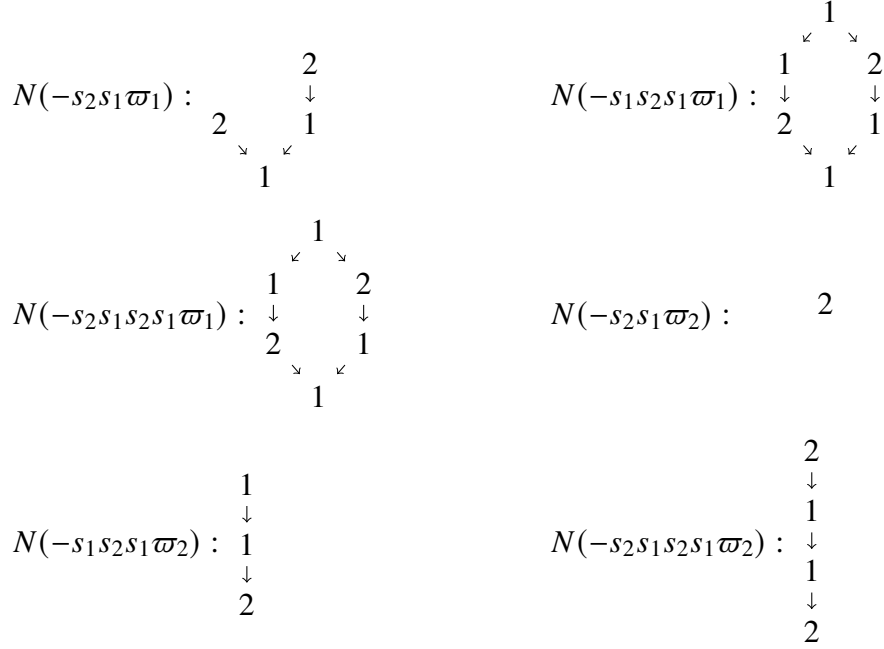
\begin{figure}
    \centering
\begin{align*}
    &N(-s_2s_1\varpi_1): \begin{tikzpicture}[auto,baseline=30pt]
    \node (1) at (0, 0.1) {$1$};
    \node (2) at (-0.6, 0.8) {$2$};
    \node (3) at (0.6, 0.8) {$1$};
    \node (4) at (0.6, 1.6) {$2$};
    \draw[->] (2) to (1);
    \draw[->] (3) to (1);
    \draw[->] (4) to (3);
    \end{tikzpicture}
    &&N(-s_1s_2s_1\varpi_1): \begin{tikzpicture}[auto,baseline=30pt]
    \node (1) at (0, 0.1) {$1$};
    \node (2) at (-0.6, 0.8) {$2$};
    \node (3) at (0.6, 0.8) {$1$};
    \node (4) at (0.6, 1.6) {$2$};
    \node (5) at (0, 2.3) {$1$};
    \node (6) at (-0.6, 1.6) {$1$};
    \draw[->] (2) to (1);
    \draw[->] (3) to (1);
    \draw[->] (4) to (3);
    \draw[->] (5) to (6);
    \draw[->] (6) to (2);
    \draw[->] (5) to (4);
    \end{tikzpicture}\\
    &N(-s_2s_1s_2s_1\varpi_1): \begin{tikzpicture}[auto,baseline=30pt]
    \node (1) at (0, 0.1) {$1$};
    \node (2) at (-0.6, 0.8) {$2$};
    \node (3) at (0.6, 0.8) {$1$};
    \node (4) at (0.6, 1.6) {$2$};
    \node (5) at (0, 2.3) {$1$};
    \node (6) at (-0.6, 1.6) {$1$};
    \draw[->] (2) to (1);
    \draw[->] (3) to (1);
    \draw[->] (4) to (3);
    \draw[->] (5) to (6);
    \draw[->] (6) to (2);
    \draw[->] (5) to (4);
    \end{tikzpicture}
    &&N(-s_2s_1\varpi_2): \begin{tikzpicture}
    \node (0) at (0, 0) {};
    \node (1) at (1, 0) {$2$};
    \end{tikzpicture}\\
    &N(-s_1s_2s_1\varpi_2): \begin{tikzpicture}[auto,baseline=20pt]
    \node (1) at (0, 0) {2};
    \node (2) at (0, 0.8) {1};
    \node (3) at (0, 1.6) {1};
    \draw[->] (2) to (1);
    \draw[->] (3) to (2);
    \end{tikzpicture}
    &&N(-s_2s_1s_2s_1\varpi_2): \begin{tikzpicture}[auto,baseline=30pt]
    \node (1) at (0, 0) {2};
    \node (2) at (0, 0.8) {1};
    \node (3) at (0, 1.6) {1};
    \node (4) at (0, 2.4) {2};
    \draw[->] (2) to (1);
    \draw[->] (3) to (2);
    \draw[->] (4) to (3);
    \end{tikzpicture}
\end{align*}
    \caption{Stable modules for type $\mathsf{B}_2$}
    \label{fig:stab}
\end{figure}
	\end{example}
	If $C$ is of finite type, we have the following:
	\begin{theorem}\label{projcover}
	The $\Pi$-module $N(\varpi_i)$ is isomorphic to the projective cover of a simple module $S_i$ for each $i\in Q_0$.
	\begin{proof}
	Note that the weight $w_0\varpi_i$ is anti-dominant. We have the canonical inclusion $\widetilde{N}(w_0\varpi_i)\hookrightarrow\widetilde{N}(\varpi_i)$ in $\rep{\affPi}$. In particular, we have an exact sequence
	$$0\rightarrow \widetilde{N}(w_0\varpi_i)\rightarrow \widetilde{N}(\varpi_i)\rightarrow N(\varpi_i)\rightarrow0.$$
	Applying the functor $\Hom{\affPi}{-}{E_j}$, we obtain
	\begin{align}
	\Hom{\affPi}{\widetilde{N}(w_0\varpi_i)}{E_j}\rightarrow \Ext{1}{\affPi}{N(\varpi_i)}{E_j}\rightarrow\Ext{1}{\affPi}{\widetilde{N}(\varpi_i)}{E_j}.\label{eq:seq1}
	\end{align}
	Here, both of two extremal terms in \eqref{eq:seq1} are equal to $0$ by Corollary \ref{cor:dimNhomext}. So, we obtain $\Ext{1}{\Pi}{N(\varpi_i)}{E_j}=\Ext{1}{\affPi}{N(\varpi_i)}{E_j}=0.$
	We remind that $E_j$ is a uniserial module. Namely, $E_j$ has a unique filtration $0=S^{(0)}_j\subsetneq S^{(1)}_j \subsetneq\cdots \subsetneq S^{(c_j)}_j=E_j$ of $\Pi$-modules. We have the canonical short exact sequence 
	$$0\rightarrow S_j^{(c_j-1)}\hookrightarrow E_j \rightarrow S_j\rightarrow 0.$$
	Then, we apply the functor $\Hom{\affPi}{N(\varpi_i)}{-}$ and obtain an exact sequence
	\begin{align}
	\Hom{\Pi}{N(\varpi_i)}{S_j}\rightarrow \Ext{1}{\Pi}{N(\varpi_i)}{S_j^{(c_j-1)}}\rightarrow\Ext{1}{\Pi}{N(\varpi_i)}{E_j}.\label{eq:seq2}
	\end{align}
	In the case of $i\neq j$, both of two extremal terms of \eqref{eq:seq2} are equal to $0$ by the above argument and the uniseriality of $E_j$. In fact, we have exact sequences for $X\in \rep{\Pi}$:
	$$0\rightarrow \Hom{\Pi}{X}{S_j^{(c_j-1)}} \rightarrow \Hom{\Pi}{X}{E_j} \rightarrow \Hom{\Pi}{X}{S_j}$$
	and
	$$0\rightarrow \Hom{\Pi}{X}{S_j^{(k)}}\rightarrow \Hom{\Pi}{X}{S_j^{(k+1)}}\rightarrow \Hom{\Pi}{X}{S_j}$$
	for $1 \leq k\leq c_j-1$.
	Thus, an induction on $k$ and $S_j^{(1)}\cong S_j$ gives $\Hom{\Pi}{X}{S_j}=0$ if $\Hom{\Pi}{X}{E_j}=0$. By a similar argument, we obtain $\Ext{1}{\Pi}{N(\varpi_i)}{S_j}=0$.
	
	In the case of $i=j$, the middle morphism of
	\begin{align*}
	\Hom{\Pi}{N(\varpi_i)}{E_i}&\rightarrow \Hom{\Pi}{N(\varpi_i)}{S_i}\\&\rightarrow \Ext{1}{\Pi}{N(\varpi_i)}{S_i^{(c_i-1)}}\rightarrow\Ext{1}{\Pi}{N(\varpi_i)}{E_i}=0
	\end{align*}
	is surjective. Since we have $\dim\Hom{\Pi}{N(\varpi_i)}{S_i}=1$, it is found that $\Ext{1}{\Pi}{N(\varpi_i)}{S_i^{(c_i-1)}}=0$. So, we obtain $\Ext{1}{\Pi}{N(\varpi_i)}{S_i}=0$ by induction. Thus, $N(\varpi_i)$ is a projective cover of a simple module. 
	On the other hand, we have an equality
	\begin{align*}
	\dim\Hom{\Pi}{N(\varpi_i)}{S_j}&=\dim\Hom{\affPi}{N(\varpi_i)}{S_j}\\
	&=\left(\dim\Hom{\affPi}{\widetilde{N}(\varpi_i)}{E_j}\right)\bigg/c_i=\langle \alpha_j, \varpi_i\rangle.
	\end{align*}So, we obtain our assertion.
	\end{proof}
	\end{theorem}
	We prepare the following modules and develop properties of them:
	\begin{definition}\label{def:M&V}
	    Let $w\in W$ and $w=s_{i_1}\cdots s_{i_l}$ be a reduced expression. We record this reduced expression by $\mathbf{i}\coloneqq (i_1, \dots i_l)$. We define $V_{\mathbf{i}, k} \coloneqq N\left(-s_{i_1}\cdots s_{i_k}\varpi_{i_k}\right)$ and $M_{\mathbf{i}, k}\coloneqq V_{\mathbf{i}, k}/V_{\mathbf{i}, k^-}$ for $1\leq k \leq l$ where $k^- =\max\{0, 1\leq s\leq k-1\mid i_s= i_k\}$.
	\end{definition}
	These kinds of modules have been studied in the works of Geiss-Leclerc-Schr\"oer~\cite{GLS7} and Buan-Iyama-Reiten-Scott~\cite{BIRS} for classical cases. We adopt the argument in the work of Baumann-Kamnitzer-Tingley~\cite[Example 5.12, 5.14]{BKT} and obtain the following theorem about these modules:
	\begin{theorem}\label{VM_iso}
	In the setting of Definition \ref{def:M&V}, we have isomorphisms
	$$V_{\mathbf{i}, k} \cong ((\Pi/I_{s_{i_1}\cdots s_{i_k}})\tens{\Pi} \Pi e_{i_k})^*$$
	and
	$$M_{\mathbf{i}, k} \cong (I_{s_{i_1}\cdots s_{i_{k-1}}} \tens{\Pi} E_{i_k})^* \cong (I_{s_{i_1}\cdots s_{i_{k-1}}}/I_{s_{i_1}\cdots s_{i_k}})^*.$$
	\begin{proof}
	First, we prove that $M_{\mathbf{i}, k} \cong (I_{s_{i_1}\cdots s_{i_{k-1}}} \tens{\Pi} E_{i_k})^*$. Since $\widetilde{N}(-s_{i_k}\varpi_{i_k}) \cong \Sigma_{i_k} \widetilde{N}(-\varpi_{i_k})$ and $\widetilde{N}(-\varpi_{i_k}) = E_{i'_k}$,
	we have the following exact sequence \eqref{eq:seq3} in $\rep{\affPi}$:
	\begin{align}
	0\rightarrow
	\widetilde{N}(-\varpi_{i_k})\rightarrow \widetilde{N}(-s_{i_k}\varpi_{i_k})\rightarrow E_{i_k}\rightarrow 0.\label{eq:seq3}
	\end{align}
	Since these three modules appearing in the sequence \eqref{eq:seq3} have trivial $i_{k-1}$-factors, the functor $\Sigma_{i_{k-1}}$ is exact on the above short exact sequence by Theorem \ref{refl}:
	$$0\rightarrow \widetilde{N}(-s_{i_{k-1}}\varpi_{i_k})\rightarrow \widetilde{N}(-s_{i_{k-1}}s_{i_k}\varpi_{i_k})\rightarrow \Sigma_{i_{k-1}}E_{i_k}\rightarrow0.$$
	Since $\widetilde{N}(-s_{i_{k-1}}\varpi_{i_k})$ and $\widetilde{N}(-s_{i_{k-1}}s_{i_k}\varpi_{i_k})$ have trivial $i_{k-2}$-factors by Corollary \ref{cor:dimNhomext} \ref{dimhomext2}, so does $\Sigma_{i_{k-1}}E_{i_k}$. Since $\mathbf{i}$ is a reduced expression, we can repeat these procedures by Corollary \ref{cor:dimNhomext} \ref{dimhomext2}. We obtain
	$$0\rightarrow \widetilde{N}(-s_{i_1}\cdots s_{i_{k-1}}\varpi_{i_k})\rightarrow \widetilde{N}(-s_{i_{1}}\cdots s_{i_k}\varpi_{i_k})\rightarrow \Sigma_{i_1}\cdots\Sigma_{i_{k-1}}E_{i_k}\rightarrow0.$$
	Now, $\widetilde{N}(-s_{i_{r}}\varpi_{i_k}) = \Sigma_{i_r}\widetilde{N}(-\varpi_{i_k})\cong \widetilde{N}(-\varpi_{i_k})$ for any $k^{-}<r<k$. So, we obtain
	$$0\rightarrow \widetilde{N}(-s_{i_1}\cdots s_{i_{k^-}}\varpi_{i_k})\rightarrow \widetilde{N}(-s_{i_{1}}\cdots s_{i_k}\varpi_{i_k})\rightarrow \Sigma_{i_1}\cdots\Sigma_{i_{k-1}}E_{i_k}\rightarrow0.$$
	The map $\widetilde{N}(-s_{i_{1}}\cdots s_{i_k}\varpi_{i_k})\rightarrow \Sigma_{i_1}\cdots\Sigma_{i_{k-1}}E_{i_k}$ sends the submodule $\widetilde{N}(-\varpi_{i_k})$ to $0$. So, we have the following short exact sequence in $\rep{\Pi}$:
	$$0\rightarrow V_{\mathbf{i}, k^-}\rightarrow V_{\mathbf{i}, k} \rightarrow \Sigma_{i_1}\cdots\Sigma_{i_{k-1}}E_{i_k} \rightarrow 0.$$
	So, we obtain 
	\begin{align*}
	M_{\mathbf{i}, k} =& V_{\mathbf{i}, k}/V_{\mathbf{i}, k^-} \\
	\cong& \Sigma_{i_1}\cdots\Sigma_{i_{k-1}}E_{i_k}\\
	\cong &(\Sigma_{i_1}^- \cdots\Sigma^-_{i_{k-1}}E_{i_k})^* \,\,\quad\textit{($E_i^* \cong E_i$)}\\
	\cong&  (I_{s_{i_1}\cdots s_{i_{k-1}}} \tens{\Pi} E_{i_k})^* \,\,\quad\textit{(Theorem \ref{refl})}.
	\end{align*}
	Then, we prove that
	$V_{\mathbf{i}, k} \equiv N(-s_{i_1}\cdots s_{i_k}\varpi_{i_k}) \cong ((\Pi/I_{s_{i_1}\cdots s_{i_k}})\tens{\Pi} \Pi e_{i_k})^*$. We pay attention to that $\Efilt{\Pi}$ is closed under taking quotients by submodules in $\Efilt{\Pi}$.
	For a reduced expression $w=s_{i_1}\cdots s_{i_l}$, we have a descending chain of ideals $\Pi\supsetneq I_{s_{i_1}}\supsetneq \cdots \supsetneq I_{s_{i_1}\cdots s_{i_l}}$. This filtration induces the filtration $\Pi/I_w \supsetneq I_{s_{i_1}}/I_w \supsetneq \cdots \supsetneq I_{s_{i_1}\cdots s_{i_l}}/I_w=0$. The $k$-th successive quotient in this filtration is isomorphic to $I_{s_{i_1}\cdots s_{i_{k-1}}}/I_{s_{i_1}\cdots s_{i_k}} \cong I_{s_{i_1}\cdots s_{i_{k-1}}}\tens{\Pi} (\Pi/I_{i_k})\cong I_{s_{i_1}\cdots s_{i_{k-1}}} \tens{\Pi} E_{i_k}$. By applying the exact functor $-\tens{\Pi}{\Pi e_i}$, all of the successive quotients $I_{s_{i_1}\cdots s_{i_{k-1}}}/I_{s_{i_1}\cdots s_{i_k}}$ in this filtration such that $i_k\neq i$ are killed by Theorem \ref{mutationI}. Thus, we obtain 
	\begin{align*}
	\rankvec{((\Pi/I_w)\tens{\Pi}{\Pi e_i})}=& \sum_{i_k=i} I_{s_{i_1}\cdots s_{i_{k-1}}} \tens{\Pi} E_{i_k} \\=&\sum_{i_k=i} s_{i_1}\cdots s_{i_{k-1}}\alpha_{i_k}\,\,\,\textit{(Proposition \ref{prop:I_wrankvec})} \\=&\varpi_i -w\varpi_i
	\end{align*}
	by a bit of calculation. Now, we have $\fac{\Pi e_i}\cong E_i$, so that $\fac{((\Pi/I_w)\tens{\Pi}{\Pi e_i})}$ is isomorphic to $0$ or $E_i$. So, the module $((\Pi/I_{s_{i_1}\cdots s_{i_k}})\tens{\Pi} \Pi e_{i_k})^*$ must coincide with $N(-s_{i_1}\cdots s_{i_k}\varpi_{i_k})$ by Theorem \ref{thm_stabmod}.
	\end{proof}
	\end{theorem}
\subsection{g-vectors and wall-chamber structures}\label{subsec:gvec}
\begin{figure}
            \centering
        \begin{tikzpicture}
        \node (1) at (0, 0) {$\tiny{
	\begin{array}{rrr}
		1 & 0 & 0\\
	   0 & 1 & 0\\
	   0 & 0 & 1
	\end{array}
	}$};
	\node (2) at (-2.5, -2) {$\tiny{
	\begin{array}{rrr}
		-1 & 0 & 0\\
	   1 & 1 & 0\\
	   0 & 0 & 1
	\end{array}
	}
	$};
	\node (3) at (0, -2) {$\tiny{
	\begin{array}{rrr}
		1 & 1 & 0\\
	   0 & -1 & 0\\
	   0 & 1 & 1
	\end{array}
	}
	$};
	\node (4) at (2.5, -2) {$\tiny{
	\begin{array}{rrr}
		1 & 0 & 0\\
	   0 & 1 & 1\\
	   0 & 0 & -1
	\end{array}
	}
	$};
	\node (5) at (-5, -4) {$\tiny{
	\begin{array}{rrr}
		-1 & -1 & 0\\
	   1 & 0 & 0\\
	   0 & 1 & 1
	\end{array}
	}
	$};
	\node (6) at (-2.5, -4) {$\tiny{
	\begin{array}{rrr}
		0 & 1 & 0\\
	   -1 & -1 & 0\\
	   0 & 1 & 1
	\end{array}
	}
	$};
	\node (7) at (0, -4) {$\tiny{
	\begin{array}{rrr}
		-1 & 0 & 0\\
	   1 & 1 & 1\\
	   0 & 0 & -1
	\end{array}
	}
	$};
	\node (8) at (2.5, -4) {$\tiny{
	\begin{array}{rrr}
		1 & 1 & 1\\
	   0 & -1 & -1\\
	   0 & 1 & 0
	\end{array}
	}
	$};
	\node (9) at (5, -4) {$\tiny{
	\begin{array}{rrr}
		1 & 1 & 0\\
	   0 & 0 & 1\\
	   0 & -1 & -1
	\end{array}
	}
	$};
	\node (10) at (-6.25, -6) {$\tiny{
	\begin{array}{rrr}
		0 & -1 & 0\\
	   -1 & 0 & 0\\
	   1 & 1 & 1
	\end{array}
	}
	$};
	\node (11) at (-3.75, -6) {$\tiny{
	\begin{array}{rrr}
		-1 & -1 & -1\\
	   1 & 0 & 0\\
	   0 & 1 & 0
	\end{array}
	}
	$};
	\node (12) at (-1.25, -6) {$\tiny{
	\begin{array}{rrr}
		0 & 1 & 1\\
	   -1 & -1 & -1\\
	   1 & 1 & 0
	\end{array}
	}
	$};
	\node (13) at (1.25, -6) {$\tiny{
	\begin{array}{rrr}
		-1 & -1 & 0\\
	   1 & 1 & 1\\
	   0 & -1 & -1
	\end{array}
	}
	$};
	\node (14) at (3.75, -6) {$\tiny{
	\begin{array}{rrr}
		1 & 1 & 1\\
	   0 & 0 & -1\\
	   0 & -1 & 0
	\end{array}
	}
	$};
	\node (15) at (6.25, -6) {$\tiny{
	\begin{array}{rrr}
		0 & 1 & 0\\
	   0 & 0 & 1\\
	   -1 & -1 & -1
	\end{array}
	}
	$};
	\node (16) at (-5, -8) {$\tiny{
	\begin{array}{rrr}
		0 & -1 & -1\\
	   -1 & 0 & 0\\
	   1 & 1 & 0
	\end{array}
	}
	$};
	\node (17) at (-2.5, -8) {$\tiny{
	\begin{array}{rrr}
		-1 & -1 & -1\\
	   1 & 1 & 0\\
	   0 & -1 & 0
	\end{array}
	}
	$};
	\node (18) at (0, -8) {$\tiny{
	\begin{array}{rrr}
		0 & 0 & 1\\
	   -1 & -1 & -1\\
	   1 & 0 & 0
	\end{array}
	}
	$};
	\node (19) at (2.5, -8) {$\tiny{
	\begin{array}{rrr}
		0 & -1 & 0\\
	   0 & 1 & 1\\
	   -1 & -1 & -1
	\end{array}
	}
	$};
	\node (20) at (5, -8) {$\tiny{
	\begin{array}{rrr}
		0 & 1 & 1\\
	   0 & 0 & -1\\
	   -1 & -1 & 0
	\end{array}
	}
	$};
	\node (21) at (-2.5, -10) {$\tiny{
	\begin{array}{rrr}
		0 & 0 & -1\\
	   -1 & -1 & 0\\
	   1 & 0 & 0
	\end{array}
	}
	$};
	\node (22) at (0, -10) {$\tiny{
	\begin{array}{rrr}
		0 & -1 & -1\\
	   0 & 1 & 0\\
	   -1 & -1 & 0
	\end{array}
	}
	$};
	\node (23) at (2.5, -10) {$\tiny{
	\begin{array}{rrr}
		0 & 0 & 1\\
	   0 & -1 & -1\\
	   -1 & 0 & 0
	\end{array}
	}
	$};
	\node (24) at (0, -12) {$\tiny{
	\begin{array}{rrr}
		0 & 0 & -1\\
	   0 & -1 & 0\\
	   -1 & 0 & 0
	\end{array}
	}
	$};
	\draw[->](1) to (2);
	\draw[->](1) to (3);
	\draw[->](1) to (4);
	\draw[->](2) to (5);
	\draw[->](2) to (7);
	\draw[->](3) to (6);
	\draw[->](3) to (8);
	\draw[->](4) to (7);
	\draw[->](4) to (9);
	\draw[->](5) to (10);
	\draw[->](5) to (11);
	\draw[->](6) to (10);
	\draw[->](6) to (12);
	\draw[->](7) to (13);
	\draw[->](8) to (12);
	\draw[->](8) to (14);
	\draw[->](9) to (14);
	\draw[->](9) to (15);
	\draw[->](10) to (16);
	\draw[->](11) to (16);
	\draw[->](11) to (17);
	\draw[->](12) to (18);
	\draw[->](13) to (17);
	\draw[->](13) to (19);
	\draw[->](14) to (20);
	\draw[->](15) to (20);
	\draw[->](15) to (19);
	\draw[->](16) to (21);
	\draw[->](17) to (22);
	\draw[->](18) to (21);
	\draw[->](18) to (23);
	\draw[->](19) to (22);
	\draw[->](20) to (23);
	\draw[->](21) to (24);
	\draw[->](22) to (24);
	\draw[->](23) to (24);
        \end{tikzpicture}
        \caption{$g$-matrices for type $\mathsf{A}_3$}
        \label{fig:A3}
        \end{figure}
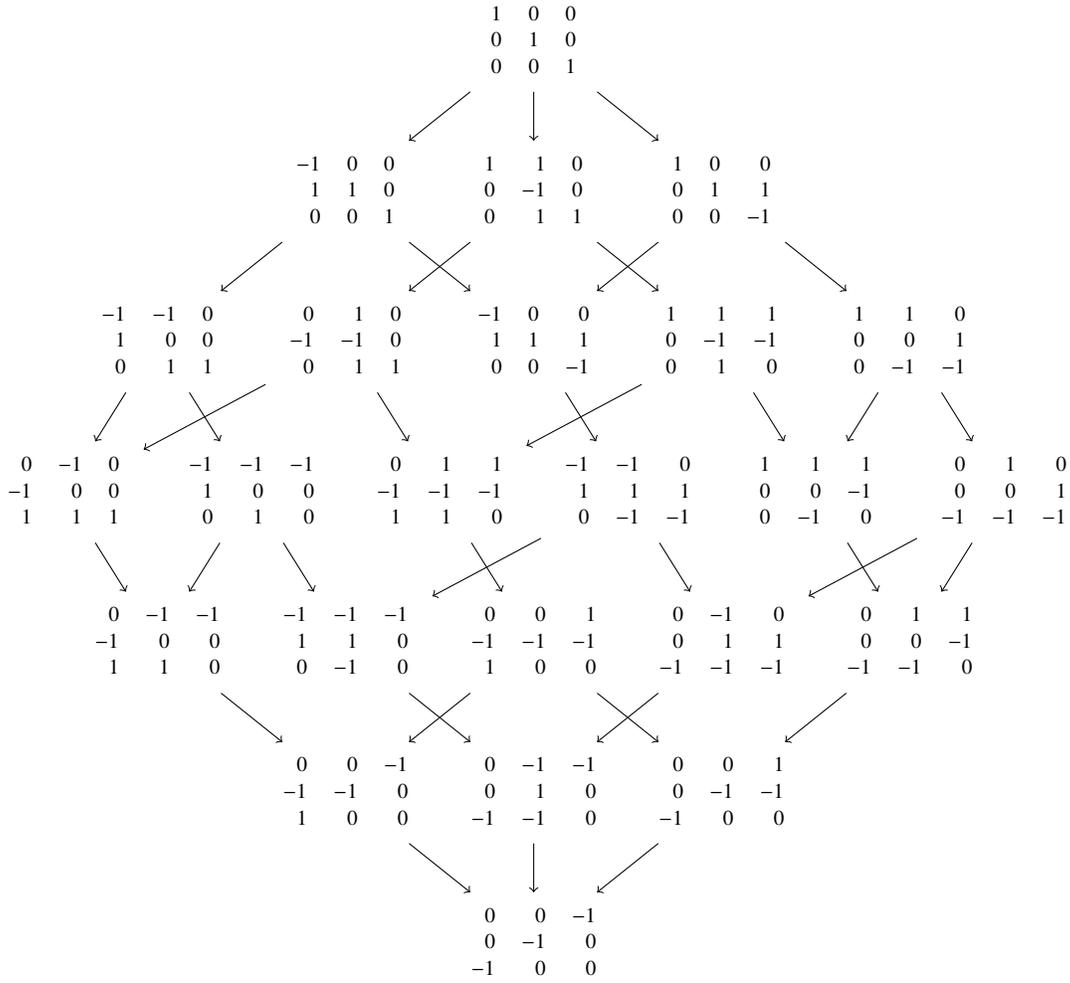

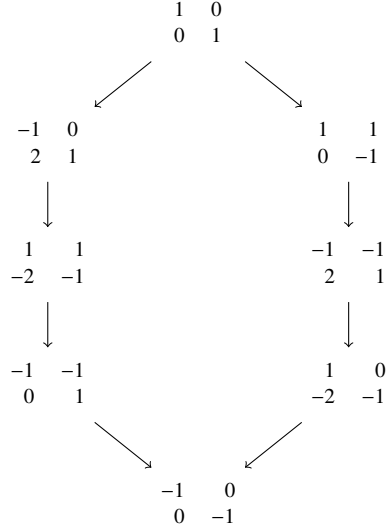
\begin{figure}
        \centering
    \begin{tikzpicture}
    \node (1) at (0, 0) {$\tiny{
	\begin{array}{rr}
		1 & 0 \\
	   0 & 1 
	\end{array}
	}
	$};
	\node (2) at (-2, -1.6) {$\tiny{
	\begin{array}{rr}
		-1 & 0 \\
	   2 & 1 
	\end{array}
	}
	$};
	\node (3) at (2, -1.6) {$\tiny{
	\begin{array}{rr}
		1 & 1 \\
	   0 & -1 
	\end{array}
	}
	$};
	\node (4) at (-2, -3.2) {$\tiny{
	\begin{array}{rr}
		1 & 1 \\
	   -2 & -1 
	\end{array}
	}
	$};
	\node (5) at (2, -3.2) {$\tiny{
	\begin{array}{rr}
		-1 & -1 \\
	   2 & 1 
	\end{array}
	}
	$};
	\node (6) at (-2, -4.8) {$\tiny{
	\begin{array}{rr}
		-1 & -1 \\
	   0 & 1 
	\end{array}
	}
	$};
	\node (7) at (2, -4.8) {$\tiny{
	\begin{array}{rr}
		1 & 0 \\
	   -2 & -1 
	\end{array}
	}
	$};
	\node (8) at (0, -6.4) {$\tiny{
	\begin{array}{rr}
		-1 & 0 \\
	   0 & -1 
	\end{array}
	}
	$};
	\draw[->] (1) to (2);
	\draw[->] (1) to (3);
	\draw[->] (2) to (4);
	\draw[->] (3) to (5);
	\draw[->] (4) to (6);
	\draw[->] (5) to (7);
	\draw[->] (6) to (8);
	\draw[->] (7) to (8);
    \end{tikzpicture}
    \caption{$g$-matrices for type $\mathsf{B}_2$}
        \label{fig:B2}
    \end{figure}

In this subsection, we develop $g$-vectors for finite dimensional generalized preprojective algebras in order to realize MV polytopes. Our discussion is an analogue of the work of Mizuno~\cite{Miz}. In this subsection, $\Pi$ means a finite dimensional generalized preprojective algebra, and $\hatPi$ means the completed generalized preprojective algebra (\cf \cite[Proposition 12.3]{GLS1}) of un-twisted affine type associated with $\Pi$. Here, this completion is with respect to the two-sided ideal generated by $\bigoplus_{(i,j)\in \overline{\Omega}} {}_i H_j$. In particular, $\hatPi$ is semi-perfect, and $\proj{\hatPi}$ is a Krull-Schmidt category. In this subsection, we identify the standard basis of $K_0\left(\proj{\Pi}\right)\cong \Z^n$ with $\{[\Pi e_i]\mid i\in Q_0\}$.
\begin{definition}
    Let $M \in \rep{\Pi}$ and $P_1(M) \rightarrow P_0(M) \rightarrow M \rightarrow 0$ be the minimal projective presentation of $M$. Then, we define the \textit{$g$-vector of $M$} by $$g(M) \coloneqq (g_1(M), \dots, g_n(M))^\trp = [P_0(M)]-[P_1(M)] \in K_0\left(\proj{\Pi}\right).$$
\end{definition}
A significance of $g$-vectors and $\tau$-tilting pairs is characterized by the following theorem:
\begin{prop}[{\cite[Theorem 5.1]{AIR}}]\label{prop:gvec-finite}
 Let $(M, P)$ be a support $\tau$-tilting pair for a finite dimensional $K$-algebra $\Lambda$ with $M=\bigoplus_{i=1}^l M_i$ and $P=\bigoplus_{i=l+1}^n P_i$ with $M_i$ and $P_i$ indecomposable. Then, $g(M_1), \dots, g(M_l)$, $g(P_{l+1}), \dots, g(P_n)$ form a basis of the Grothendieck group $K_0(\proj{\Lambda})$.
\end{prop}
Let $\affI{w} \in \Rep{\hatPi}$ denote the tilting ideal for $w \in \widehat{W}$ in Theorem \ref{bijection_W_tilt}. This module induces an automorphism on the bounded homotopy category $\mathcal{K}^b ({\proj{\hatPi}})$. In the work of Fu-Geng~\cite{FG}, they considered the action of affine Weyl group $\widehat{W}=W(\widehat{C})$ on the Grothendieck group of the bounded homotopy category $K_0 (\mathcal{K}^b(\proj{\hatPi}))\cong K_0 (\proj{\hatPi})\cong \mathbb{Z}^{n+1}$ by using tilting complexes.
    \begin{prop}[{\cite[Proof of Theorem 4.7]{FG}}] \label{prop:affgmat}
        Let $\widehat{C}=(c_{ij})$. We have a faithful group action $\sigma^* \colon \widehat{W} \rightarrow  \mathrm{GL}_{n+1}(\mathbb{Z}^{n+1})$ defined by
        $$s_i \mapsto [\affI{i} \Ltens{\hatPi}-] \curvearrowright K_0 (\proj{\hatPi})$$
        \[[\hatPi e_j] \mapsto \begin{cases}
        [\hatPi e_i]-\sum_{k \in \widehat{Q}_0} c_{kj}[\hatPi e_k]\,\,&(i=j)\\
        [\hatPi e_j]\,\,&(i\neq j).
        \end{cases}\]
    \end{prop}
    We introduce $g$-matrices for generalized preprojective algebras in view of above Proposition~\ref{prop:gvec-finite} and Proposition~\ref{prop:affgmat}.
\begin{definition}
    Let $\affI{w}\in \hatPi$ be the tilting ideal corresponding to $w \in \widehat{W}$ and let $\sigma$ be the Nakayama permutation of $\Pi$ from the self-injectivity of $\Pi$.
    \begin{enumerate}
    \item We define the $g$-matrix attached to $w \in \widehat{W}$ by $\hat{g}(w) = (g(\affI{w} e_0), \dots , g(\affI{w}e_n)) \in M_{n+1}(\Z)$.
    \item We define the $g$-matrix attached to $w \in W$ by $g(w) = (g_1, \dots , g_n) \in M_n(\Z)$, where $g_i \coloneqq g(I_w e_i) - (0, \dots, \delta_{I_w e_i, 0}, \dots , 0)$. (Here, $\delta_{I_w e_i, 0}$ is on the $\sigma(i)$-th coordinate.)
    \end{enumerate}
\end{definition}
Namely, $g(w)$ is the matrix whose column vectors correspond to $g$-vectors of indecomposable direct summands of $I_w$ and $P_w$ for the $\tau$-tilting pair $(I_w, P_w)$. The aim of this subsection is to compare $\hat{g}(w)$ with $g(w)$ for $w \in W \subset \widehat{W}$, and calculate $g(w)$. Recall that $\pdim{\hatPi}{\affI{w}} \leq 1$.
\begin{theorem}\label{resl_Iw}
Let $w \in W$ and $0 \rightarrow P_1 \xrightarrow{f_1} P_0 \xrightarrow{f_0} \affI{w}e_i \rightarrow 0$ be the minimal projective resolution for $\affI{w}e_i$. Applying 
$\overline{(-)} \coloneqq \Pi \tens{\hatPi}-$
to this sequence, we obtain an exact sequence 
$$0 \rightarrow \nu^{-1}(\Pi e_i/ I_w e_i) \rightarrow \overline{P_1} \xrightarrow{\overline{f_1}} \overline{P_0} \xrightarrow{\overline{f_0}} I_w e_i \rightarrow 0,$$
where $\nu$ is the Nakayama functor for $\rep{\Pi}$, and the sequence $$\overline{P_1} \xrightarrow{\overline{f_1}} \overline{P_0} \xrightarrow{\overline{f_0}} I_w e_i \rightarrow 0$$ is the minimal projective presentation of $I_w e_i\in \rep{\Pi}$ if $I_w e_i \neq 0$.
\begin{proof}
First, we remind that $\Pi \cong \hatPi /\langle e_0 \rangle$ and $\Pi \tens{\hatPi}{\affI{w}e_i} \cong (\affI{w}e_i)/ (\langle e_0 \rangle \affI{w}e_i) \cong (\affI{w}e_i)/(\langle e_0 \rangle e_i) \cong I_w e_i$, where we use $\affI{w}\supset \langle e_0 \rangle$. We apply the tensor functor $\Pi\tens{\hatPi}{(-)}$ to the minimal projective resolution of $\affI{w}e_i$. Then, we obtain the following exact sequence in $\rep{\hatPi}$:
$$0 \rightarrow \Tor{1}{\hatPi}{\Pi}{I_w e_i}  \rightarrow \overline{P_1} \xrightarrow{\overline{f_1}} \overline{P_0} \xrightarrow{\overline{f_0}} I_w e_i \rightarrow 0.$$
Since $\overline{f}_0$ is minimal by construction, we have enough to show that  
$\Tor{1}{\hatPi}{\Pi}{I_w e_i}$ is not projective. Since any $\hatPi$-module with finite projective dimension is locally free  and $\pdim{\hatPi}{\affI{w}}\leq 1$ by Proposition \ref{Pi_projdim} and Theorem \ref{bijection_W_tilt}, we have the following isomorphisms:
\begin{align*}
\Tor{1}{\hatPi}{\Pi}{I_w e_i}
&\cong \D\Ext{1}{\hatPi}{I_w e_i}{\D\Pi}&&\textit{(\cite[Appendix  Proposition 4.11]{ASS})}\\ &=\D\Ext{2}{\hatPi}{\hatPi e_i/\affI{w}e_i}{\D\Pi}\\
&=\Hom{\hatPi}{\D\Pi}{\hatPi e_i /\affI{w}e_i}&&\textit{(Proposition \ref{PiExtduality} and locally-freeness)}\\
&=\Hom{\Pi}{\D\Pi}{\Pi e_i / I_w e_i}\\
&=\nu^{-1}(\Pi e_i/ I_w e_i).&&\textit{(Definition of Nakayama functor)}
\end{align*}
Note that we have an exact sequence 
$$0\rightarrow I_w e_i \rightarrow \Pi e_i \rightarrow \Pi e_i/ I_w e_i \rightarrow 0,$$
and we apply the functor $\Hom{\affPi}{-}{\D\Pi}$ to obtain
\begin{align}
\Ext{1}{\hatPi}{\hatPi e_i}{\D \Pi} &\rightarrow \Ext{1}{\hatPi}{\affI{w} e_i}{\D \Pi}\rightarrow \Ext{2}{\hatPi}{\hatPi e_i/ \affI{w} e_i}{\D \Pi} \rightarrow \Ext{2}{\hatPi}{\hatPi e_i}{\D \Pi}.\label{eq:seq4}
\end{align}
Since both of two extremal terms in the sequence \eqref{eq:seq4} are equal to $0$, we have $\Ext{1}{\hatPi}{\affI{w} e_i}{\D \Pi} \cong \Ext{2}{\hatPi}{\hatPi e_i/ \affI{w} e_i}{\D \Pi}$.
Since $\Pi$ is a self-injective algebra and $I_w e_i\neq 0$, it follows that $\Pi e_i/ I_w e_i$ is not a projective indecomposable object. In particular, $\nu^{-1}(\Pi e_i/ I_w e_i)$ is not projective. So, we obtain our assertion.
\end{proof}
\end{theorem}
We obtain the following Corollary~\ref{cor:gmatforget} by the similar argument as in \cite[Proposition 3.4]{Miz}.
\begin{cor}\label{cor:gmatforget}
Let $w\in W$. The matrix $g(w)$ is obtained from $\hat{g}(w)$ by forgetting the $0$-th row and the $0$-th column. 
\begin{proof}
We take the minimal projective presentation of $\affI{w}e_i \in \Rep{\hatPi}$:
$$0 \xrightarrow{f_1} P_1 \xrightarrow{f_0} P_0 \rightarrow \affI{w}e_i \rightarrow 0.$$
We decompose $P_0 = P'_0 \oplus (\hatPi e_0)^{\oplus l}$ and $P_1=P'_1 \oplus (\hatPi e_0)^{\oplus m}$, where $P'_0$ and $P'_1$ do not belong to $\add{\hatPi e_0}$. For the $i$-th column  $\hat{g}^i(w)=\left(g_0 (\affI{w} e_i), \dots g_n (\affI{w} e_i )\right)^\trp$ of the $g$-matrix $\hat{g}(w) = \left(g(\affI{w} e_0), \dots , g(\affI{w}e_n)\right)$, we find that $\left(g_1(\affI{w} e_i), \dots , g_n (\affI{w}e_i)\right)$ depends only upon $P'_0$ and $P'_1$ by the definition of $g$-vectors. So, we have enough to show that it coincides with the $i$-th column of the $g$-matrix $g(w) = (g(I_w e_1), \dots , g(I_w e_n ))$. Now, we have the following exact sequence:
$$0 \rightarrow \nu^{-1}(\Pi e_i/ I_w e_i) \rightarrow \overline{P_1} \xrightarrow{\overline{f_1}} \overline{P_0} \xrightarrow{\overline{f_0}} I_w e_i \rightarrow 0.$$
If we have the equality $P'_0 \neq 0$, then our assertion follows from Theorem \ref{resl_Iw}. On the other hand, if we have $P'_0 = 0$, then it holds that $\overline{P'_1} \cong \nu^{-1}(\Pi e_i) \cong \Pi e_{\sigma(i)}$, and so we have $P'_1 \cong \hatPi e_{\sigma(i)}$. So, we obtain our assertion by the definition of $g(w)$.
\end{proof}
\end{cor}
Under these preparations, we can realize Tits cones by using $g$-matrices. 
    \begin{theorem}\label{thm:gmatNak}
    We have an equality $g(w)= \sigma^*(w)$ for each $w \in W$.
    \begin{proof}
    We can obtain our assertion by combining our propositions and the fact~\cite[Proposition 3.5]{FG} that the minimal projective resolution of $\affI{i} e_i$ is 
    $$0 \rightarrow \hatPi e_i \rightarrow \bigoplus_{j\in \overline{\Omega}(i, -)}(\hatPi e_j)^{|c_{ji}|} \rightarrow \affI{i} e_i \rightarrow 0.$$ Note that $\affI{i} e_j = \hatPi e_j$ if $i\neq j$ by Corollary \ref{cor:gmatforget}. Thus, we can prove our assertion by induction on length of $w \in W$.
    \end{proof}
    \end{theorem}
    Note that the Nakayama permutation for $\Pi$ is characterized by
    $$w_0(\varpi_i) = -\varpi_{\sigma(i)}$$
    as a direct consequence of Theorem \ref{thm:gmatNak} \footnote{The author was motivated by some questions about this from Y. Mizuno and S. Kato to use $g$-vectors in study of MV polytopes. B. Leclerc was kind to point out the author's error about a description of the Nakayama permutation. He thanks for their comments.}. In particular, the Nakayama permutation for the generalized preprojective algebra coincides with the trivial permutation if we consider one of cases of $\mathsf{B, C, D_{2n}, E_7, E_8, F_4, G_2}$. Thus, we realize Weyl chambers in terms of $g$-matrices. In particular, we obtain the following by the similar argument as in \cite[Theorem 3.9]{Miz}.
    \begin{theorem}\label{Tits}
        We define the cone $C(w)$ associated with $w \in W$ by
        $\{ a_1 g_1(w)+\cdots + a_n g_n(w) \mid a_i \in \mathbb{R}_{>0}\}$. Then, this cone is identified with a Weyl chamber in the \S \ref{pre:rootsys} by our choice of basis of $\Z^n$.
       \begin{proof}
        We obtain the assertion by computing
       \begin{align*}
        C(w)&=\{ a_1 g_1(w)+\cdots + a_n g_n(w) \mid a_i \in \mathbb{R}_{>0}\}\\
        & =\{ a_1 \sigma^*(w)(\varpi_1)+\cdots + a_n \sigma^*(w)(\varpi_n) \mid a_i \in \mathbb{R}_{>0}\}\\
        &= \sigma^*(w)(C_0).
        \end{align*}
       Thus, we obtain the realization.
       \end{proof}
    \end{theorem}
    \begin{example}\begin{enumerate}
        \item[(1)] Let $\Pi$ be the generalized preprojective algebra of type $\mathsf{A}_3$ with any symmetrizer. $g$-matrices for $\Pi$ are given by the lattice corresponding to mutations of support $\tau$-tilting modules as Figure \ref{fig:A3}. (The same example for the classical case has already appeared in the work of Mizuno~\cite[Example 3.10 (b)]{Miz}.):
        
    \item[(2)] Let $\Pi$ be the generalized preprojective algebra of type $\mathsf{B}_2$ with any symmetrizer. $g$-matrices for $\Pi$ are given by the lattice corresponding to mutations of support $\tau$-tilting modules as Figure \ref{fig:B2}.
    \end{enumerate}
    \end{example}
    \begin{rem}
	    By definition, any $V_{\mathbf{i}, k}$ is an indecomposable $\tau$-rigid module. For $\mathsf{A, D, E}$ and minimal symmetrizer cases, the set of isoclasses of layer modules coincides with the set of isoclasses of bricks. In general, there exists a set $\{B_1, \dots, B_n\}$ of isoclasses of bricks whose dimension vectors coincide with the set of column vectors $\{c_1(w), \dots, c_n(w)\}$ of a $c$-matrix $c(w)\coloneqq (g(w)^{-1})^\trp$ defined by Fu~\cite{Fu} up to sign (\cf Treffinger~\cite{Tre}). Note that these column vectors of the $c$-matrix are naturally identified with some elements of $\Delta(C^\trp)$ by Theorem \ref{thm:gmatNak}. Let $\mathcal{B}_w$ be the set of isoclasses of bricks such that $\Hom{\Pi}{B}{B'}=0$ for any non-isomorphic $B, B'\in \mathcal{B}_w$. Then, the smallest torsion class which contains $\mathcal{B}_w$ coincides with $\Fac{I_w}$. In addition, $\mathcal{B}_w$ is given by $\ind (I_w/ \rad_{\Gamma} I_w)$, where let $\Gamma \coloneqq {\End{\Pi}{I_w}}^{\mathsf{op}}$ and let $\ind (I_w/ \rad_{\Gamma} I_w)$ denote the set of isoclasses of indecomposable summands of $I_w/ \rad_{\Gamma} I_w$ (\cf Asai~\cite{Asa2}). Then, $\mathcal{B}_w$ bijectively corresponds to the set $\{c_i(w) \in \Delta(C^\trp) \mid i\in Q_0 \}\cap \Delta^+(C^\trp)$ by taking dimension vectors. A similar kind of correspondence between bricks and positive Schur roots for Langlands dual datum has been considered in the work of Gei\ss-Leclerc-Schr\"oer~\cite{GLS6} for $H(C, D, \Omega)$ (see Definition \ref{def of H&Pi}). So, the combinatorial structures about torsion-free classes are expected to be described in terms of the root system of $\mathfrak{g}( C^\trp)$ as in Asai~\cite{Asa} or Enomoto~\cite{Eno}. 
	\end{rem}
        \afterpage{\clearpage}
\newpage
\section{Nilpotent varieties and MV polytopes}\label{section-MV}
Keep the setting of the previous sections.
\subsection{On crystal structures on irreducible components and PBW parametrizations}

In this subsection, we realize Lusztig's PBW parametrizations via nilpotent varieties of Gei\ss-Leclerc-Schr\"oer. First, we recall definitions and basic properties about crystals.
\begin{definition}
    Let $C$ be a symmetrizable GCM of finite type and let $P\coloneqq \bigoplus_{i\in Q_0}\Z \varpi_i$ be the associated integral weight lattice. A crystal is a tuple $(B, \wt{-}, \tilde{e}_i, \tilde{f}_i, \varepsilon_i, \varphi_i)$ of a set $B$, maps $\wt{-}\colon B\rightarrow P$, $\tilde{e}_i, \tilde{f}_i\colon B\rightarrow B\cup \{\varnothing \}$ and $\varepsilon_i, \varphi_i\colon B\rightarrow \Z$ which satisfy the following conditions for any $i\in Q_0$ and for any $b\in B$:
    \begin{enumerate}
        \item[\textsf{(cr1)}] $\varphi_i(b)=\varepsilon_i(b)+ \left\langle \wt{b}, \alpha_i \right\rangle$;
        \item[\textsf{(cr2)}]
        $\varphi_i(\tilde{e}_i(b))=\varphi_i(b)+1$, $\varepsilon_i(\tilde{e}_i)(b)=\varepsilon_i(b)-1$, $\wt{\tilde{e}_i(b)}=\wt{b}+\alpha_i$\,\,\,($\tilde{e}_i(b)\neq \varnothing$);
        \item[\textsf{(cr3)}] Two equalities $\tilde{f}_i(b)=b'$ and $\tilde{e}_i(b')=b$ are equivalent.
    \end{enumerate}
    A lowest weight crystal is a crystal with an element $b_-\in B$ such that
    \begin{enumerate}
    \item[\textsf{(cr4)}]
    there exists a sequence $(i_1, \dots, i_t)$ with $i_k\in Q_0$ for any $1\leq k \leq t$ such that $b_-=\tilde{f}_{i_1}\cdots \tilde{f}_{i_t}(b)$ for each $b\in B$;
    \item[\textsf{(cr5)}]We have $\varphi_i(b)=\max \{m\mid \tilde{f}_i^m(b)\neq \varnothing \}$.
    \end{enumerate}
\end{definition}
For a lowest weight crystal $B$, the functions $\wt{-}$, $\tilde{f}_i$, $\varepsilon_i$ and $\varphi_i$ are determined by the $\tilde{e}_i$ and $\wt{b_-}$. We mainly deal with infinity crystals $B(-\infty)$. Tingley-Webster~\cite[Proposition 1.4]{TW} gave a criterion whether a lowest weight crystal is isomorphic to the infinity crystal $B(-\infty)$ as a reformulation of the well-known criterion by Kashiwara-Saito~\cite[Proposition 3.2.3]{KS}. We think of the Tingley-Webster criterion as a definition of $B(-\infty)$ as in the work of \cite{GLS4}.
\begin{def-pro}[{\cite[Proposition 3.2.3]{TW}}]
Fix a set $B$ with operators $\tilde{e}_i, \tilde{f}_i, \tilde{e}_i^*, \tilde{f}_i^*\colon B\rightarrow B\cup \{\varnothing \}$, and assume $(B, \tilde{e}_i, \tilde{f}_i)$ and $(B, \tilde{e}_i^*, \tilde{f}_i^*)$ are both lowest weight crystals with the same lowest weight element $b_-$. Now, we assume that the other data are determined by setting $\wt{b_-}=0$. We have $(B, \tilde{e}_i, \tilde{f}_i )\cong (B, \tilde{e}_i^*, \tilde{f}_i^* )\cong B(-\infty)$ if they satisfy the following conditions for each $i, j \in Q_0$ and any $b\in B$:
\begin{itemize}
    \item $\tilde{e}_i(b), \tilde{e}_i^*(b) \neq \varnothing$;
    \item If $i\neq j$, then we have $\tilde{e}_i^*\tilde{e}_j(b)=\tilde{e}_j \tilde{e}_i^*(b)$;
    \item We have $\varphi_i(b)+ \varphi_i^*(b)-\langle \wt{b}, \alpha_i \rangle\geq 0$;
    \item If we have $\varphi_i(b)+ \varphi_i^*(b)-\langle \wt{b}, \alpha_i \rangle=0$, then $\tilde{e}_i(b)=\tilde{e}_i^*(b)$;
    \item If we have $\varphi_i(b)+ \varphi_i^*(b)-\langle \wt{b}, \alpha_i \rangle\geq 1$, then $\varphi_i(\tilde{e}_i^*(b) )=\varphi_i(b)$ and $\varphi_i^* \left(\tilde{e}_i(b)\right)=\varphi_i^*(b)$;
    \item If we have $\varphi_i(b)+ \varphi_i^*(b)-\langle \wt{b}, \alpha_i \rangle\geq 2$, then $\tilde{e}_i^*\tilde{e}_i(b)=\tilde{e}_i \tilde{e}_i^*(b)$.
\end{itemize}
\end{def-pro}
Gei\ss-Leclerc-Schr\"oer~\cite{GLS4} generalized a geometric realization of infinity crystals $B(-\infty)$ by Kashiwara-Saito~\cite{KS}. They define a special class of $\E$-filtered modules, called crystal modules, and they construct a crystal structure on the set of special irreducible components of varieties of $\E$-filtered modules with specific rank vectors.
\begin{definition}
    Let $M \in \Efilt{\Pi}$, and $0 \rightarrow K_i(M) \rightarrow M \rightarrow \fac_i(M) \rightarrow 0$ and $0 \rightarrow \sub_i(M) \rightarrow M \rightarrow C_i(M) \rightarrow 0$ be canonical short exact sequences. We say that $M$ is a crystal module if and only if:
    \begin{enumerate}
    \item $\{0\}$ is a crystal module.
    \item $K_i(M)$ and $C_i(M)$ are crystal modules for any $i \in Q_0$, and $\fac_i{M}$ and $\sub_i{M}$ are locally free modules for any $i\in Q_0$.
    \end{enumerate}
\end{definition}
Let $\mathbf{i}\coloneqq (i_1, \dots , i_l)$ be a reduced expression of the longest element of $W$. Note that any series of inclusions of torsion classes, which has length $\length{w_0}$ (called \textit{maximal green sequence}) corresponds to a reduced expression of the longest element of $W$ by the $\tau$-tilting theory. Namely, we can associate the maximal green sequence $\rep{\Pi}= \Fac{I_e} \supsetneq \Fac{I_{i_1}} \supsetneq \cdots \supsetneq \Fac{I_{s_{i_1}\cdots s_{i_l}}} = \{0\}$ to $\mathbf{i}$. We characterize objects in $\Fac(I_w)$ and $\Sub{\Pi/I_w}$.
\begin{theorem}\label{thm:HNseq}
Let $M\in \rep{\Pi}$ be a crystal module and $w \in W$ such that $\length{ws_i}= \length{w} + 1$. Then, $T_w/T_{ws_i} \cong \left(I_w/I_{ws_i}\right)^{\oplus a}$ for some $a \in \mathbb{Z}_{\geq 0}$, where $T_w$ is the torsion submodule of $M$ associated to the torsion class $\Fac{I_{w}}$ for $w \in W$.
\begin{proof}
We put $X = \Hom{\Pi}{I_w}{M}$. Then, we have the canonical exact sequence from Theorem \ref{refl}:
$$0\rightarrow \Sigma_i^{-}\Sigma_i X \rightarrow X \rightarrow\fac_i X\rightarrow 0.$$
 By applying $I_w\tens{\Pi}(-)$, we obtain an exact sequence:
$$0\rightarrow I_{ws_i}\tens{\Pi}{\Hom{\Pi}{I_{ws_i}}{M}}\rightarrow I_w\tens{\Pi}{\Hom{\Pi}{I_w}{M}}\rightarrow I_w\tens{\Pi}{{E_i}^{\oplus m}}\rightarrow 0,$$
where $c_i m=\dim{\fac_i(X)}$.
Here, we use $\Tor{\Pi}{1}{I_w}{E_i} = 0$ for $\length{ws_i}>\length{w}$ by Proposition \ref{prpty-I_w}. 
So, we obtain our assertion from Theorem \ref{VM_iso}.
\end{proof}
\end{theorem}

\begin{cor}\label{HNcor}
Let $T_0\coloneqq M \in \rep{\Pi}$ and $T_j$ be the torsion submodule of $T_{j-1}$ associated to the torsion class $\Fac{I_{s_{i_1}\cdots s_{i_j}}}$ ($1\leq j \leq l$). Then, we have a filtration $M=T_0 \supsetneq T_1 \supsetneq \cdots \supsetneq T_l =0$ such that $T_{j-1}/T_{j} \cong (I_{s_{i_1}\cdots s_{i_{j-1}}}/I_{s_{i_1}\cdots s_{i_j}})^{\oplus a_j}$ for some $a_j \in \mathbb{Z}_{\geq 0}$.
\end{cor}
We refer to the filtration $M=T_0 \supsetneq T_1 \supsetneq \cdots \supsetneq T_l =0$ in Corollary \ref{HNcor} as a \textit{Harder-Narasimhan filtration} for $\mathbf{i}=(i_1, \dots, i_l)$ and refer to the modules $I_{s_{i_1}\cdots s_{i_{j-1}}}/I_{s_{i_1}\cdots s_{i_j}}\,\,(1\leq j \leq l)$ as \textit{layer modules}. The aim of this section is to prove that it gives rise to a stratification of irreducible components of nilpotent varieties from these data.
\begin{rem}
For non-Dynkin cases, idempotent ideals corresponding to elements of affine Weyl group never coincide with $\{0\}$. However, the similar short exact sequences as in the Proof of Theorem \ref{thm:HNseq} can be obtained also for non-Dynkin cases, though generalized preprojective algebras of non-Dynkin type are always infinite dimensional. In addition, we have Harder-Narasimhan filtrations determined by stabilities (\cf \cite[\S 3.2]{BKT}) also for non-Dynkin cases, though the $\tau$-tilting theory does not work.
\end{rem}
\begin{example}\label{ex:layerB2}
    Let $C=\left(
	\begin{array}{rr}
		2 & -1 \\
	   -2 & 2 \\
	\end{array}
	\right)
	$ and
	$D=\diag{2, 1}$. In this $\mathsf{B}_2$ case, the lattice of mutations of $\sttilt{\Pi}$ is illustrated as Figure \ref{fig:sttiltB2}.
	\begin{figure}
	    \centering
	\begin{tikzpicture}
    \node (1) at (0, 0) {$\begin{tikzpicture}[auto,baseline=-3pt]
\node (25) at (0, -0.45) {$1$};
\node (26) at (0.3, -0.15) {$2$};  
\node (27) at (-0.3, -0.15) {$1$};
\node (28) at (-0.3, 0.15) {$2$};
\node (29) at (0.3, 0.15) {$1$};
\node (30) at (0, 0.45) {$1$};
\end{tikzpicture} 
	\oplus
	\begin{tikzpicture}[auto,baseline=-3pt]
\node (31) at (0, -0.45) {$2$};
\node (32) at (0, -0.15) {$1$};
\node (33) at (0, 0.15) {$1$};
\node (34) at (0, 0.45) {$2$};
\end{tikzpicture}$};
	\node (2) at (-2, -2) {$\begin{tikzpicture}[auto,baseline=-3pt]
\node (35) at (0, -0.45) {$1$};
\node (36) at (0.3, -0.15) {$2$};  
\node (37) at (-0.3, -0.15) {$1$};
\node (38) at (-0.3, 0.15) {$2$};
\end{tikzpicture} 
	\oplus
	\begin{tikzpicture}[auto,baseline=-3pt]
\node (39) at (0, -0.45) {$2$};
\node (40) at (0, -0.15) {$1$};
\node (50) at (0, 0.15) {$1$};
\node (51) at (0, 0.45) {$2$};
\end{tikzpicture}$};
	\node (3) at (2, -2) {$\begin{tikzpicture}[auto,baseline=-3pt]
\node (52) at (0, -0.45) {$1$};
\node (53) at (0.3, -0.15) {$2$};  
\node (54) at (-0.3, -0.15) {$1$};
\node (55) at (-0.3, 0.15) {$2$};
\node (56) at (0.3, 0.15) {$1$};
\node (57) at (0, 0.45) {$1$};
\end{tikzpicture} 
	\oplus
	\begin{tikzpicture}[auto,baseline=-3pt]
\node (58) at (0, -0.45) {$2$};
\node (59) at (0, -0.15) {$1$};
\node (60) at (0, 0.15) {$1$};
\end{tikzpicture}$};
	\node (4) at (-2, -4) {$\begin{tikzpicture}[auto,baseline=-3pt]
\node (61) at (0, -0.45) {$1$};
\node (62) at (0.3, -0.15) {$2$};  
\node (63) at (-0.3, -0.15) {$1$};
\node (64) at (-0.3, 0.15) {$2$};
\end{tikzpicture} 
	\oplus
	\begin{tikzpicture}[auto,baseline=-3pt]
\node (65) at (0, 0) {$2$};
\end{tikzpicture}$};
	\node (5) at (2, -4) {$\begin{tikzpicture}[auto,baseline=-3pt]
\node (66) at (0, -0.15) {$1$};
\node (67) at (0, 0.15) {$1$};
\end{tikzpicture} 
	\oplus
	\begin{tikzpicture}[auto,baseline=-3pt]
\node (68) at (0, -0.45) {$2$};
\node (69) at (0, -0.15) {$1$};
\node (70) at (0, 0.15) {$1$};
\end{tikzpicture}
	$};
	\node (6) at (-2, -5.5) {$\begin{tikzpicture}[auto]
	\node (71) at (0, 0) {$2$};
	\end{tikzpicture}
	$};
	\node (7) at (2, -5.5) {$
	\begin{tikzpicture}[auto]
	\node (72) at (0, 0) {$1$};
	\end{tikzpicture}$};
	\node (8) at (0, -6.5) {$\{0\}
	$};
	\draw[->] (1) to (2);
	\draw[->] (1) to (3);
	\draw[->] (2) to (4);
	\draw[->] (3) to (5);
	\draw[->] (4) to (6);
	\draw[->] (5) to (7);
	\draw[->] (6) to (8);
	\draw[->] (7) to (8);
    \end{tikzpicture}
    \caption{$\sttilt{\Pi}$ of type $\mathsf{B}_2$}
	    \label{fig:sttiltB2}
	\end{figure}
	All layer modules are illustrated as the following Loewy series:
     \begin{alignat*}{4}
        \Pi/I_{s_1}&=
        \begin{tikzpicture}[auto,baseline=-3pt]
\node (x) at (0, 0.15) {$1$};
\node (a) at (0, -0.15) {$1$};
\end{tikzpicture};
&I_{s_1}/I_{s_1 s_2}&=
\begin{tikzpicture}[auto,baseline=-3pt]
\node (x) at (0, 0.30) {$2$};
\node (a) at (0, 0) {$1$}; 
\node (b) at (0, -0.30) {$1$};  
\end{tikzpicture};
&I_{s_1 s_2}/I_{s_1 s_2 s_1}&= \begin{tikzpicture}[auto,baseline=-3pt]
\node (x) at (0, 0) {};
\node (a) at (0, -0.45) {$1$}; 
\node (b) at (0.3, -0.15) {$2$};  
\node (c) at (-0.3, -0.15) {$1$};
\node (e) at (-0.3, 0.15) {$2$};
\end{tikzpicture};
&I_{s_1 s_2 s_1}/I_{s_1 s_2 s_1 s_2}&=\begin{tikzpicture}[auto,baseline=-3pt]
\node (x) at (0, 0) {2};
\end{tikzpicture};\\
\Pi/I_{s_2}&=\begin{tikzpicture}[auto,baseline=-3pt]
\node (x) at (0, 0) {$2$};
\end{tikzpicture};
&I_{s_2}/I_{s_2 s_1}&= \begin{tikzpicture}[auto,baseline=-3pt]
\node (x) at (0, 0) {};
\node (b) at (0.3, -0.15) {$2$};  
\node (d) at (0.3, 0.15) {$1$};
\node (e) at (-0.3, 0.15) {$2$};
\node (f) at (0, 0.45) {$1$};
\end{tikzpicture};
&I_{s_2 s_1}/I_{s_2 s_1 s_2}&=
\begin{tikzpicture}[auto,baseline=-3pt]
\node (b) at (0, -0.3) {$2$};  
\node (x) at (0, 0) {$1$};
\node (d) at (0, 0.3) {$1$};
\end{tikzpicture};
&I_{s_2 s_1 s_2}/I_{s_2 s_1 s_2 s_1}&=
\begin{tikzpicture}[auto,baseline=-3pt]
\node (b) at (0, -0.15) {$1$};  
\node (d) at (0, 0.15) {$1$};
\end{tikzpicture}.
\end{alignat*}
\end{example}
\begin{definition}
     A module $M\in Z$ for $Z \in \Irr{\Pi(\mathbf{r})}$ is called \textit{generic} if and only if $M$ is contained in a sufficiently small open dense subset of $Z$ defined by a finite set of suitable open conditions which are suitable for contexts. In particular, all the upper semi-continuous constructible functions on $Z$ defined below take its minimal value at $M$. For example, we can take an open dense subset $U_Z$ for each $Z\in \Irr{\Pi(\mathbf{r})}$ such that $\sub_{k}(M) \cong \sub_{k}(N)$ and $\fac_{k}(M) \cong \fac_{k}(N)$ for any $M, N \in U_Z$ by upper semi-continuity. So, we can define $\fac_{k}(Z)\coloneqq \fac_{k}(M)$ and $\sub_{k}(Z)\coloneqq \sub_{k}(M)$ up to isomorphism. We demand $M\in Z$ with $\sub_{k}(M) \cong \sub_{k}(Z)$ and $\fac_{k}(M) \cong \fac_{k}(Z)$ in this context.
     We say that $Z \in \Irr{\Pi(\mathbf{r})}$ is crystal if $Z$ contains an open dense subset of crystal modules. 
\end{definition}

\begin{theorem}[{\cite[Theorem 4.1]{GLS4}}]\label{thm:maxcomp}
    For each $Z\in \Irr{\Pi(\br)}$, we have
    $$\dim Z \leq \dim H(\br).$$
\end{theorem}

We use the following criterion given by \cite{GLS4}.
\begin{prop}[{\cite[Proposition 4.4]{GLS4}}]\label{prop:dimofcrycom}
 For $Z\in \Irr{\Pi(\br)}$, the irreducible component $Z$ is crystal if and only if $\dim Z = \dim H(\br) = \dim G(\br) - q_{DC}(\br)$.
\end{prop}
Let $\Irrmax{\Pi(\mathbf{r})}$ denote the set of crystal irreducible components of $\Pi(\mathbf{r})$.
Gei\ss-Leclerc-Schr\"oer~\cite{GLS4} realized $B(-\infty)$ for general symmetrizable GCMs and for general symmetrizers as follows:
\begin{def-thm}[\cite{GLS4}]\label{defthm:crystal}
We give the following:
\begin{enumerate}
    \item Let $\mathcal{B} \coloneqq \coprod_{\mathbf{r}\in {\mathbb{N}}^n} \Irrmax{\Pi(\mathbf{r})}$. We write \begin{align*}
    &{\Pi(\br)}^{k, p}\coloneqq \{M \in \Pi(\br) \mid \fac_k(M) \cong E_k^p \}; \\ &{\Pi(\br)}_{k, p}\coloneqq \{M \in \Pi(\br) \mid \sub_k(M) \cong E_k^p \}.
    \end{align*}
    Let $\varphi_i, \varphi_i^* \colon \Efilt{\Pi} \rightarrow \mathbb{Z}$ denote functions 
    \begin{align*}
    \varphi_i \colon M \mapsto [\sub_i(M)\colon E_i],
    &&\varphi_i^* \colon M\mapsto [\fac_i(M)\colon E_i],
    \end{align*}
    where $[M: E_i]$ denotes the maximal number $p \geq 0$ such that there exists a direct sum decomposition $M= U \oplus V$ with $U \cong E_i^p$ for $M \in \rep{H_i}$.
    \item Let $\mathrm{wt}(Z) \coloneqq \br$ for $Z \in \Irrmax{\br}$. We write 
    \begin{align*}
    \varphi_i(Z) \coloneqq \min \{ \varphi_i(M) \mid M \in Z \},
    &&\varepsilon_i(Z) \coloneqq \varphi_i(Z) - \langle \wt{Z}, \alpha_i \rangle.
    \end{align*}
     Dually, we write
     \begin{align*}
     \varphi_i^*(Z) \coloneqq \min\{ \varphi_i^*(M) \mid M \in Z \},
     &&\varepsilon_i^*(Z) \coloneqq \varphi_i^*(Z) - \langle \wt{Z}, \alpha_i \rangle.
     \end{align*}
    \item Let $r \coloneqq p-q \geq 1$. Then, we have a bijection 
    $$f_{k, r} \colon\djuni{\Irrmax{\Pi(\br)_{k, p}}} \rightarrow \djuni{\Irrmax{\Pi(\br -r \alpha_k)_{k, q}}}.$$
    Let $e_{k, r}\coloneqq (f_{k, r})^{-1}$ denote the inverse map of $f_{k, r}$. Dually, we have a bijection $$f_{k, r}^* \colon\djuni{\Irrmax{\Pi(\br)^{k, p}}} \rightarrow \djuni{\Irrmax{\Pi(\br -r \alpha_k)^{k, q}}}.$$
    Let $e_{k, r}^*\coloneqq (f_{k, r}^*)^{-1}$ denote the inverse map of $f_{k, r}^*$. 
    
    \item Under above notation, we have a crystal isomorphism 
    $$(\mathcal{B}, \wt{-}, \tilde{e}_i, \tilde{f}_i, \varepsilon_i, \varphi_i ) \cong (\mathcal{B}, \wt{-}, \tilde{e}_i^*, \tilde{f}_i^*, \varepsilon_i^*, \varphi_i^* ) \cong B(-\infty)$$
    by assigning the unique irreducible component $Z_{-}$ of $\Pi(0)$ to the lowest weight element $b_-$,
    where we write 
    \begin{align*}
    \tilde{f}_i(Z) \coloneqq \begin{cases} \overline{f_{i, 1}(Z)}\,\,&(\varphi_i(Z)\geq 1)\\
    \varnothing\,\,&(\varphi_i(Z)=0)
    \end{cases}
     && \tilde{e}_i(Z) \coloneqq \overline{e_{i, 1}(Z)},
     \end{align*}
     and 
     \begin{align*}
     \tilde{f}_i^*(Z) \coloneqq \begin{cases}
     \overline{f_{i, 1}^*(Z)}\,\,&(\varphi_i^*(Z)\geq 1)\\
    \varnothing\,\,&(\varphi_i^*(Z)=0)
    \end{cases}
     && \tilde{e}_i^*(Z) \coloneqq \overline{e_{i, 1}^*(Z)}.
     \end{align*}
\end{enumerate}
\end{def-thm}
In this setting, we define the Kashiwara $*$-operator as follows:
\begin{definition}[{\cite[\S 5.4]{GLS4}}]
    Let $M=(M(\varepsilon_i), M(\alpha_{ij}))\in \nil_{\E}(\Pi, d)$ and let $M^* \coloneqq (M(\varepsilon_i)^*, M(\alpha_{ij})^*)\in \nil_{\E}(\Pi, d)$, where $M(\varepsilon_i)^* \coloneqq M(\varepsilon_i)^{\trp}$ and $M(\alpha_{ij})^* \coloneqq M(\alpha_{ji})^\trp$. For each dimension vector $d$, we obtain
    an automorphism $*_d$ of the variety $\nil_{\E}(\Pi, d)$ defined as $*_d(M)\coloneqq M^*$. This yields an automorphism $*_{\br}$ of $\Pi(\br)$. This $*_{\br}$ yields an involution $*_\br \colon \calB_\br \rightarrow \calB_\br$. This induces an involution $*\colon \calB\rightarrow \calB$.
    In this setting, we have
    \begin{align*}
        &*\tilde{e}_i*=\tilde{e}_i^* 
        && *\tilde{e}^*_i*=\tilde{e}_i,
        &*\tilde{f}_i*=\tilde{f}_i^*,
        &&*\tilde{f}^*_i*=\tilde{f}_i.
    \end{align*}
\end{definition}

The aim of this subsection is to prove the following theorem by the help of the argument in the work of Baumann-Kamnitzer-Tingley~\cite[\S 4]{BKT}. This kind of theorem for a classical situation has appeared in the work of Geiss-Leclerc-Schr\"oer~\cite[\S 14.2]{GLS7}.
\begin{theorem}\label{strata-thm}
We choose a reduced expression $\mathbf{i}=(i_1, \dots , i_r)$ of $w\in W$. Let $\Pi_{\mathbf{i}}^{\mathbf{a}}$ be the set of crystal modules which have a Harder-Narasimhan filtration $M= T_0 \supsetneq T_1 \supsetneq \cdots \supsetneq T_r =0$ such that $T_{j-1}/ T_{j} \cong (I_{s_{i_1}\cdots s_{i_{j-1}}}/I_{s_{i_1}\cdots s_{i_j}})^{\oplus a_j}$ for  $\mathbf{a} \coloneqq (a_1, \dots, a_r) \in \Z^r$. Then, $\Pi_{\mathbf{i}}^{\mathbf{a}}$ is an irreducible constructible set and  $Z_{\mathbf{i}}^{\mathbf{a}}\coloneqq\overline{\Pi_{\mathbf{i}}^{\mathbf{a}}}$ is a maximal irreducible component.
\end{theorem}

We give a proof of this theorem after Proposition \ref{strata-1step}. We define $S\coloneqq \prod_{i\in Q_0}H_i$. We fix a torsion pair $(\mathcal{T}, \mathcal{F})$ in $\rep{\Pi}$. We consider the set $\T(\br)$ (\resp $\F(\br)$) of all maximal irreducible components of $\Pi(\br)$ whose generic points belong to $\calT$ (\resp $\calF$). We define $\T \coloneqq \coprod_{\br\in \Z^n_{\geq 0}} \T(\br)$ (\resp $\F \coloneqq \coprod_{\br\in \Z^n_{\geq 0}} \F(\br)$). Let $(\br_t, \br_f)$ be a tuple of $\N^n$, and let $\br= \br_t + \br_f$. We define
\begin{align*}
\Pi^{\calT}(\br_t)\coloneqq \{ M_t\in \Pi(\br_t)\mid M_t \in \calT \},
&&\Pi^{\calF}(\br_f)\coloneqq \{ M_f\in \Pi(\br_f)\mid M_f \in \calF \}.
\end{align*}
Since $\calF=\{M\in \rep{\Pi} \mid \Hom{\Pi}{I_w}{M}=0\}$ for some $w\in W$, it follows that $\Pi^{\calF}(\br)$ is an open subset of $\Pi(\br)$ as $\dim\Hom{\Pi}{I_w}{-}$ is an upper semi-continuous function. Since the set of torsion pairs bijectively corresponds to the Weyl group $W(C)$, any torsion class has a description of \eqref{eq:tors_w}, so that $\Pi^{\calT}(\br)$ is also open similarly. Finally, we define a quasi-affine algebraic variety $E(\br_t, \br_f )$ as the set of tuples $(M, M_t, M_f, f, g)$ such that $M\in \Pi(\br)$ and take $(M_t, M_f )\in \Pi^{\calT}(\br_t)\times \Pi^{\calF}(\br_f)$ with an exact sequence
$$0\rightarrow M_t \xrightarrow{f} M \xrightarrow{g} M_f \rightarrow 0$$ in $\rep{\Pi}$.
Then, we consider the diagram:
$$\Pi^{\calT}(\br_t)\times \Pi^{\calF}(\br_f)\xleftarrow{p} E(\br_t, \br_f )\xrightarrow{q} \Pi(\br),$$
where $p$ and $q$ are canonical projections.
\begin{lemma}\label{bdl_lem}
In the above setting, the projection $p$ is a locally trivial fibration with a smooth and connected fiber of dimension $\dim G(\M)-(\br_t, \br_f )$. The image of $q$ is equal to the set of $M\in \Pi(\br)$ such that $\rankvec{tM}=\br_t$. The non-empty fibers of $q$ are isomorphic to $G(\br_t)\times G(\br_f)$.
\begin{proof}
  Our assertion for $q$ follows from the uniqueness of the canonical short exact sequence for a $\Pi$-module $M$ and the torsion pair $(\calT, \calF)$ (\cf \cite[VI. Proposition 1.5]{ASS}). In fact, $q$ is a principal $G(\br_t)\times G(\br_f)$-bundle by considering the ordinary action on $E(\br_t, \br_f)$. We take points $(M_a), (M_{t, a} ), (M_{f, a})$ in $\Pi(\br)$, $\Pi^{\calT}(\br_t)$ and $\Pi^{\calF}(\br_f)$ respectively, and define locally free $\Pi$-module structures on $\bigoplus_{i\in Q_0} (H_i)^{r_i}$, $\bigoplus_{i\in Q_0} (H_i)^{{r_t}_i}$ and $\bigoplus_{i\in Q_0} (H_i)^{{r_f}_i}$, where $H_i\coloneqq K[\varepsilon_i]/(\varepsilon^{c_i})$. We consider the following complex from \cite[\S 12.2]{GLS1} which can be seen as a fiberwise description of tautological vector bundles on $\Pi^{\calT}(\br_t)\times \Pi^{\calF}(\br_f)$ over a point $((M_{t, a}), (M_{f, a}))$:
\begin{align*}
0\rightarrow &\bigoplus_{i\in Q_0}\Hom{H_i}{(H_i)^{{r_f}_i}}{(H_i)^{{r_t}_i}}\\ \xrightarrow{d^0_{M_f, M_t}} &\bigoplus_{(i, j)\in\overline{\Omega}}\Hom{H_i}{{}_i H_j\tens{H_j}(H_j)^{{r_f}_j}}{{(H_i)^{{r_t}_i}}}\\ \xrightarrow{d^1_{M_f, M_t}}&\bigoplus_{i\in Q_0}\Hom{H_i}{(H_i)^{{r_f}_i}}{(H_i)^{{r_t}_i}},
\end{align*}
where
\begin{align*}
(d^0_{M_f, M_t}(\phi_k)_{k\in Q_0}))_{(i, j)}&=(M_t)_{ij}\circ (\id_{{}_i  {}H_j}\tens{j}\phi_j)-\phi_i \circ {(M_f)}_{ij},\\
(d^1_{M_f, M_t}((\psi_{ij})_{(i, j)\in \overline{\Omega}}))_k
&= \sum_{j\in \overline{\Omega}(-, k)}\sgn{j}{k}(M_t)_{kj}\circ \ad_{jk}(\psi_{jk})-\psi_{kj}\circ\ad_{jk}((M_f)_{jk}).
\end{align*}
Note that the map $d^0_{M_f, M_t}$ has rank $\dim\Hom{S}{M_f}{M_t}-\dim\Hom{\Pi}{M_f}{M_t}$, and that $\Ext{1}{\Pi}{M_f}{M_t}\cong \Ker{d^1_{M_f, M_t}}/\Img{d^0_{M_f, M_t}}$. By the generalized Crawley-Boevey formula in Proposition \ref{PiExtduality} and the definition of torsion pairs, we compute:
\begin{align*}
\dim \Ker d^1_{M_f, M_t}&=\dim \Ext{1}{\Pi}{M_f}{M_t}+\rk d^0_{M_f, M_t}\\
=&\dim \Hom{\Pi}{M_f}{M_t}+ \dim\Hom{\Pi}{M_t}{M_f}\\
&-(\br_f, \br_t) +\dim\Hom{S}{M_f}{M_t}-\dim\Hom{\Pi}{M_f}{M_t}\\
=& \dim \Hom{S}{M_f}{M_t}-(\br_f, \br_t).
\end{align*}
In particular, $\dim \Ker{d^1_{M_f, M_t}}$ only depends on $\br_f$ and $\br_t$. Let $E$ be the set of all $S$-submodules of $E(\br_t, \br_f)$ that yields exact sequences isomorphic to $0\rightarrow M_f \xrightarrow{f} M \xrightarrow{g} M_t \rightarrow 0$ of $S$-modules. This is a homogeneous space for the action of $G(\M)$ and the stabilizer of a point of $(f, g)\in E$ is equal to $\{\id_M + fhg \mid h\in \Hom{S}{M_f}{M_t} \}$, and it is a smooth connected variety of dimension $\dim G(\M)- \dim \Hom{S}{M_f}{M_t}$. Now, we consider the fiber of $p$ over $(M_{t, a}, M_{f, a})\in \Pi^{\calT}(\br_t)\times \Pi^{\calF}(\br_f)$. The fiber consists of $(f, g)$ in $E$ and $M_a \in \Pi(\br)$ which satisfy compatibilities in the definition of morphisms of $\Pi$-modules. We determine $(M_a)$ from given $((M_{t, a}), (M_{f, a}))$ and $(f, g)$. Then, the space of possible choices for $(M_a)$ is isomorphic to $\Ker{d^1_{M_f, M_t}}$, and in fact defines a subsheaf of a tautological vector bundle on $\Pi^{\calT}(\br_t)\times \Pi^{\calF}(\br_f)$. The linear map $d^1_{M_f, M_t}$ depends smoothly on $((M_{a, t}), (M_{a, f}))$, and has constant rank by the above equality. In particular, the fiber of the sheaf defined by $\Ker{d^1_{M_f, M_t}}$ has constant rank. Therefore, our $p$ is a vector bundle by Nakayama's lemma. In particular, $p$ is a locally trivial fibration. Finally, we conclude that the dimension of fibers of $p$ is equal to $\dim E + \dim \Ker{d^1_{M_f, M_t}}=\dim G(\M)-(\br_{f}, \br_t)$.
\end{proof}
\end{lemma}
Now, we take a tuple of maximal irreducible components $(Z_t, Z_f)\in \T(\br_t) \times \F(\br_f)$. Then, $p^{-1}(Z_t \times Z_f)$ is an irreducible component of $E(\br_t, \br_f)$ by Lemma \ref{bdl_lem}. We consider the irreducible subset $q(p^{-1}(Z_t \times Z_f))$ of $\Pi(\br)$, which has dimension 
$$\dim(Z_t \times Z_f) + (\dim G(\M)-(\br_t, \br_f)) - \dim(G(\br_t) \times G(\br_f)).$$
Since we take $Z_t$ and $Z_f$ as maximal components, this dimension is equal to $G(\M) -q_{DC}(\br)$ by Proposition \ref{prop:dimofcrycom}. Namely, it is equal to the common dimension of maximal irreducible components of $\Pi(\br)$. Thus, $\overline{q(p^{-1}(Z_t \times Z_f))}$ is a maximal irreducible component of $\Pi(\br)$. Since the canonical short exact sequence for $M$ associated with a torsion pair is unique up to isomorphism, this ensures that crystal modules are closed under taking torsion parts or torsion-free parts. In particular, we obtain a map $\Xi\colon \T\times \F \rightarrow \B$ by running $(\br_t, \br_f )$.
\begin{theorem}
The map $\Xi\colon \T\times \F\rightarrow \B$ is a bijection.
\begin{proof}
For any $(\br_t, \br_f)$, the map $\Xi$ induces an injection from $\T(\br_t)\times \F(\br_f)$ to the set of all irreducible components of $\overline{q(E(\br_t, \br_f))}$. Now, we decompose $\Pi(\br)$ as a disjoint union of constructible subsets:
$$\Pi(\br) =\coprod_{\br_t+ \br_f=\br} q(E(\br_t, \br_f)).$$
Thus, any crystal irreducible component is contained in one and only one of maximal dimensional $\overline{q(E(\br_t, \br_f))}$. In particular, we obtain a bijection $\Xi\colon \coprod_{\br_t+ \br_f=\br} \T(\br_t)\times \F(\br_f) \rightarrow \B(\br)$.
\end{proof}
\end{theorem}
Let $Z\coloneqq \overline{q(p^{-1}(Z_t \times Z_f) )}$, and $M$ be a generic point of $Z$. Then, our construction shows that the torsion submodule $tM$ of $M$ has the rank-vector $\br_t$ and $(tM, M/tM)$ is a generic point in $Z_t \times Z_f$. In fact, $p^{-1}(U_t, U_f )$ is an open dense subset of $p^{-1}(Z_t \times Z_f)$ for a $G(\br_t)$-invariant open dense subset $U_t \subseteq Z_t$ and a $G(\br_f)$-invariant open dense subset $U_f \subseteq Z_f$. So, the subset $U\coloneqq q(p^{-1}(U_t, U_f))$ is open dense in $Z$, and this subset is constructible by Chevalley's theorem. If we have $M\in U$, then $\rankvec{tM}=\br_t$ and $(tM, M/tM)$ belongs to $U_t \times U_f$. 

Now, we consider the torsion pair $(\Fac{I_w}, \Sub{\Pi/I_w})$ for $w \in W$. Let $\T^w \times \F^w$ denote the set of pairs $(Z_1, Z_2)$ of irreducible components $Z_1, Z_2$, where a generic point of $Z_1$ (\resp a generic point of $Z_2$) belongs to $\Fac{I_{w}}$ (\resp $\Sub{\Pi/I_w}$). Then, we have the bijection $\Xi_w$ from $\T^w \times \F^w$ to $\B$.
\begin{prop}\label{strata-1step}
    Let $w\in W$ and $i\in Q_0$ such that $\length{ws_i}> \length{w}$. Then, the map $\Xi_{ws_i}$ restricts to a bijection $\T^{ws_i}\times \left(\F^{ws_i}\cap \T^w \right)\rightarrow \T^w$, and the map $\Xi_w$ restricts to a bijection $\left(\F^{ws_i}\cap \T^w \right)\times \F^w \rightarrow \F^{ws_i}$. In particular, we have the following commutative diagram:
    $$\begin{tikzpicture}[auto]
\node (a) at (0, 0) {$\T^{ws_i}\times \F^{ws_i}$}; 
\node (b) at (0, 1.8) {$\T^{ws_i}\times \left(\F^{ws_i}\cap \T^{w} \right)\times \F^{w}$};  
\node (c) at (6.4, 0) {$\mathfrak{B}$};
\node (d) at (6.4, 1.8) {$\T^{w}\times \F^{w}$};
\draw[->] (b) to node[swap] {\footnotesize{$\id \times \Xi_w$}} (a);
\draw[->] (b) to node {\footnotesize{$\Xi_{ws_i}\times\mathrm{id}$}} (d);
\draw[->] (a) to node[swap] {\footnotesize{$\Xi_{ws_i}$}} (c);
\draw[->] (d) to node {\footnotesize{$\Xi_{w}$}} (c);
\end{tikzpicture}.$$
    \begin{proof}
    Let $(Z_1, Z_2)\in \T^{ws_i}\times \F^{ws_i}$ and $Z\coloneqq \Xi_{ws_i}(Z_1, Z_2)$. We take a generic point $M$ of $Z$ and let $T_{ws_i}$  be the torsion submodule of $M$ with respect to $\left(\Fac{I_{ws_i}}, \Sub{\Pi/I_{ws_i}} \right)$ as in Theorem~\ref{thm:HNseq}. Then, $T_{ws_i}$ belongs to $\Fac{I_w}$ and $\left(T_{ws_i}, M/T_{ws_i} \right)$ is generic in $(Z_1, Z_2)$. Since any torsion class is closed under taking quotients and extensions, $M$ belongs to $\Fac{I_w}$ if and only if $M/T_{ws_i}$ belongs to $\Fac{I_w}$. Namely, $Z$ belongs to $\T^{w}$ if and only if $Z_2$ belongs to $\T^{w}$. Thus, $\Xi_{ws_i}$ restricts to a bijection $\T^{ws_i}\times \left(\F^{ws_i}\cap \T^w \right)\rightarrow \T^w$. Our assertion about $\Xi_w$ can be proved similarly to $\Xi_{ws_i}$. Let $(Z_1, Z_2, Z_3) \in \T^{ws_i}\times \left(\F^{ws_i}\cap \T^w \right)\times \F^w$, and let $Z_4\coloneqq \Xi_{ws_i}(Z_1, Z_2)$ and $Z\coloneqq \Xi_{w}(Z_4, Z_3)$. Then, the point $(T_w, M/T_w)$ is generic in $Z_4 \times Z_3$. In particular, $\left(T_{ws_i}, T_{w}/T_{ws_i}, M/T_w \right)$ is generic in $Z_1 \times Z_2 \times Z_3$ because $T_{ws_i}$ is the torsion submodule of $T_w$ with respect to the torsion pair $\left(\Fac{I_{ws_i}}, \Sub{\Pi/I_{ws_i}}\right)$. A similar argument for $(Z'_1, Z'_2, Z'_3)\coloneqq \left(\Xi_{ws_i}\circ \left(\id \times \Xi_{w} \right) \right)^{-1}(Z)$ shows that $\left(T_{ws_i}, T_{w}/T_{ws_i}, M/T_w \right)$ is generic also in $(Z'_1, Z'_2, Z'_3)$. This implies the commutativity of the diagram.
    \end{proof}
\end{prop}
\subsubsection*{Proof of Theorem \ref{strata-thm}}
In view of Proposition \ref{strata-1step} and its proof, we deal with any maximal green sequence. Namely, $\T^{s_{i_1}\cdots s_{i_{l-1}}} \times (\F^{s_{i_1}\cdots s_{i_{l-1}}}\cap \T^{s_{i_1}\cdots s_{i_{l-2}}} ) \times \cdots \times (\F^{s_{i_1}s_{i_{2}}}\cap \T^{s_{i_1}}) \times \F^{s_{i_1}}$ bijectively corresponds to $\B$. In particular, the bijection is given by the correspondence of generic points $$(T_{s_1\cdots s_{l-1}}, T_{s_1\cdots s_{l-2}}/T_{s_1 \cdots s_{l-1}}, \dots, M/T_{s_1})$$ and $M$. This shows Theorem \ref{strata-thm} under consideration of Corollary \ref{HNcor}.\qed

As an application of Theorem \ref{strata-thm}, we construct Saito reflections in our setting. In particular, we describe the PBW parametrizations of canonical bases in our setting. For classical cases, a dual argument to ours has appeared in the work of Jiang~\cite{Jia}.
\begin{theorem}\label{strata-refl}
Choose a reduced expression $\mathbf{i}\coloneqq (i_1, \dots , i_l)$ of $w \in W$. Let $\Pi_{\mathbf{i}}^{\mathbf{a}}$ be the constructible set of modules $M$ which have the HN filtration $M=M_0 \supsetneq M_1 \supsetneq \cdots \supsetneq M_l =0$ such that $M_{j-1}/M_{j} \cong (I_{s_{i_1}\cdots s_{i_{j-1}}}/I_{s_{i_1}\cdots s_{i_j}})^{\oplus a_j}$ for  $\mathbf{a} \coloneqq (a_1, \dots, a_l) \in \Z^r_{\geq 0}$. Then $\Sigma_{i_1} M \in \Pi_{\mathbf{i}'}^{\mathbf{a}'}$, for $\mathbf{i}'\coloneqq (i_2, \dots , i_l)$ and $\mathbf{a}' \coloneqq (a_2, \dots, a_l)$.
\begin{proof}
Let $M \in \Pi_{\mathbf{i}}^{\mathbf{a}}$. We have a HN filtration $M=M_0 \supsetneq M_1 \supsetneq \cdots \supsetneq M_l =0$ such that $M_{j-1}/M_{j} \cong (I_{s_{i_1}\cdots s_{i_{j-1}}}/I_{s_{i_1}\cdots s_{i_{j}}})^{\oplus a_j}$ for some $a_j \in \mathbb{Z}_{\geq 0}$. Let $N\coloneqq M_{1}$ be the torsion submodule of $M$ with respect to a torsion pair $(\Fac{I_{s_{i_1}}}, \Sub{\Pi/ I_{s_{i_1}}})$. Since we have an exact sequence $0 \rightarrow N \rightarrow M \rightarrow E_{i_1}^{\oplus a_1} \rightarrow 0$ by Theorem \ref{refl}, it is found that $\Sigma_{i_1} M\cong \Sigma_{i_1} N$ by $\Sigma_{i_1} E_{i_1} = 0$. Then, we have a filtration $N=N_0 \supsetneq N_1 \supsetneq \cdots \supsetneq N_{l-1}=0$ such that $N_{j-1}/N_{j}\cong M_{j}/M_{j+1} \cong (I_{s_{i_1}\cdots s_{i_{j}}}/I_{s_{i_1}\cdots s_{i_{j+1}}})^{\oplus a_{j+1}}\cong (I_{s_{i_1}\cdots s_{i_{j}}}\tens{\Pi}{E_{i_{j+1}}})^{\oplus a_{j+1}}$ for $j=1\dots, l-1$. Thus, $N_j/N_{j+1}$ has a trivial $i_1$-factor for $j>1$. Thus, we have a filtration $N'=N'_0 \supsetneq N'_1 \supsetneq \cdots \supsetneq N'_{l-1}=0$ where $N'_j\coloneqq \Sigma_{i_1} N_j$ and $N'_j/N'_{j+1} \cong \Sigma_{i_1} (I_{s_{i_1}\cdots s_{i_{j}}}\tens{\Pi}{E_{i_{j+1}}})^{\oplus a_{j+1}}\cong (I_{s_{i_2}\cdots s_{i_{j}}}\tens{\Pi}{E_{i_{j+1}}})^{\oplus a_{j+1}}$ because $I_{s_{i_2}\cdots s_{i_{j}}}\tens{\Pi}{E_{i_{j+1}}}$ has a trivial $i_1$-sub. Thus, we obtain that $\Sigma_{i_1} M \cong \Sigma_{i_1}N \in \Pi_{\mathbf{i}'}^{\mathbf{a}'}$.
\end{proof}
\end{theorem}
Saito~\cite[Corollary 3.4.8]{Sai} has defined mutually bijections between certain subsets of crystals $B(-\infty)$, called Saito reflections. Namely,
we define $\tilde{f}_{i}^{\max}(b) \coloneqq \tilde{f}_{i}^{\varphi_i (b)}(b)$ (\resp ${\tilde{f}_{i}^*}{}^{\max}(b) \coloneqq {\tilde{f}^*_{i}}{}^{\varphi_i^* (b)}(b)$) and $S_{i}(b)\coloneqq \tilde{e}_i^{\varepsilon_i^*(b)}{\tilde{f}^{*}_{i}}{}^{\max}(b)$ (\resp $S_i^*(b)\coloneqq (\tilde{e}_i^*)^{\varepsilon_i(b)}{\tilde{f}_{i}}{}^{\max}(b)$). Then, these $S_i$ and $S_i^*$ give mutual bijections
\begin{equation*}
\{b\in B(-\infty) \mid \varphi_i(b)=0\}
\xymatrix{\ar@<1ex>[r]^{S_i}&\ar@<1ex>[l]^{S_i^{*}}}
\{b\in B(-\infty) \mid \varphi_i^*(b)=0\}.
\end{equation*}
We compare the Saito reflection $S_i$ and our reflection functor $\Sigma_i$ based on the idea of Baumann-Kamnitzer~\cite{BK} for a classical setting. We prepare the following Propositions as analogies of the work of \cite{BK}:
\begin{prop}\label{prop:reflopendense}
 Let $Z\in \Irrmax{\Pi(\br))}$ and set $c=\varphi_i (Z)$ and $\br'=\br -c \alpha_i$. Let $V$ be an open dense $G(\br)$-invariant subset of $\tilde{f}_i^{\max}Z$. Then, $\{M\in Z \mid [\Sigma_i^- \Sigma_i M] \subseteq V\}$ contains an open dense subset in $Z$.
 \begin{proof}
 Note that $\Sigma_i^- \Sigma_i M$ is the torsion submodule of $M$ with respect to a torsion pair $(\Fac{I_{s_i}}, \Sub{\Pi/I_{s_i}})$. If $M\in \Efilt{\Pi}$ with $\rankvec{M}=\br$ and $\rk_i \fac_i M= c$, then $\Sigma_i^- \Sigma_i M$ is the unique submodule with the rank vector $\br'$ of $M$ up to isomorphism by Proposition \ref{refl}. Let $Y(\br', i, c)$ be the set of tuples $((N_a), (M_a), g)$ of $(N_a)\in \Pi(\br')_{i, 0}$, $(M_a)\in \Pi(\br)_{i, c}$ and an isomorphism $g\colon N\rightarrow \Sigma_i^- \Sigma_i M$. We have a diagram
 $$\Pi(\br')_{i, 0}\xleftarrow{p} Y(\br', i, c)\xrightarrow{q} \Pi(\br)_{i, c},$$
 where $p$ and $q$ are the first and the second projections respectively. 
 By a similar argument as in the Proof of Lemma \ref{bdl_lem}, we find that $p$ and $q$ are principal bundles, and that the set $p(q^{-1}(M))$ is the orbit $[\Sigma_i^- \Sigma_i M]$ for each $(M_a)\in \Pi(\br)_{i, c}$.
 
 Now, $Z\cap \Pi(\br)_{i, c}$ is an open dense subset of $Z$. By definition of $\tilde{f}^{\max}$, we obtain a restriction
 $$(\tilde{f}_i^{\max}Z)\cap \Pi(\br')_{i, 0}\xleftarrow{p} Y' \xrightarrow{q} Z\cap \Pi(\br)_{i,c},$$
 where $Y'$ is a maximal dimensional irreducible component of $Y(\br', i, c)$. We put $U\coloneqq Z\cap \Pi(\br)_{i,c}$ and $U'\coloneqq (\tilde{f}_i^{\max}Z)\cap \Pi(\br')_{i, 0}$. Since $U'$ is open dense in $\tilde{f}_i^{\max} Z$, it is found that $U' \cap V$ is open dense in $U'$. Thus, $q(p^{-1}(U'\cap V))$ is open dense in $U$. We observe that
 \begin{align*}
     \{M \in Z\mid [\Sigma_i^- \Sigma_i M]\subseteq V\} &=
     \{ M\in Z \mid [\Sigma_i^- \Sigma M] \,\,\textrm{meets}\,\, V\}\\
     &\subseteq \{M\in U \mid p(q^{-1}(M)) \,\,\textrm{meets}\,\,U'\cap V\}\\
     &= q(p^{-1}(U'\cap V))
 \end{align*}
 Since we take $V$ as $G(\br')$-invariant, we obtain our assertion.
 \end{proof}
\end{prop}
Let $b \in B(-\infty)$ and let $Z_b$ be the corresponding maximal irreducible component under Definition-Theorem \ref{defthm:crystal}. We take $\br \in R_+$. We define the set $\Theta(\br, i)$ of triples $((M_a), (N_a), h)$ such that $(M_a) \in \Pi(\br)_{i, 0}$, $(N_a)\in \Pi(s_i \br)_{i, 0}^*$ and $h\colon N \rightarrow \Sigma_i M$ is an isomorphism  between $\E$-filtered modules. We have a diagram
\begin{align}
    \Pi(\br)_{i, 0}\xleftarrow{r} \Theta(\br, i) \xrightarrow{s} (\Pi(s_i \br)_{i, 0})^*, \label{eq:bdldual}
\end{align}
where $r$ and $s$ are principal bundles over the images with structural groups $G(s_i \br)$ and $G(\br)$ respectively. In fact, for $(M_a)\in \Pi(\br)_{i, 0}$, the $H_i$-linear map $M_{i, \mathrm{in}}$ is full rank and the image is locally free. Thus, $(M_a)\mapsto \Ker{M_{i, \mathrm{in}}}$ is a continuous map from $\Pi(\br)_{i, 0}$ to a Grassmannian of the $\widetilde{M}_i$ in a similar manner to \cite[\S 3.4]{GLS4}. Thus, by taking an open covering of $\Pi(\br)_{i, 0}$, we can prove $r$ is a principal $G(s_i \br)$-bundle. By a similar argument, $s$ is a principal $G(\br)$-bundle.

In the crystal side, let $b\in B(-\infty)$ of weight $\br$ and $\varphi_i(b)=0$, and let $c= \varphi_i(b^*)$; $d=\varepsilon_i(b^*)$ and $\br'=\wt{(\tilde{f}_i^*)^{\max} b}$. Then, $\wt{b}=\br= \br'+ c\alpha_i$. A bit of calculation yields that $d=c -\langle \alpha_i, \br  \rangle \geq 0$ and that $\varphi_i((\tilde{f}_i^*)^{\max}b)=0$, which imply
\begin{align*}
   & d= \varphi_i(S_i b);
   &&c+d= 2c- \langle\alpha_i, \br' + c\alpha_i \rangle = -\langle \alpha_i, \br' \rangle.
\end{align*}
In view of this observation, we obtain the following proposition as an analogue of \cite[Theorem 5.3]{BK}.
\begin{prop}\label{prop:crydual}
 Let $b \in B(-\infty)$ of weight $\br$ and such that $\varphi_i (b)=0$. We have
 $$r^{-1}(Z_b \cap \Pi(\br)_{i, 0})=s^{-1}(Z_{S_i b} \cap (\Pi (s_i \br)_{i, 0})^*).$$
 \begin{proof}
 Let $M$ (\resp N) be an $\E$-filtered module such that $\fac_i(M)=0$ (\resp $\sub_i(N)=0$) and $N\cong \Sigma_i M$. We have an equality
 $$\rk_i M_i -\rk_i \sub_i M = \rk_i N_i -\rk_i \fac_i N.$$
 
 Let $\br' \in R^+$ be such that $\br= \br' + c\alpha_i$ and $s_i \br =\br' + d\alpha_i$ for $c, d \in \Z$. Note that $c+d= -\langle \alpha_i, \br' \rangle$. We take two $\E$-filtered modules $M, N$ such that
 \begin{align*}
     &\rankvec{M} = \br, 
     &&\fac_i{M} = 0;\\
     &\rankvec{N}=s_i \br,
     &&\sub_i{N}=0;
     &\Sigma_i M&\cong N
 \end{align*}
 Note that we have $\rk_i{\sub_i{M}}=c$ if and only if $\rk_i{\fac_i{N}}=d$. We restrict the diagram \eqref{eq:bdldual} to a suitable locally closed subset $\Theta(\br, i)_{c, d} \subseteq \Theta (\br, i)$. Namely, we obtain a diagram
 \begin{align}
     \Pi(\br)_{i, 0}\cap (\Pi(\br)_{i, c})^* \xleftarrow{r} \Theta (\br, i)_{c, d} \xrightarrow{s} (\Pi(s_i \br)_{i, 0})^* \cap \Pi(s_i \br)_{i, d}. \label{eq:bij_saitorefl}
 \end{align}
 For $\br \in R^+$, a vertex $i\in Q_0$ and $c\in \Z_{\geq 0}$, we set $\Pi(\br)_{i, c}^{\times}\coloneqq \Pi(\br)_{i, c}\cap (\Pi(\br)_{i, 0})^*$.
 
 Let $\Psi$ be the set of tuples $((L_a), (M_a), (N_a)), f, g, h)$ such that \begin{align}
     &(L_a)\in \Pi(\br')_{i, 0}^{\times};
     &&(M_a)\in (\Pi(\br)_{i, c}^{\times})^*;
     &(N_a)\in \Pi(s_i \br)_{i, d}^{\times}
     \end{align}
     and that $f\colon M \twoheadrightarrow L$, $g\colon L \hookrightarrow N$ and $h\colon N \xrightarrow{\cong} \Sigma_i M$ in $\Efilt{\Pi}$ which satisfy that $hgf \colon M\rightarrow \Sigma_i M$ is the canonical map in Remark \ref{rem:canmap}.
     
     $$\begin{tikzpicture}[auto]
\node (a) at (0, 0) {$\Pi(\br')_{i, 0}^{\times}$}; 
\node (b) at (4, 0) {$Y(\br', i, d)^{\times}$};
\node (c) at (8, 0) {$\Pi(s_i \br)_{i, d}^{\times}$};
\node (d) at (4, 2) {$\Psi$};
\node (e) at (8, 2) {$\Theta(\br, i)_{c, d}$};
\node (f) at (0, 4) {$(\Pi(\br')_{i, 0}^{\times})^*$};
\node (g) at (4, 4) {$(Y(\br', i, c)^{\times})^*$};
\node (h) at (8, 4) {$(\Pi(\br)_{i, c}^{\times})^*$};
\node (i) at (4.5, 2.5) {$\llcorner$};
\node (i) at (4.5, 1.5) {$\ulcorner$};
\draw[->] (b) to node[swap] {\footnotesize{$p^{\times}$}} (a);
\draw[->] (b) to node {\footnotesize{$q^{\times}$}} (c);
\draw[->] (d) to (b);
\draw[->] (d) to (e);
\draw[->] (d) to (g);
\draw[->] (g) to node[swap] {\footnotesize{$*p^{\times}*$}} (f);
\draw[->] (g) to node {\footnotesize{$*q^{\times}*$}} (h);
\draw[double distance=2pt] (f) to (a);
\draw[->] (e) to node[swap] {$r$} (h);
\draw[->] (e) to node {$s$} (c);
\end{tikzpicture}.$$

Note that $\Psi$ defines principal bundles, with structural groups $G(s_i \br)$, $G(\br')$ and $G(\br)$ in the above diagram as an analogue of \cite[\S 5.2]{BK}. In particular, all the maps are locally trivial fibrations with a smooth connected fiber in the diagram, so that they induce bijections between the sets of irreducible components.

By definition, we have $\varphi_i((S_i b)^*)=0$. The commutativity of the above commutative diagram implies that $Z_b \cap (\Pi(\br)_{i, c}^{\times})$ and $Z_{S_i b}\cap \Pi(s_i \br)_{i, d}^{\times}$ correspond in the bijection defined by \eqref{eq:bij_saitorefl}. Since we take $Z_{b}$ and $Z_{S_i b}$ as crystal components, Theorem \ref{thm:maxcomp} yields our assertion by taking closures in $\Pi(\br)_{i, 0}$ and $(\Pi(s_i \br)_{i, 0})^*$.
 \end{proof}
\end{prop}
We have the following analogue of \cite[Proposition 5.5]{BK} to our setting.
\begin{prop}\label{prop:saitopendense}
 Let $b\in B(-\infty)$ and set
 \begin{align*}
     &\br=\wt{b}; &&c=\varphi_i(b);
     &&U=Z_b\cap \Pi(\br)_{i, c};
     &&b'=S_i (\tilde{f}_i^{\max} b);
     &&\br'=\wt{b'}
 \end{align*}
 Let $V$ be an open dense $G(\br')$-invariant subset of $Z_{b'}$. Then, $\{(M_a)\in U\mid [\Sigma_i M] \subseteq V\}$ contains an open dense subset of $U$.
 \begin{proof}
 If $c=0$, we have $b'= S_i b$ and $\br' = s_i \br$. Then, we use the diagram \eqref{eq:bdldual} and Proposition \ref{prop:crydual} to obtain our assertion by an analogue of Proposition \ref{prop:reflopendense}. In general case, we use $\Sigma_i M \cong \Sigma_i(\Sigma_i^- \Sigma_i M)$ and Proposition \ref{prop:reflopendense} to obtain our assertion from $c=0$ case.
 \end{proof}
\end{prop}
As a conclusion of this subsection, we have the following Theorem about Saito reflections:
\begin{theorem}[Saito reflections]\label{Saito-refl}
Let $\overline{\Pi_{\mathbf{i}}^{\mathbf{a}}}=Z_{\mathbf{i}}^{\mathbf{a}}$. We put $\mathbf{i'}=(i_2, \dots, i_l)$ and $\mathbf{a'}=(a_2, \dots, a_l)$. Then, we have $ S_{i_1} \tilde{f}_{i_1}^{\max}(Z_{\mathbf{i}}^{\mathbf{a}})=Z_{\mathbf{i}'}^{\mathbf{a}'}$.
\begin{proof}
Keep with the notation in Theorem~\ref{strata-thm}. For any $M\in \Pi_{\mathbf{i}}^{\mathbf{a}}$, it is found that $T_{j-1}/T_{j} \cong (I_{s_{i_1}\cdots s_{i_{j-1}}}/I_{s_{i_1}\cdots s_{i_{j}}})^{\oplus a_j}$ has a trivial $i_1$-factor for $j=2, \dots, l-1$. Thus, we have $\fac_{i_1}M\cong E_{i_1}^{\oplus a_1}$. Now, by our construction, the irreducible component $Z_{\mathbf{i}}^{\mathbf{a}}$ has an open dense subset contained in $\Pi_{\mathbf{i}}^{\mathbf{a}}$, so that we obtain $\varphi^*_{i_1}Z_{\mathbf{i}}^{\mathbf{a}}=a_1$.
Now, we take a generic point $M\in \Pi_{\mathbf{i}}^{\mathbf{a}}$ and assume $M$ belongs to an open dense subset $U$. We find that the set of $M_a\in U$ such that $\Sigma_{i_1}M_a$ is generic in $S_{i_1} {\tilde{f}_{i_1}}^{\max}(Z_{\mathbf{i}}^{\mathbf{a}})$, contains an open dense subset in $U$ by Proposition \ref{prop:saitopendense}. Thus, we obtain our assertion by Theorem \ref{strata-refl}.
\end{proof}
\end{theorem}
In particular, we obtain the PBW parametrization as a corollary of Theorem \ref{Saito-refl}.
\begin{cor}
${\tilde{f}_{i_r}}^{\max}S_{i_{r-1}} {\tilde{f}_{i_{r-1}}}^{\max}\cdots S_{i_1} {\tilde{f}_{i_1}}^{\max}Z_{\mathbf{i}}^{\mathbf{a}}= Z_0$.\qed
\end{cor}
\begin{rem}
    In Theorem~\ref{strata-thm}, if we take an element $w\in W$, we can obtain an injection $\psi_{\mathbf{i}}\colon \Z^r_{\geq 0} \rightarrow \mathcal{B}$ given by $\mathbf{a}\mapsto Z_{\mathbf{i}}^{\mathbf{a}}$ and the image of $\psi_{\mathbf{i}}$ does not depend on the choice of $\mathbf{i}$. In particular, we obtain the PBW parametrization of canonical bases of the quantum unipotent subalgebra $U_q(\mathfrak{n}(w))$ with respect to $w$ of any symmetrizable finite Dynkin type.
\end{rem}
\subsection{MV polytopes from generalized preprojective algebras}
 
In this subsection, we give a description of a version of Harder-Narasimhan polytopes for $\Pi$. Namely, we give a description of \emph{generalized Harder-Narasimhan polytopes} (see Definition~\ref{def:gHN}) by using $g$-matrices developed in \S \ref{subsec:gvec}. We prove that they realize Mirkovi\'c-Vilonen polytopes. 
First, we recall the Auslander-Reiten inner product for the Grothendieck group:
\begin{def-thm}[{Auslander-Reiten~\cite[Theorem 1.4]{AR}}]\label{defthm:ARprod} 
    Let $M$ and $N$ be $\Pi$-modules. Then, we have 
    $$\langle g(M), \dimvc{N} \rangle = \dim \Hom{\Pi}{M}{N} - \dim \Hom{\Pi}{N}{\tau M},$$
    where the left hand side is the standard inner product on $\mathbb{Z}^n \times \mathbb{Z}^n$.
\end{def-thm}
This pairing is extended to a group homomorphism $\theta \colon K_0(\rep{\Pi})\rightarrow \R$ for $\theta\in \R^n \cong K_0\left(\proj{\Pi} \right)\tens{\Z}\R$. We have a realization of Weyl chambers via $g$-matrices for $\tau$-tilting pairs by Theorem \ref{Tits}. In our setting, we want to consider rank vectors instead of dimension vectors, so that we introduce the following inner product:
\begin{def-pro}\label{def-pro:ARproduct}
    Let $M$ be a locally free $\Pi$-module with $\alpha \coloneqq \rankvec{M}$ and let $(I_w, P_w)$ be the $\tau$-tilting pair with respect to $w\in W$ and $g(w)=(g_1, \dots, g_n)$ be the $g$-matrix for this pair. Let $\gamma =\sum_{i\in Q_0} a_i g_i\in \Z^n$ with $a_i \in \Z_{\geq 0}$. If we regard $\alpha$ as a positive root and $\gamma$ as an integral weight $\sum_{i\in Q_0}a_iw \varpi_i$ by Theorem \ref{Tits}, then the pairing $\langle \gamma, \alpha \rangle \coloneqq \langle \gamma, \alpha\rangle_C$ (see \S\ref{pre:rootsys})
    coincides with the natural standard inner product given by 
    \begin{align*}
    \langle \gamma, \alpha \rangle =&\langle  \sum_{i\in Q_0} a_i g_i, \alpha \rangle \\
    =&\sum_{I_w e_i \neq 0} \frac{a_i}{c_i}(\dim \Hom{\Pi}{I_w e_i}{M}-\dim \Hom{\Pi}{M}{\tau (I_w e_i)})\\
    &-\sum_{I_w e_j = 0}\frac{a_j}{c_j}(\dim \Hom{\Pi}{P_w e_j}{M}),
    \end{align*}
    where $\diag{c_1, \dots, c_n}$ is a symmetrizer of $C$.
    \begin{proof}
    This is a direct consequence of the fact that 
    $$\rankvec{M}=(\rk_1 (e_1 M), \dots, \rk_n (e_n M))= ((\dim e_1 M)/c_1, \dots, (\dim e_n M)/c_n),$$
    Theorem \ref{Tits}, and Definition-Proposition \ref{defthm:ARprod}. Note that we have $\tau P=0$ for any projective module $P$. We obtain an equality
    \begin{align*}
    \langle  \sum_{i\in Q_0} a_i g_i, \alpha \rangle 
    =&\sum_{I_w e_i \neq 0} \frac{a_i}{c_i}(\dim \Hom{\Pi}{I_w e_i}{M}-\dim \Hom{\Pi}{M}{\tau (I_w e_i)})\\
    &-\sum_{I_w e_j = 0}\frac{a_j}{c_j}(\dim \Hom{\Pi}{P_w e_j}{M}-\dim \Hom{\Pi}{M}{\tau (P_w e_j)})\\
    =&\sum_{I_w e_i \neq 0} \frac{a_i}{c_i}(\dim \Hom{\Pi}{I_w e_i}{M}-\dim \Hom{\Pi}{M}{\tau (I_w e_i)})\\
    &-\sum_{I_w e_j = 0}\frac{a_j}{c_j}(\dim \Hom{\Pi}{P_w e_j}{M}).\qedhere
    \end{align*}
    \end{proof}
\end{def-pro}
\begin{definition}\label{def:D}
    Let $w\in W$ and let $\gamma \in P^+$. Let $M \in \rep{\Pi}$. We write $\pm w\gamma \coloneqq \pm \sum_{i\in Q_0} \langle \gamma, \alpha_i \rangle w\varpi_i \in \Z^n$.
    $$D_{\pm w\gamma}(M)\coloneqq \sum_{i\in Q_0} \frac{\langle \gamma, \alpha_i \rangle}{c_i}\dim\Hom{\Pi}{N(\pm w\varpi_i)}{M}\,\,(\geq 0).$$
\end{definition}
We have the following propositions by similar arguments as in \cite[Proposition 4.1, 4.2, 4.3]{BK}:

\begin{prop}\label{prop:Dformula}
 Let $\gamma$ be a weight such that $\langle \gamma, \alpha_i \rangle \leq 0$ for $ i\in Q_0$. Then, we have the following equalities for a crystal module $T\in \Efilt{\Pi}$:
 \begin{align}
     D_\gamma(T)&= D_{s_i \gamma}(\Sigma_i T);\label{eq:Dformula1}\\
     D_{s_i \gamma}(T)&=D_{\gamma}(\Sigma_i^- T)- \langle \gamma, \rankvec{(\sub_i T)}\rangle ; \label{eq:Dformula2}\\
     D_{\gamma}(T)&= D_{\gamma}(\Sigma_i^- \Sigma_i T).\label{eq:Dformula3}
 \end{align}
 \begin{proof}
 Note that we can regard $T$ as a module in $\rep{\Pi}$ or $\rep{\affPi}$. By Definition \ref{def_reflection}, this difference does not matter when we consider the module $\Sigma _{i} T$. The first equality \eqref{eq:Dformula1} follows from the following isomorphisms:
 \begin{align*}
     \Hom{\Pi}{N(\gamma)}{T}
     &\cong \Hom{\affPi}{\widetilde{N}(\gamma)}{T}\\
     &\cong \Hom{\affPi}{\Sigma_i^- \widetilde{N}(s_i \gamma)}{T}\\
     &\cong \Hom{\affPi}{\widetilde{N}(s_i \gamma)}{\Sigma_i T}.
 \end{align*}
 
 Next, we consider the second equality \eqref{eq:Dformula2}. In the above isomorphisms, we obtain, by substituting $T= \Sigma_i^- T$,
 $$D_{\gamma}(\Sigma_i^-T)= D_{s_i \gamma}(\Sigma_i \Sigma_i^- T).$$
 Now, Corollary \ref{cor:dimNhomext} says that
 \begin{align*}
     &\frac{1}{c_i}\dim\Hom{\affPi}{\widetilde{N}(s_i \gamma)}{E_i}=-\langle \gamma, \alpha_i \rangle; &&\frac{1}{c_i}\dim \Ext{1}{\affPi}{\widetilde{N}(s_i\gamma)}{E_i}=0.
 \end{align*}
 Since $T$ is a crystal module, the second equality \eqref{eq:Dformula3} is obtained by applying the functor $\Hom{\affPi}{\widetilde{N}(s_i \gamma)}{-}$ to the short exact sequence
 $$0\rightarrow \sub_i{T} \rightarrow T \rightarrow \Sigma_i \Sigma_i^- T \rightarrow 0.$$
Finally, we consider the third equality \eqref{eq:Dformula3}. We have an exact sequence
\begin{align}
0\rightarrow \Sigma_i^- \Sigma_i T \rightarrow T \rightarrow \fac_i T \rightarrow 0 \label{eq:seq5}
\end{align}
and an equality
$$\frac{1}{c_i} \dim \Hom{\Pi}{\widetilde{N}(\gamma)}{T}=0.$$
Thus, we obtain the third equality by applying the functor $\Hom{\Pi}{\widetilde{N}(\gamma)}{-}$ to \eqref{eq:seq5}. 
 \end{proof}
\end{prop}

\begin{prop}[edge relations]\label{prop:edge}
 Let $T\in \Efilt{\Pi}$ be a crystal module, and let $i\in Q_0$ and $w\in W$. We have an inequality
 \begin{equation}
 D_{-w\varpi_i}(T)+D_{-ws_i \varpi_i}(T)+\sum_{j\neq i}c_{ji}D_{-w\varpi_j}(T)\geq 0.\label{eq:edgerel}
 \end{equation}
 \begin{proof}
 Note that an equality $s_i \varpi_i =\varpi_i -\alpha_i$ is equivalent to
 \begin{equation}
     \varpi_i + s_i \varpi_i +\sum_{j\neq i}c_{ji} \varpi_j =0.
 \end{equation}
 Let $L_w(T)$ be the left hand side of \eqref{eq:edgerel}. Since $L_{s_i w}(T)=L_{w}(T)$, we may assume $\length{ws_i}>\length{w}$ without loss of generality. We take $w$ with $\length{w}> 0$, and write $w = s_k v$ with $k\in Q_0$ and $\length{v}= \length{w}-1$. Note that we have an inequality of length
 \begin{align}
     &\length{s_k v}> \length{v},
     &&\length{s_k v s_i}> \length{w}= \length{v}+1 \geq \length{vs_i}.
 \end{align}
 We calculate as follows:
 \begin{align*}
     L_w(T)&=D_{-w\varpi_i}(T)+D_{-ws_i \varpi_i}(T)+\sum_{j\neq i}c_{ji}D_{-w\varpi_j}(T)\\
     &=D_{-s_k vs_i \varpi_i}(T) + D_{-s_k v\varpi_i}(T)+ \sum_{j\neq i}c_{ji} D_{-s_k v \varpi_j}(T)\\
     &=D_{-vs_i \varpi_i}(\Sigma_k^- T) + D_{-v\varpi_i}(\Sigma_k^- T)+ \sum_{j\neq i}c_{ji} D_{-v \varpi_j}(\Sigma_k^- T)\\
     &=L_{v}(\Sigma_k^- T).
 \end{align*}
 By calculating inductively, for a reduced expression $w= s_{i_1}\cdots s_{i_r}$, we have
 $$L_w(T)=L_{1}(\Sigma_{i_r}^- \cdots \Sigma_{i_1}^- T).$$
 Since $D_{\gamma}(M)=0$ for any anti-dominant weight $\gamma$ and for any $M\in \rep{\Pi}$, we obtain
 \begin{equation*}
 L_w(T)=D_{-s_i \varpi_i}(\Sigma_{i_r}^- \cdots \Sigma_{i_1}^- T)\geq 0. \qedhere
 \end{equation*}
 \end{proof}
\end{prop}

\begin{prop}\label{D-Dprop}
    Let $M\in \rep{\Pi}$ be an $\E$-filtered module. We have $$D_{\gamma}(M)-D_{-\gamma}(M^*)=\langle \gamma, \rankvec{M}\rangle.$$
    \begin{proof}
    For each $i\in Q_0$ and each $M\in \Efilt{\Pi}$, the number $D_{\varpi_i}(M)$ is equal to the multiplicity of $E_i$ in a filtration in Definition \ref{def:Efilt}. By definition, this is equal to $\langle \varpi_i, \rankvec{M} \rangle$. Since we have $D_{-\varpi_i}(M)=0$, we have a desired equality for $\gamma=\varpi_i$.
    We fix an anti-dominant weight $\lambda$ and let $\gamma\coloneqq w\lambda$. Then, we put $F_w(M)=D_{w\lambda}(M)-D_{-w\lambda}(M^*)-\langle \gamma, \rankvec{M}\rangle$.
    We assume that $\length{s_i w}>\length{w}$. Now, we have $\Sigma_i^- \widetilde{N}(-w\lambda)\cong \widetilde{N}(-s_i w\lambda)$.
    Then, we have 
    \begin{align*}
        \Hom{\affPi}{\widetilde{N}(-s_i w\lambda)}{M^*}
        &\cong \Hom{\affPi}{\widetilde{N}(-w\lambda)}{\Sigma_i M^*}\\
        &\cong \Hom{\affPi}{\widetilde{N}(-w\lambda)}{\left(\Sigma^-_i M\right)^*},
    \end{align*}
    so that $D_{-s_i w\lambda}(M^*)=D_{-w\lambda}((\Sigma^-_i M)^*)$. 
    Now, we use \eqref{eq:Dformula2} and Proposition \ref{prop:I_wrankvec} to obtain
    $$\rankvec{M}= \rankvec{\Sigma_i \Sigma_i^- M}+ \rankvec{\sub_i M}=s_i(\rankvec{\Sigma_i^- M})+\rankvec{\sub_i M}.$$
    This shows that the equality $F_{s_i w}(M)=0$ is equivalent to $F_{w}(\Sigma_i^- M)=0$. Thus, we obtain our assertion by induction on $\length{w}$.
    \end{proof}
\end{prop}
Then, we review the definition and basic properties of Harder-Narasimhan polytopes of objects in Abelian categories developed in the work of Baumann-Kamnitzer-Tingley~\cite{BKT}.
\begin{definition}
    Let $K_0(\mathcal{A})$ be the Grothendieck group of an Abelian category $\mathcal{A}$. For $T \in \mathcal{A}$, we consider the convex hull $\Pol{T}$ in $K_0(\mathcal{A})\tens{\mathbb{Z}}\mathbb{R}$ of all of points $[X]$ defined by subobjects $X \subseteq T$. We refer to this $\Pol{T}$ as the \emph{HN (=Harder-Narasimhan) polytope} of $T$.
\end{definition}
\begin{theorem}[{\cite[\S 3.2]{BKT}}] \label{thm:HNproperty}
Under the setting of the above definition, we have the following:
\begin{enumerate}
    \item The HN polytope $\Pol{T}$ has $P_\theta=\{ x\in \Pol{T}\mid \langle \theta, x \rangle =\psi_{\Pol{T}}(\theta)\}$ as faces, where $\theta \in K_0(\mathcal{A}\tens{\Z}\R)^*$ and $\psi_{\Pol{T}}$ is a linear function $\R^n \rightarrow \R$ which $\theta$ maps to its maximum on $\Pol{T}$.
    \item If $x$ is a vertex of the HN polytope $\Pol{T}$, then there is a unique subobject $X$ of $T$ such that $[X]=x$ up to isomorphism.
    \item Let $T=T_0 \supsetneq T_1 \supsetneq \cdots \supsetneq T_l =0$ be a HN filtration in $\mathcal{A}$. Then, each $[T_i]$ is a vertex of $\Pol{T}$.
\end{enumerate}
\end{theorem}
We define torsion classes and torsion-free classes determined by stability conditions as follows:
\begin{definition}
    Let $\theta \colon K_0(\mathcal{A})\rightarrow \R$ for $\theta\in \R^n$ be a group homomorphism. Then, we define
    \begin{align*}
    \calT_\theta&\coloneqq \left\{M \in \mathcal{A} \mid \langle \theta, [N] \rangle > 0 \,\, \textrm{for any quotient $N$ of $M$}\right\}\\  \overline{\calT_\theta} &\coloneqq \{M \in \mathcal{A} \mid \langle \theta, [N] \rangle \geq 0 \,\, \textrm{for any quotient $N$ of $M$}\}\\ 
   \calF_\theta &\coloneqq \{M \in \mathcal{A}\mid \langle \theta, [L] \rangle < 0\,\, \textrm{for any submodule $L$ of $M$}\}\\ \overline{\calF_\theta} &\coloneqq \{M \in \mathcal{A}\mid \langle \theta, [L] \rangle \leq 0\,\, \textrm{for any submodule $L$ of $M$}\}
    \end{align*}
    We refer to objects in $\mathcal{R}_\theta \coloneqq \overline{\calT_\theta} \cap \overline{\calF_\theta}$ as \emph{$\theta$-semistable objects}.
\end{definition}
It is well-known that $(\overline{\calT}_\theta, \calF_\theta)$ and $(\mathcal{T}_\theta, \overline{\calF}_\theta)$ are torsion pairs with $\calT_{\theta}\subseteq \overline{\calT}_{\theta}$ in $\mathcal{A}$ for any $\theta\in \R^n$. Here, we review a connection between stability conditions and $g$-vectors.  By the work of Br\"ustle-Smith-Treffinger~\cite[Proposition 3.27, Remark 3.28]{BST} (see also Yurikusa~\cite{Yur}) which developed stability conditions for general finite dimensional algebras by the $\tau$-tilting theory, we know the following proposition by applying their theorem to our setting by Definition-Proposition \ref{def-pro:ARproduct}:
\begin{prop}\label{gmatrix_tors}
    Let $w\in W$ and let $\theta \in C(w)$. We define
    $$\overline{\calT_\theta} \coloneqq \{M \in \rep{\Pi} \mid \langle \theta, [N] \rangle \geq 0 \,\, \textrm{for any quotient $N$ of $M$}\}$$
    and $$\calF_\theta \coloneqq \{M \in \rep{\Pi}\mid \langle \theta, [L] \rangle < 0\,\, \textrm{for any submodule $L$ of $M$}\}.$$ Then, $\overline{\calT}_\theta = \Fac{I_w}$ and $\calF_\theta=\Sub{\Pi/I_w}$ for $\theta\in C(w)$.\qed
\end{prop}
We can consider a three-step filtration $0\subseteq \Tmin{\theta}\subseteq \Tmax{\theta}\subseteq M$ of torsion submodules $\Tmin{\theta} \in \calT_{\theta}$ and $\Tmax{\theta} \in \overline{\calT}_\theta$ for $M\in \rep{\Pi}$. We compare this operation and Weyl group combinatorics in the case of $\mathcal{A}=\rep{\Pi}$. Now, we relate these torsion pairs with the Weyl group.
\begin{definition}
    For $J\subseteq Q_0$, we define $$F_J=\left\{\theta\in \R^n\mid \langle \theta, \alpha_j\rangle =0 \,\,\text{for any $j\in J$ and} \,\,\left\langle \theta, \alpha_i \right\rangle > 0 \,\,\text{for any $i\in Q_0 \backslash J$}\right\}.$$
    We say that $w\in W$ is $J$-reduced on the right if $\length{ws_j}>\length{w}$ for any $j\in J$. We define a parabolic subgroup $W_J \coloneqq \langle s_j \mid j\in J \rangle$ of $W$. Let $w_J$ be the longest element of $W_J$.
\end{definition}
The following Lemma is an analogue of a result of Buan-Iyama-Reiten-Scott~\cite[Proof of Theorem III.3.5]{BIRS}.
\begin{lemma}\label{lem:submodule}
Let $J\subseteq Q_0$. Let $C'$ and $D'$ respectively be the GCM  and its symmetrizer associated with the quiver determined by $Q_0 \backslash J$. Let $\Pi' \coloneqq \Pi(C', D')$. Then, we have an algebra isomorphism $\Pi' \cong \Pi/ I_{w_J}$. In particular, $\Sub{\Pi/I_{w_J}}= \rep{\Pi'}$.
\begin{proof}
Let $I_{Q'}\coloneqq \Pi (\sum_{i\in Q_0 \backslash J} e_i)\Pi$. Namely, we have $\Pi/I_{Q'} \cong \Pi'$. We enough to prove $I_{w_J}=I_{Q'}$. By definition of $I_{Q'}$, the two-sided ideal $I_{Q'}$ contains all two-sided ideals $I$ such that any composition factor of $\Pi/ I$ is $S_i$ for $i\in J$. Since $w_J$ is the longest element of $W_J$, we have $\length{w_J s_i}< \length{w_J}$. We have $I_{w_J}I_i= I_{w_J}$ for any $i \in J$ by Theorem \ref{ppoofsttiltofpi}. On the other hand, any composition factor of $\Pi/ I_{w_J}$ is $S_i$ for $i\in J$ by our choice of $w_J$ and definition of $I_{w_J}$. Thus, we obtain our assertion.
\end{proof}
\end{lemma}
The following result is an analogue of the work of Baumann-Kamnitzer-Tingley~\cite[Theorem 5.18 and Remark 5.19]{BKT} to our setting:
\begin{prop}\label{Jred_tors}
 Keep the notation in Theorem~\ref{thm:HNseq}. Let $\theta\in F_J$ and let $w\in W$ such that $w$ is $J$-reduced on the right.
\begin{enumerate}
    \item For any crystal module $M$, we have $\Tmin{w\theta} =T_{ww_{J}}$. In particular, we have an equality $$(\calT_{w\theta}, \overline{\calF}_{w\theta})=(\Fac{I_{ww_J}}, \Sub{\Pi/I_{ww_J}}).$$
    \item For any $w\in W$, we have $\mathcal{R}_{w\theta}\cong \rep{\Pi/I_{w_J}}$. In particular, $\mathcal{R}_{w\theta}$ is in common for all possible $w$.
\end{enumerate}
\begin{proof}
For any crystal module $M$, we have enough to show that $T_{w w_{J}} \in \calT_{w\theta}$ and $M/T_{w w_{J}}\in \overline{\calF}_{w\theta}$ by Theorem \ref{strata-1step}. We have $T_{ww_J} = I_w \tens{\Pi}{X}$ for some locally free module $X\in \calF_w \cap \Fac{I_{w_J}}$ by arguments in \S \ref{subsec:refl}. If we assume $T_{w w_{J}}\neq0$, then $X\notin \Sub{\Pi/I_{w_J}}$. So, $\rankvec{X}$ cannot be written as a linear combination of elements in $\{\alpha_j\mid j\in J\}$ , and $\langle w\theta, [T_{ww_J}] \rangle=(w\theta, \rankvec{T_{w w_J}})_C= ( \theta, \rankvec{X})_C >0$ by $\theta \in F_J$ and Proposition \ref{prop:I_wrankvec}. Thus, we obtain $T_{w w_J}\in \calT_{w\theta}$. On the other hand, $M/T_{w w_J}$ has a rank vector in $\mathbb{N}$-linear combination of roots in $N_{w w_J}=\{s_{i_1}\cdots s_{i_{k-1}}\alpha_{i_k}\mid 1\leq k\leq r\}$ for a reduced expression $ww_J=s_{i_1}\cdots s_{i_r}$ by Theorem~\ref{HNcor} and Theorem~\ref{VM_iso}. We have $\langle w\theta, [M/T_{w w_J}]\rangle = ( w\theta, \rankvec{M/T_{w w_J}})_C =(ww_J\theta, \rankvec{M/T_{w w_J}})_C \leq 0$. So, we have our assertion. The second assertion follows from the equivalence $\Sub{\Pi/I_{w_J}}\cong \Sub{\Pi/I_{ww_J}}\cap \Fac{I_w}$ for $\length{u}+\length{v}=\length{uv}$ by Proposition \ref{prop:intersecTF}. Thus, we obtain $\mathcal{R}_\theta \cong \Sub{\Pi/I_{w_J}} \cong \rep{\Pi/I_{w_J}}$ by the first assertion, Proposition \ref{gmatrix_tors} and Lemma \ref{lem:submodule}.
\end{proof}
\end{prop}

\begin{rem}\label{rem:invHN}
By the $K$-duality, we obtain an involution of HN polytopes determined by
$$\overline{\calT}_{w\theta}=\Fac{I_w}\mapsto \Sub{\Pi/I_{ww_0}}=\overline{\calF}_{-w\theta},$$
$$\calF_{w\theta}=\Sub{\Pi/I_w}\mapsto \Fac{I_{ww_0}}=\calT_{-w\theta},$$
and
$$\mathcal{R}_{w\theta}\mapsto \mathcal{R}_{-w\theta}.$$
\end{rem}

We define a \textit{generalized Harder-Narasimhan polytope} as follows:
\begin{definition}\label{def:gHN}
    For each maximal irreducible component $Z$ of $\Irrmax{\Pi(\mathbf{r})}$ and a generic module $T$ of this component, we define $$P(T)\coloneqq \left\{ \left(x_1/c_1, \dots, x_n/c_n\right)\mid x=(x_1, \dots, x_n) \in \Pol{T}\right\},$$ where $\diag{c_1, \dots, c_n}$ is a symmetrizer of $C$. We call this polytope the \emph{generalized HN polytope} of $T$. Let $\psi_{P(T)}$ be a linear function $\R^n \rightarrow \R$ which $\theta$ maps to its maximum on $P(T)$.
\end{definition}
By definition, generalized HN polytopes can be thought as a kind of HN polytopes for the category of locally free modules. By our construction, $P(T)$ for a crystal generic module $T$ is GGMS (=Gelfand-Goresky-MacPherson-Serganova) in the sense of \cite[\S 2.6]{BKT}. In particular, we have for any $\theta \in \R^n$ $$\psi_{P(T)}(\theta)=\langle \theta, \rankvec{T_\theta^{\min}}\rangle=\langle \theta, \rankvec{T_\theta^{\max}}\rangle.$$
We have the following analogue of \cite[Corollary 5.20]{BKT}.
\begin{theorem}\label{supfun}
Let $J\subseteq Q_0$, let $\theta \in F_J$, and let $w\in W$ be $J$-reduced. We assume $\theta \coloneqq  \sum_{i\in Q_0} \langle \theta, \alpha_i \rangle \varpi_i$ is an dominant integral weight. Then, we have an
equality
$$D_{\pm w\theta}(M)=\psi_{P(M)}(\pm w\theta). $$
\begin{proof}
We put $\psi'_{P(M)}( w\theta) \coloneqq ( \theta, \rankvec{T_\theta^{\max}})_C$ and
$$D'_{w\theta}(M)\coloneqq \sum_{i\in Q_0} \langle \theta, \alpha_i \rangle\dim\Hom{\Pi}{N( w\varpi_i)}{M}.$$
We calculate the $D'_{w\theta}(M)$:
\begin{align*}
    &\sum_{i\in Q_0} \langle \theta, \alpha_i\rangle\dim \Hom{\Pi}{N( w\varpi_i)}{M}\\
    =&\sum_{i\in Q_0}\langle \theta, \alpha_i \rangle \dim \Hom{\Pi}{I_w\tens{\Pi}{\Pi e_i}}{M}\\
    =&\sum_{i\in Q_0} \left\langle \theta, \alpha_i \right\rangle \dim \Hom{\Pi}{\Pi e_i}{\Hom{\Pi}{I_w}{M}}\\
    =&\sum_{i\in Q_0}\frac{1}{c_i} \left( \theta, \alpha_i \right)_C \dim e_i \Hom{\Pi}{I_w}{M}\\
    =&(\theta, \rankvec{\Hom{\Pi}{I_w}{M}} )_C\\
    =&( w\theta, \rankvec{I_w\tens{\Pi}{\Hom{\Pi}{I_w}{M}}} )_C\,\,\quad\textit{(Proposition \ref{prop:I_wrankvec})}\\
    =&(w\theta, \rankvec{M_w})_C \,\,\quad\textit{(Remark \ref{rem:torfF_w})}\\
    =& (w\theta, \rankvec{M^{\max}_{w\theta}})_C\,\,\quad \textit{(Proposition \ref{gmatrix_tors})}\\
    =&\psi'_{P(M)}(w\theta)\,\,\quad\textit{(cf. \cite[Proposition 2.16]{BKT})},
\end{align*}
where $M_w$ is the torsion submodule of $M$ with respect to the torsion class $\Fac{I_w}$. We obtain our assertion $D_{w\theta}(M)=\psi_{P(M)}( w\theta)$ by Definition~\ref{def:D} and Definition-Proposition~\ref{def-pro:ARproduct}. Note that the fifth equality in the above equalities is derived as follows: Since $M$ is a crystal module, $I_w \tens{\Pi}\Hom{\Pi}{I_w}{M}$ is an $\E$-filtered module in $\Fac{I_w}$. We have applied the functor $\Hom{\Pi}{I_w}{-}$ and have used Proposition \ref{prop:I_wrankvec} and the equivalence in the Proof of Proposition \ref{prop:dualtors}.

Now, we apply the functor  $\Hom{\Pi}{-}{M^*}$ for an exact sequence $$0\rightarrow I_w\tens{\Pi}{\Pi e_i}\rightarrow \Pi e_i \rightarrow (\Pi/I_w)\tens{\Pi}{\Pi e_i}\rightarrow 0,$$
so that we have:
$$0\rightarrow \Hom{\Pi}{(\Pi/I_w)\tens{\Pi}{\Pi e_i}}{M^*}\rightarrow \Hom{\Pi}{\Pi e_i}{M^*}$$$$\rightarrow \Hom{\Pi}{I_w\tens{\Pi}{\Pi e_i}}{M^*} \rightarrow \Ext{1}{\Pi}{(\Pi/I_w)\tens{\Pi}{\Pi e_i}}{M^*}\rightarrow 0.$$
By using the generalized Crawley-Boevey formula and the above sequence, we have:
\begin{align*}
    &\dim \Hom{\Pi}{M^*}{(\Pi/I_w)\tens{\Pi}{\Pi e_i}}\\
    =&-\dim \Hom{\Pi}{(\Pi/I_w)\tens{\Pi}{\Pi e_i}}{M^*} + \Ext{1}{\Pi}{(\Pi/I_w)\tens{\Pi}{\Pi e_i}}{M^*} + (\rankvec{(\Pi/I_w)\tens{\Pi}{\Pi e_i}}, \rankvec{M^*})_C\\
    =&\dim \Hom{\Pi}{I_w\tens{\Pi}{\Pi e_i}}{M^*} - \Hom{\Pi}{\Pi e_i}{M^*}+(\rankvec{(\Pi/I_w)\tens{\Pi}{\Pi e_i}}, \rankvec{M^*})_C\\
    =&\dim \Hom{\Pi}{I_w \tens{\Pi}{\Pi e_i}}{M^*}-\dim e_i M^* +(\varpi_i-w\varpi_i, \rankvec{M^*})_C.
\end{align*}
Here, we used an equality $\rankvec{((\Pi/I_w)\tens{\Pi}{\Pi e_i})} =\varpi_i -w \varpi_i$ obtained in the Proof of Theorem \ref{VM_iso} for the third equality. Now, by multiplying $\frac{1}{c_i}\langle \theta, \alpha_i \rangle$ and summing up for $i\in Q_0$, we obtain
\begin{align*}
    &\sum_{i\in Q_0} \frac{1}{c_i} \langle \theta, \alpha_i\rangle\dim \Hom{\Pi}{N(w\varpi_i)}{M}\\
    =&\sum_{i\in Q_0} \frac{1}{c_i}\langle \theta, \alpha_i \rangle \dim \Hom{\Pi}{M^*}{(\Pi/I_w)\tens{\Pi}{\Pi e_i}}\,\,\quad \textit{(Theorem \ref{VM_iso})}\\
    =&\left(\sum_{i\in Q_0} \frac{1}{c_i}\langle \theta, \alpha_i \rangle \dim \Hom{\Pi}{N(w\varpi_i)}{M^*}\right)-\langle \theta, \rankvec{M^*} \rangle+ \langle \theta - w\theta, \rankvec{M^*}\rangle\\
    =& \psi_{P(M^*)}(w\theta)-\langle w\theta, \rankvec{M}\rangle.
\end{align*}
Thus, we obtain the assertion because $P(M^*)$ is the image of $P(M)$ by an involution $\mathbf{r}\mapsto \rankvec{M}-\mathbf{r}$ derived from Remark \ref{rem:invHN}.
\end{proof}
\end{theorem}
By our argument Proposition \ref{Jred_tors} and Theorem \ref{supfun}, we obtain the following corollary:
\begin{cor}\label{face_HNpoly}
Let $J\subseteq Q_0$, let $\theta \in F_J$, and let $w\in W$ be $J$-reduced. Let $T$ be a crystal module and let $X\coloneqq \Hom{\Pi}{I_w}{\Tmax{w\theta}/\Tmin{w\theta}}\in \rep{\Pi/I_{w_J}}$. Then, we have
an equality
$$\left\{x\in P(T)\mid \langle w\theta, x\rangle = D_{w\theta}(T)\right\}=\rankvec{\Tmin{w\theta}}+ wP(X).$$
\end{cor}
Now, we compare our generalized HN polytopes and Kamnitzer's Mirkovi\'c-Vilonen polytopes~\cite{Kam1,Kam2}. We recall the definition of MV polytopes:
\begin{definition}
    Let $V\coloneqq R\tens{\mathbb{Z}}{\R}$ for root lattice $R$ associated with the GCM $C$, and let $V^*$ be its $\R$-dual space. By Hahn-Banach theorem, a non-empty convex subset $P$ of $V$ can be written as $\{v \in V\mid \langle v, \alpha \rangle \leq \psi_{P}(\alpha)\,\,\text{for any $\alpha\in V^*$}\}$ where $\psi_{P}\colon V^*\rightarrow \R$ maps $\alpha$ to the maximal value such that $\alpha$ takes on P. We say that a polytope $P$ is a pseudo-Weyl polytope if $\psi_P$ is linear on each Weyl chamber. We write $A_{\gamma}\coloneqq \psi_{P}(\gamma)$ for $\gamma \in \Gamma$. We say that a pseudo-Weyl polytope $P$ is an MV polytope if the tuple of data $(A_{\gamma})_{\gamma\in \Gamma}$ satisfies the following Berenstein-Zelevinsky data \textsf{(BZ1)-(BZ3)}:
    \begin{description}
        \item[\sf{(BZ1)}] Each $A_{\gamma}$ is an integer and each $A_{-\varpi_i}=0$.
        \item[\sf{(BZ2)}] We have the edge inequalities $A_{-w\varpi_i}+A_{-ws_i\varpi_i}+\sum_{j\in I, j\neq i}c_{ji}A_{-w\varpi_j}\geq 0$.
        \item[\sf{(BZ3)}] All of $A_{\gamma}$ satisfy tropical Pl\"ucker relations.
    \end{description}
    Here, we say that $(A_{\gamma})_{\gamma\in \Gamma}$ satisfies tropical Pl\"ucker relations at $(w, i, j)$ for $w\in W$ and $i, j\in Q_0$ if $c_{ij}=0$,
    \begin{enumerate}
        \item if $c_{ij}=c_{ji}=-1$, then
        $$A_{-ws_i\varpi_i}+A_{-ws_j\varpi_j}=\min(A_{-w\varpi_i}+A_{-ws_is_j\varpi_j}, A_{-ws_js_i\varpi_i }+ A_{-w\varpi_j});$$
        \item if $c_{ij}=-1, c_{ji}=-2$, then
        \begin{align*}        &A_{-ws_j\varpi_{j}}+A_{-ws_i s_j\varpi_j}+A_{-ws_i \varpi_i}\\
        =&\min(2A_{-ws_i s_j\varpi_j}+A_{-w\varpi_i}, 2A_{-w\varpi_j}+A_{-ws_is_j s_i\varpi_i}, A_{-w\varpi_j}+A_{-ws_js_is_j\varpi_j}+A_{-ws_i\varpi_i})
        ,\\
        &A_{-ws_js_i\varpi_i}+2A_{-ws_is_j\varpi_j}+A_{-ws_i\varpi_i}+A_{-ws_i\varpi_i}\\=&\min(2A_{-w\varpi_j}+2A_{-ws_is_js_i\varpi_i}, 2A_{-ws_js_is_j\varpi_j}+2A_{-ws_i\varpi_i}, A_{-ws_is_js_i\varpi_i}+ 2A_{-ws_is_j\varpi_j}+A_{-w\varpi_i});
        \end{align*}
        \item if $c_{ij}=-2, c_{ji}=-1$, then
        \begin{align*}
            &A_{-ws_j s_i\varpi_i}+A_{-ws_i\varpi_i}+A_{-ws_is_j\varpi_j}\\
            =&\min(2A_{-ws_i \varpi_i}+A_{-ws_js_is_j\varpi_j}, 2A_{-ws_is_js_i\varpi_i}+A_{-w\varpi_j}, A_{-ws_i s_j s_i\varpi_i}+2A_{-ws_is_j\varpi_j}+A_{-w\varpi_i}),\\
            &A_{-ws_j\varpi_j}+2A_{-ws_i\varpi_i}+A_{-ws_is_j\varpi_j}\\
            =&\min(2A_{-ws_i s_j s_i\varpi_i}+2A_{-w\varpi_j}, 2A_{-w\varpi_i}+2A_{-ws_is_j\varpi_j}, A_{-w\varpi_j}+2A_{-ws_i\varpi_i}+A_{-ws_js_is_j\varpi_j}).
        \end{align*}
    \end{enumerate}
\end{definition}
Note that any maximal irreducible component $Z \in \Irrmax{\Pi(\mathbf{r})}$ contains an open dense subset on which $D_\gamma$ takes a constant value, where $\gamma \in \Gamma$. Let $D_\gamma (Z)$ denote this constant value.
\begin{theorem}\label{nilp_MV}
For each maximal irreducible component $Z \in \Irrmax{\Pi(\mathbf{r})}$ and a generic module $T$ of this component, $(D_{\gamma}(Z))_{\gamma\in \Gamma}$ describes BZ data. In particular, 
$$P(T)=\{v\in \mathbb{R}^n \mid \langle \gamma, v \rangle \leq D_{\gamma}(T), \gamma \in \Gamma\}$$
is an MV polytope associated with the Langlands dual root datum of $C$ except for type $\mathsf{G}_2$.
\begin{proof}
If we take $D$ as the minimal symmetrizer of $C$, then any generalized preprojective algebra associated with $C$ is isomorphic to $\fgpre{C}{kD}$ for some $k\in \Z_{>0}$. Then, we have the reduction functor $\rep{(\fgpre{C}{kD})}\ni M\mapsto M/\varepsilon M \in \rep{(\fgpre{C}{D})}$, where $\varepsilon\coloneqq \sum_{i\in Q_0}\varepsilon_i^{c_i}$ (\cf Geiss-Leclerc-Schr\"oer~\cite{GLS2}). Through this reduction functor, it is enough to consider minimal symmetrizer cases.
The condition \textsf{(BZ1)} is clear from Theorem \ref{thm_stabmod}, and \textsf{(BZ2)} is a consequence of Proposition \ref{prop:edge}. By Corollary \ref{face_HNpoly}, we have enough to show that each 2-face of a generalized HN polytope describes tropical Pl\"ucker relations. Namely, we enough to consider HN filtrations of generic modules in the category of $\Pi/I_{w_J}$-modules, where we put $J=\{i, j\}$ with $i\neq j$ by Lemma~\ref{lem:submodule} and Proposition~\ref{Jred_tors}. Since the case of $c_{ij} =0$ is clear, we consider cases of $c_{ij}\neq 0$. We check that the multiplicities of layer modules in our HN filtrations satisfy the following relations case by case by using classification theorems of generic modules:
\begin{enumerate}
    \item If $c_{ij}=c_{ji}=-1$, then any generic module has a form $X_1^{\oplus a} \oplus X_2^{\oplus b} \oplus X_3^{\oplus c}$. Let $(a_1, a_2, a_3)$ (\textit{resp.} $(a'_1, a'_2, a'_3)$) be the multiplicities of layer modules in the HN filtration with respect to a reduced expression $\mathbf{i}=(i, j, i)$ (\textit{resp.} $\mathbf{i}'=(j, i, j)$). Then, we can prove they satisfy the following list \ref{TP1a}--\ref{TP1b} of equalities by calculation, where $P_i$ (\textit{resp.} $E_i$) is the indecomposable projective (\textit{resp.} generalized simple) module with respect to a vertex $i$. In particular, they satisfy relations \eqref{eq:TP1}:
    \begin{align}
    a'_1=a_2+a_3-p, a'_2=p, a'_3=a_1+a_2-p,
    \label{eq:TP1}
    \end{align}
    where $p=\min(a_1, a_3)$.
    \begin{enumerate}
        \item \label{TP1a} $(X_1, X_2, X_3)=(P_1, P_2, E_1)$
       \begin{align*} 
        (a_1, a_2, a_3)&=(a+c, b, a)\\
        (a'_1, a'_2, a'_3)&=(b, a, b+c)\\
        p&=a
        \end{align*}
        \item \label{TP1b} $(X_1, X_2, X_3)=(P_1, P_2, E_2)$
        \begin{align*}
            (a_1, a_2, a_3)&=(a, b, a+c)\\
        (a'_1, a'_2, a'_3)&=(b+c, a, b)\\
        p&=a
        \end{align*}
    \end{enumerate}
    \item If $c_{ij}=-2, c_{ji}=-1$, then any generic module has a form $X_1^{\oplus a}\oplus X_2^{\oplus b}\oplus X_3^{\oplus c}\oplus X_4^{\oplus d}$ where $X_i$ is one of the 8 Loewy series of modules below ($i=2, j=1$) and $(X_1, X_2, X_3, X_4)$ is one of the below list (\cf \cite[\S 8.2.2]{GLS4}).
    \begin{alignat*}{4}
        P_1&=
        \begin{tikzpicture}[auto,baseline=-3pt]
\node (x) at (0, 0) {};
\node (a) at (0, -0.45) {$1$}; 
\node (b) at (0.3, -0.15) {$2$};  
\node (c) at (-0.3, -0.15) {$1$};
\node (d) at (0.3, 0.15) {$1$};
\node (e) at (-0.3, 0.15) {$2$};
\node (f) at (0, 0.45) {$1$};
\end{tikzpicture};
&P_2&=
\begin{tikzpicture}[auto,baseline=-3pt]
\node (x) at (0, 0) {};
\node (a) at (0, -0.45) {$2$}; 
\node (b) at (0, -0.15) {$1$};  
\node (d) at (0, 0.15) {$1$};
\node (f) at (0, 0.45) {$2$};
\end{tikzpicture};
&E_1&=\begin{tikzpicture}[auto,baseline=-3pt]
\node (x) at (0, 0) {};
\node (b) at (0, -0.15) {$1$};  
\node (d) at (0, 0.15) {$1$};
\end{tikzpicture};
&E_2&=\begin{tikzpicture}[auto,baseline=-3pt]
\node (x) at (0, 0) {$2$};
\end{tikzpicture};\\
T_1&=
\begin{tikzpicture}[auto,baseline=-3pt]
\node (b) at (0, -0.3) {$2$};  
\node (x) at (0, 0) {$1$};
\node (d) at (0, 0.3) {$1$};
\end{tikzpicture};
&T_2&=
\begin{tikzpicture}[auto,baseline=-3pt]
\node (b) at (0, -0.3) {$1$};  
\node (x) at (0, 0) {$1$};
\node (d) at (0, 0.3) {$2$};
\end{tikzpicture};
&T_3&=
        \begin{tikzpicture}[auto,baseline=-3pt]
\node (x) at (0, 0) {};
\node (b) at (0.3, -0.15) {$2$};  
\node (d) at (0.3, 0.15) {$1$};
\node (e) at (-0.3, 0.15) {$2$};
\node (f) at (0, 0.45) {$1$};
\end{tikzpicture};
&T_4&=
        \begin{tikzpicture}[auto,baseline=-3pt]
\node (x) at (0, 0) {};
\node (a) at (0, -0.45) {$1$}; 
\node (b) at (0.3, -0.15) {$2$};  
\node (c) at (-0.3, -0.15) {$1$};
\node (e) at (-0.3, 0.15) {$2$};
\end{tikzpicture}.
    \end{alignat*}
    
    Now, let $(a_1, a_2, a_3, a_4)$ (\textit{resp.} $(a'_1, a'_2, a'_3, a'_4)$) be the multiplicities of layer modules in the HN filtration from a reduced expression $\mathbf{i}=(j, i, j, i)$ (\textit{resp.} $\mathbf{i}'=(i, j, i, j)$). Then, they satisfy the following list \ref{TP2a}--\ref{TP2f} of equalities by calculating case by case (\cf Example \ref{ex:layerB2}). In particular, they satisfy relations \eqref{eq:TP2}:
    \begin{align}
    a'_1=a_2+a_3+a_4-p_1; a'_2=2p_1-p_2;
    a'_3=p_2-p_1; a'_4=a_1+2a_2+a_3-p_2,\label{eq:TP2}
    \end{align}
    where 
    $$p_1=\min(a_1+a_2, a_1+a_4, a_3+a_4),$$
    $$p_2=\min(a_1+2a_2, a_1+2a_4, a_3+2a_4).$$
    
    \begin{enumerate}
        \item \label{TP2a} $(X_1, X_2, X_3, X_4)=(P_1, P_2, E_1, T_1)$
        \begin{align*}
            (a_1, a_2, a_3, a_4)&=(b, a, b+d, a+c)\\
            (a'_1, a'_2, a'_3, a'_4)&=(a+c+d, b, a, b+d)\\
            (p_1, p_2)&=(a+b, 2a+b)
        \end{align*}
        \item $(X_1, X_2, X_3, X_4)=(P_1, P_2, E_1, T_2)$
        \begin{align*}
            (a_1, a_2, a_3, a_4)&=(b+d, a, b, a+c+d)\\
            (a'_1, a'_2, a'_3, a'_4)&=(a+c, b+d, a, b)\\
            (p_1, p_2)&=(a+b+d, 2a+b+d)
            \end{align*}
        \item $(X_1, X_2, X_3, X_4)=(P_1, P_2, E_2, T_3)$
        \begin{align*}
            (a_1, a_2, a_3, a_4)&=(b+c, a+d, b, a)\\
            (a'_1, a'_2, a'_3, a'_4)&=(a+d, b, a, b+c+2d)\\
            (p_1, p_2)&=(a+b, 2a+b)
            \end{align*}
        \item $(X_1, X_2, X_3, X_4)=(P_1, P_2, E_2, T_4)$
        \begin{align*}
            (a_1, a_2, a_3, a_4)&=(b+c+2d, a, b, a+d)\\
            (a'_1, a'_2, a'_3, a'_4)&=(a, b, a+d, b+c)\\
            (p_1, p_2)&=(a+b+d, 2a+b+2d)
            \end{align*}
        \item $(X_1, X_2, X_3, X_4)=(P_1, P_2, T_1, T_3)$
        \begin{align*}
            (a_1, a_2, a_3, a_4)&=(b, a+d, b+c, a)\\
            (a'_1, a'_2, a'_3, a'_4)&=(a+c+d, b, a, b+c+2d)\\
            (p_1, p_2)&=(a+b, 2a+b)
            \end{align*}
        \item \label{TP2f} $(X_1, X_2, X_3, X_4)=(P_1, P_2, T_2, T_4)$
        \begin{align*}
            (a_1, a_2, a_3, a_4)&=(b+c+2d, a, b, a+c+d)\\
            (a'_1, a'_2, a'_3, a'_4)&=(a, b+c, a+d, b)\\
            (p_1, p_2)&=(a+b+c+d, 2a+b+c+2d)
            \end{align*}
    \end{enumerate}

    \item If $c_{ij}=-1, c_{ji}=-2$, then any generic module has a form $X_1^{\oplus a}\oplus X_2^{\oplus b}\oplus X_3^{\oplus c}\oplus X_4^{\oplus d}$ where $X_i$ is one of the 8 Loewy series of modules below ($i=1, j=2$) and $(X_1, X_2, X_3, X_4)$ is one of the below list (\cf \cite[\S 8.2.2]{GLS4}).
    \begin{alignat*}{4}
        P_1&=
        \begin{tikzpicture}[auto,baseline=-3pt]
\node (x) at (0, 0) {};
\node (a) at (0, -0.45) {$1$}; 
\node (b) at (0.3, -0.15) {$2$};  
\node (c) at (-0.3, -0.15) {$1$};
\node (d) at (0.3, 0.15) {$1$};
\node (e) at (-0.3, 0.15) {$2$};
\node (f) at (0, 0.45) {$1$};
\end{tikzpicture};
&P_2&=
\begin{tikzpicture}[auto,baseline=-3pt]
\node (x) at (0, 0) {};
\node (a) at (0, -0.45) {$2$}; 
\node (b) at (0, -0.15) {$1$};  
\node (d) at (0, 0.15) {$1$};
\node (f) at (0, 0.45) {$2$};
\end{tikzpicture};
&E_1&=\begin{tikzpicture}[auto,baseline=-3pt]
\node (x) at (0, 0) {};
\node (b) at (0, -0.15) {$1$};  
\node (d) at (0, 0.15) {$1$};
\end{tikzpicture};
&E_2&=\begin{tikzpicture}[auto,baseline=-3pt]
\node (x) at (0, 0) {$2$};
\end{tikzpicture};\\
T_1&=
\begin{tikzpicture}[auto,baseline=-3pt]
\node (b) at (0, -0.3) {$2$};  
\node (x) at (0, 0) {$1$};
\node (d) at (0, 0.3) {$1$};
\end{tikzpicture};
&T_2&=
\begin{tikzpicture}[auto,baseline=-3pt]
\node (b) at (0, -0.3) {$1$};  
\node (x) at (0, 0) {$1$};
\node (d) at (0, 0.3) {$2$};
\end{tikzpicture};
&T_3&=
        \begin{tikzpicture}[auto,baseline=-3pt]
\node (x) at (0, 0) {};
\node (b) at (0.3, -0.15) {$2$};  
\node (d) at (0.3, 0.15) {$1$};
\node (e) at (-0.3, 0.15) {$2$};
\node (f) at (0, 0.45) {$1$};
\end{tikzpicture};
&T_4&=
        \begin{tikzpicture}[auto,baseline=-3pt]
\node (x) at (0, 0) {};
\node (a) at (0, -0.45) {$1$}; 
\node (b) at (0.3, -0.15) {$2$};  
\node (c) at (-0.3, -0.15) {$1$};
\node (e) at (-0.3, 0.15) {$2$};
\end{tikzpicture}.
    \end{alignat*}
    
    Now, let $(a_1, a_2, a_3, a_4)$ (\textit{resp.} $(a'_1, a'_2, a'_3, a'_4)$) be the multiplicities of layer modules in the HN filtration from a reduced expression $\mathbf{i}=(j, i, j, i)$ (\textit{resp.} $\mathbf{i}'=(i, j, i, j)$). Then, they satisfy the following list \ref{TP3a}--\ref{TP3f} of equalities by calculating case by case (\cf Example \ref{ex:layerB2}). In particular, they satisfy relations \eqref{eq:TP3}:
    \begin{align}
    a'_1=a_2+2a_3+a_4-p_2, a'_2=p_2-p_1,
    a'_3=2p_1-p_2, a'_4=a_1+a_2+a_3-p_1, \label{eq:TP3}
    \end{align}
    where
    $$p_1=\min(a_1+a_2, a_1+a_4,a_3+a_4),$$
    $$p_2=\min(2a_1+a_2, 2a_1+a_4, 2a_3+a_4).$$
    \begin{enumerate}
        \item \label{TP3a} $(X_1, X_2, X_3, X_4)=(P_1, P_2, E_1, T_1)$
        \begin{align*}
            (a_1, a_2, a_3, a_4)&=(b, a, b+d, a+c)\\
            (a'_1, a'_2, a'_3, a'_4)&=(a+c+d, b, a, b+d)\\
            (p_1, p_2)&=(a+b+d, 2a+b+d)
        \end{align*}
        \item $(X_1, X_2, X_3, X_4)=(P_1, P_2, E_1, T_2)$
        \begin{align*}
            (a_1, a_2, a_3, a_4)&=(b+d, a, b, a+c+d)\\
            (a'_1, a'_2, a'_3, a'_4)&=(a+c, b+d, a, b)\\
            (p_1, p_2)&=(a+b, 2a+b)
            \end{align*}
        \item $(X_1, X_2, X_3, X_4)=(P_1, P_2, E_2, T_3)$
        \begin{align*}
            (a_1, a_2, a_3, a_4)&=(b+c, a+d, b, a)\\
            (a'_1, a'_2, a'_3, a'_4)&=(a+d, b, a, b+c+2d)\\
            (p_1, p_2)&=(a+b+d, 2a+b+2d)
            \end{align*}
        \item $(X_1, X_2, X_3, X_4)=(P_1, P_2, E_2, T_4)$
        \begin{align*}
            (a_1, a_2, a_3, a_4)&=(a, b, a+d, b+c)\\
            (a'_1, a'_2, a'_3, a'_4)&=(b+c+2d, a, b, a+d)\\
            (p_1, p_2)&=(a+b, 2a+b)
            \end{align*}
        \item $(X_1, X_2, X_3, X_4)=(P_1, P_2, T_1, T_3)$
        \begin{align*}
            (a_1, a_2, a_3, a_4)&=(a+c+d, b, a, b+c+2d)\\
            (a'_1, a'_2, a'_3, a'_4)&=(b, a+d, b+c, a)\\
            (p_1, p_2)&=(a+b+c+d, 2a+b+c+2d)
            \end{align*}
        \item \label{TP3f} $(X_1, X_2, X_3, X_4)=(P_1, P_2, T_2, T_4)$
        \begin{align*}
            (a_1, a_2, a_3, a_4)&=(a, b+c, a+d, b)\\
            (a'_1, a'_2, a'_3, a'_4)&=(b+c+2d, a, b, a+c+d)\\
            (p_1, p_2)&=(a+b, 2a+b)
            \end{align*}
    \end{enumerate}
\end{enumerate}
In fact, these are equivalent to the tropical Pl\"ucker relations~\cite[Proposition 5.2]{Kam1}.
\end{proof}
\end{theorem}
\begin{rem}
 Since we have Theorem~\ref{Saito-refl} and Theorem~\ref{thm:HNproperty} also for type $\mathsf{G_2}$, two Harder-Narasimhan filtrations of generic modules of maximal irreducible components with respect to two reduced expressions of the longest element have Lusztig data as multiplicities of layer modules in general. In particular, generalized Harder-Narasimhan polytopes of generic modules of maximal irreducible components satisfies BZ data also for type $\mathsf{G}_2$ (\cf McNamara~\cite[\S 7]{Mc}).
 
Since $\mathsf{A}_2$ and $\mathsf{B}_2$ type generalized preprojective algebras with minimal symmetrizers are special biserial algebras, their module categories are understood by Auslander-Reiten theory~\cite{BR}. However, for type $\mathsf{G}_2$ case, this approach for explicit calculations for $2$-faces based is invalid. Especially, there exist infinitely many maximal rigid objects \footnote{The author was informed this result by C. Gei{\ss} and B. Leclerc while his stayed at Caen. He thanks for their kindness.}. But, also in this case, we expect that maximal rigid objects determined by some combinatorics of two-sided ideals $I_w$ describe BZ data for type $\mathsf{G}_2$. A general theory about this direction will appear in the author's future work.
\end{rem}
Finally, we relate the crystal structure on nilpotent varieties in \cite{GLS4} and that of the set $\MV$ of MV polytopes as a generalization of \cite[Theorem 6.3]{BK}. We recall the crystal structure on $\MV$, which has been developed by Kamnitzer~\cite[\S 3.6]{Kam2}. 
\begin{def-pro}
Let $P\in\MV$ with vertex data $(\mu_{w})_{w\in W}$ and BZ data $(A_\gamma)_{\gamma\in\Gamma}$ and let 
$$c\coloneqq A_{w\varpi_i}+A_{ws_i\varpi_i}+\sum_{j\in I, j\neq i}c_{ji}A_{w\varpi_j}.$$
Then, we define
\begin{gather*}
    \wt{P}\coloneqq \mu_{w_0}; \quad
    \varphi_{i}(P)=c; \quad
    \varepsilon_{i}(P)=c-\langle \alpha_i, \wt{P}\rangle;\\
    \wt{\tilde{e}_i P}=\wt{P}+\alpha _i;\quad
    \wt{\tilde{f}_i P}=\wt P-\alpha_i\,\,(\varphi_i(P)>0).
\end{gather*}
Then, $(\MV, \wt{(-)}, \tilde{e}_i, \tilde{f}_i, \varepsilon_i, \varphi_i)$ is isomorphic to $B(-\infty)$ as crystals.
\end{def-pro}
\begin{theorem}\label{bij_MV&B}
The map $P{(-)}\colon \mathcal{B}\rightarrow \MV$ is a crystal bijection.
\begin{proof}
Let $Z\in \mathcal{B}$ and let $T$ be a generic point of $Z$. Then, we have
\begin{align*}
    \rk_i \fac_i{T}&= \frac{1}{c_i}\dim\Hom{\Pi}{T}{E_i}\\
    &= \frac{1}{c_i}\dim\Hom{\Pi}{E_i}{T^*}\\
    &= D_{s_i \varpi_i}(T)-\langle s_i\varpi_i, \rankvec{T}\rangle\\
    &=D_{s_i\varpi_i}(T)+D_{\varpi_i}(T)+\sum_{j\in J,j\neq i}c_{ji} D_{\varpi_j}(T)
\end{align*}
Hence, $T\mapsto \rk_i{\fac_i{T}}$ takes a constant value $$D_{s_i\varpi_i}(P(Z))+D_{\varpi_i}(P(Z))+\sum_{j\in J,j\neq i}c_{ji} D_{\varpi_j}(P(Z))=\varphi_i(P(T))$$ on an open dense subset in $Z$.
Let $i \in Q_0$ and let $Z'\coloneqq \tilde{f}^{\max}_i Z$. Now, we have
$$\wt{P(Z')}=\wt{P(Z)}-\varphi_i(P(Z))\alpha_i=\wt{\tilde{f}^{\max}_i P(Z)}.$$
Now, $D_\gamma(T)=D_{\gamma}(Z)$ and $D_{\gamma}(T_{s_i})=D_{\gamma}(Z')$ for any chamber weight $\gamma$. Since $D_{\gamma}(T)=D_{\gamma}(T_{s_i})$ for $\left\langle \gamma, \alpha_i \right\rangle \leq 0$, we obtain $D_{\gamma}(Z)=D_{\gamma}(Z')$. Thus, $P(-)$ preserves $\tilde{f}_i^{\max}$. Similarly, $P(-)$ preserves $\tilde{e}_i^{\max}$. The assertions about $\varphi_i$ and $\varepsilon_i$ are directly proved from definitions.
\end{proof}
\end{theorem}
\subsection*{Acknowledgement}
The author thanks Ryo Fujita, Christof Gei\ss, Yuya Ikeda, Syu Kato, Yoshiyuki Kimura, Bernard Leclerc and Yuya Mizuno for their valuable discussions and comments on this work. Especially, the author expresses gratitude to his supervisor Syu Kato for his guidance for geometric representation theory about affine Grassmannians and quiver varieties. The author was indebted to Bernard Leclerc because he carefully read a draft of this paper and was kind to give many constructive suggestions. A part of this work about Harder-Narasimhan polytopes is motivated from a conversation with Yuya Mizuno and Yoshiyuki Kimura. He also thanks all of his colleagues for their special effort in difficult situations of the COVID-19 pandemic. Finally, the author was indebted to anonymous referees for their constructible suggestions.
\bibliographystyle{amsalpha}
\bibliography{crystal.bib}
\end{document}